\documentclass[a4paper,12pt,twoside]{article}
\usepackage{graphicx}\usepackage{a4wide}
\usepackage{amssymb,amsmath,amsthm}
\usepackage{siunitx}
\usepackage{subfigure,color,colordvi,graphicx,xcolor}
\usepackage[only,llbracket,rrbracket,interleave]{stmaryrd}
\usepackage[superscript,sort,compress]{cite}
\usepackage{enumitem}

\newcommand{\bm}[1]{\text{\boldmath $#1$\unboldmath}}

\newcommand{\abs}[1]{\lvert#1\rvert}
\newcommand{\norm}[1]{\lVert#1\rVert}

\newcommand{\vect}[1]{\mathbf{#1}}
\newcommand{\mat}[1]{\mathbf{#1}}

\newcommand{\eltwo}{\ensuremath{\mathcal{L}_2}}
\newcommand{\elinf}{\ensuremath{\mathcal{L}_\infty}}

\newcommand{\nsd}    {\texttt{n}_{\texttt{sd}}}
\newcommand{\npar}{\texttt{n}_{\texttt{pa}}}

\newcommand{\numel}{\texttt{n}_{\texttt{el}}}

\DeclareMathOperator{\adj}{adj}

\newcommand{\jump}[1]{\llbracket #1\rrbracket}
\DeclareMathOperator*{\argmin}{arg\,min}

\newcommand{\Div}{{\bm{\nabla}\!\!_{\bmu} \cdot\,}}
\newcommand{\Grad}{\bm{\nabla}\!\!_{\bmu}}

\newcommand{\bu}{\bm{u}}
\newcommand{\bs}{\bm{s}}
\newcommand{\bt}{\bm{g}_N}
\newcommand{\bn}{\bm{n}}
\newcommand{\bD}{\bm{D}}
\newcommand{\bE}{\bm{E}}
\newcommand{\Insd}{\mat{I}_{\!\nsd\!}}
\newcommand{\InsdTwo}{\mat{I}_{\!\nsd \times \nsd\!}}
\newcommand{\bF}{\bm{F}}

\newcommand{\bL}{\bm{L}}
\newcommand{\hu}{\hat{u}}
\newcommand{\bhu}{\bm{\hu}}
\newcommand{\btau}{\bm{\tau}}

\newcommand{\sVh}{\mathcal{V}^h}
\newcommand{\shVh}{\mathcal{\widehat{V}}^h}
\newcommand{\sLh}{\mathcal{L}^h}
\newcommand{\Lh}{\bm{\mathcal{L}}^h}

\newcommand{\sVhMu}{\mathcal{V}^h_{\!\! \bmu}}
\newcommand{\VhMu}{\bm{\mathcal{V}}^h_{\!\! \bmu}}
\newcommand{\WhMu}{\bm{\mathcal{W}}^h_{\!\! \bmu}}
\newcommand{\hVhMu}{\bm{\mathcal{\widehat{V}}}^h_{\!\! \bmu}}

\newcommand{\Pk}{\mathcal{P}^{\nDeg}}

\newcommand{\de}{\delta\!}

\newcommand{\Vh}{\bm{\mathcal{V}}^h}
\newcommand{\Wh}{\bm{\mathcal{W}}^h}
\newcommand{\hVh}{\bm{\mathcal{\widehat{V}}}^h}

\newcommand{\bW}{\bm{W}}
\newcommand{\bv}{\bm{v}}
\newcommand{\hv}{\hat{v}}

\newcommand{\bhv}{\bm{\hat{v}}}

\newcommand{\bx}{\bm{x}^\bmu}
\newcommand{\bX}{\bm{x}}

\newcommand{\bmu}{\bm{\mu}}

\newcommand{\intEMu}[2]{\big(#1,#2\big)_{\Omega_e^{\bmu} \times \bI }}
\newcommand{\intBEMu}[2]{\langle #1,#2   \rangle_{\partial \Omega_e^{\bmu} \times \bI }}

\newcommand{\intBNoDMu}[2]{\langle #1,#2 \rangle_{( \partial \Omega_e^{\bmu} \setminus \Gamma_D^{\bmu} ) \times \bI }}
\newcommand{\intBNoDSMu}[2]{\langle #1,#2 \rangle_{( \partial \Omega_e^{\bmu} \setminus (\Gamma_D^{\bmu} \cup \Gamma_S^{\bmu} )) \times \bI }}
\newcommand{\intBDMu}[2]{\langle #1,#2   \rangle_{( \partial \Omega_e^{\bmu} \cap \Gamma_D^{\bmu} ) \times \bI }}
\newcommand{\intBNMu}[2]{\langle #1,#2   \rangle_{( \partial \Omega_e^{\bmu} \cap \Gamma_N^{\bmu} ) \times \bI }}
\newcommand{\intBSMu}[2]{\langle #1,#2   \rangle_{( \partial \Omega_e^{\bmu} \cap \Gamma_S^{\bmu} ) \times \bI }}

\newcommand{\intExI}[2]{\big(#1,#2\big)_{\Omega_e \times \bI }}
\newcommand{\intBExI}[2]{\langle #1,#2   \rangle_{\partial \Omega_e \times \bI }}

\newcommand{\intBNoDxI}[2]{\langle #1,#2 \rangle_{( \partial \Omega_e \setminus \Gamma_D ) \times \bI }}
\newcommand{\intBNoDSxI}[2]{\langle #1,#2 \rangle_{( \partial \Omega_e \setminus (\Gamma_D \cup \Gamma_S )) \times \bI }}
\newcommand{\intBDxI}[2]{\langle #1,#2   \rangle_{( \partial \Omega_e \cap \Gamma_D ) \times \bI }}
\newcommand{\intBNxI}[2]{\langle #1,#2   \rangle_{( \partial \Omega_e \cap \Gamma_N ) \times \bI }}
\newcommand{\intBSxI}[2]{\langle #1,#2   \rangle_{( \partial \Omega_e \cap \Gamma_S ) \times \bI }}

\newcommand{\intE}[2]{\big(#1,#2\big)_{\Omega_e }}
\newcommand{\intI}[2]{\big(#1,#2\big)_{\bI}}
\newcommand{\intBE}[2]{\langle #1,#2   \rangle_{\partial \Omega_e}}

\newcommand{\intBNoD}[2]{\langle #1,#2 \rangle_{\partial \Omega_e \setminus \Gamma_D}}
\newcommand{\intBNoDS}[2]{\langle #1,#2 \rangle_{\partial \Omega_e \setminus (\Gamma_D \cup \Gamma_S) } }
\newcommand{\intBD}[2]{\langle #1,#2   \rangle_{\partial \Omega_e \cap \Gamma_D}}
\newcommand{\intBN}[2]{\langle #1,#2   \rangle_{\partial \Omega_e \cap \Gamma_N}}
\newcommand{\intBS}[2]{\langle #1,#2   \rangle_{\partial \Omega_e \cap \Gamma_S }}

\newcommand{\Lpgd}{\bL_{_{\texttt{PGD}}}}
\newcommand{\bupgd}{\bu_{_{\texttt{PGD}}}}
\newcommand{\ppgd}{p_{_{\texttt{PGD}}}}
\newcommand{\bhupgd}{\bhu_{_{\texttt{PGD}}}}
\newcommand{\rpgd}{\rho_{_{\texttt{PGD}}}}
\newcommand{\dLpgd}{\De \bL_{_{\texttt{PGD}}}}
\newcommand{\dbupgd}{\De \bu_{_{\texttt{PGD}}}}
\newcommand{\dppgd}{\De p_{_{\texttt{PGD}}}}
\newcommand{\dbhupgd}{\De \bhu_{_{\texttt{PGD}}}}
\newcommand{\drpgd}{\De \rho_{_{\texttt{PGD}}}}

\newcommand{\fL}{\bm{F}_{\!\! L}}
\newcommand{\fU}{\bm{f}_{\!\! u}}
\newcommand{\fP}{f_{\! p}}
\newcommand{\fHU}{\bm{f}_{\!\! \hu}}
\newcommand{\fR}{f_{\! \rho}}

\newcommand{\ampL}{\sigma_{\! L}}
\newcommand{\ampU}{\sigma_{\! u}}
\newcommand{\ampP}{\sigma_{\! p}}
\newcommand{\ampHU}{\sigma_{\! \hu}}
\newcommand{\ampR}{\sigma_{\! \rho}}

\newcommand{\DivX}{{\bm{\nabla}\cdot\,}}
\newcommand{\GradX}{\bm{\nabla}}
\newcommand{\nmap}{\texttt{n}_{\texttt{M}}}

\newcommand{\bJ}{\mat{J}}
\newcommand{\Jk}{\bJ^k}

\newcommand{\ndet}{\texttt{n}_{\texttt{d}}}
\newcommand{\nadj}{\texttt{n}_{\texttt{a}}}

\newcommand{\bA}{\mat{A}}
\newcommand{\mapping}{\bm{\mathcal{M}}_{\bmu}}
\newcommand{\Jaco}{\mat{J}_{\!\bmu}}

\newcommand{\nDir}{\texttt{n}_{\texttt{D}}}
\newcommand{\nNeu}{\texttt{n}_{\texttt{N}}}
\newcommand{\nSou}{\texttt{n}_{\texttt{S}}}

\newcommand{\bg}{\bm{g}}
\newcommand{\bgD}{\bg_{D}}
\newcommand{\bgN}{\bg_{N}}
\newcommand{\bgS}{\bg_{S}}

\newcommand{\bI}{\bm{\mathcal{I}}}
\newcommand{\I}{\mathcal{I}}

\newcommand{\bM}{\mat{M}}

\newcommand{\nDeg}{\ensuremath{k}}

\newcommand{\Rin} {R_{\text{in}}}
\newcommand{\Rout}{R_{\text{out}}}
\newcommand{\Rref}{R_{\text{ref}}}
\newcommand{\Rint}{R_{\text{int}}}
\newcommand{\Vin} {\Omega_{\text{in}}}
\newcommand{\Vout}{\Omega_{\text{out}}}

\newcommand{\De}{\varDelta}

\renewcommand{\emph}[1]{\textit{#1}}
\newenvironment{keywords}{\begin{quote}\emph{\textbf{Keywords:}}}{\end{quote}}

\theoremstyle{definition}

\newtheorem{remark}{Remark}

\RequirePackage{algorithm}[0.1] 
\usepackage{algorithmic} 

\usepackage{scalerel}
\usepackage{stackengine}
\stackMath
\def\hatgap{0pt}
\def\subdown{-2pt}
\newcommand\reallywidehat[2][]{
	\renewcommand\stackalignment{l}
	\stackon[\hatgap]{#2}{
		\stretchto{\scalerel*[\widthof{$#2$}]{\kern-.6pt\bigwedge\kern-.6pt}{\rule[-\textheight/2]{1ex}{\textheight}}}{0.5ex}_{\smash{\belowbaseline[\subdown]{\scriptstyle#1}}}}}
		
\makeatletter
\newcommand*\wt[2][0.2ex]{%
        \begingroup
        \mathchoice{\wt@helper{#1}{#2}{\displaystyle}{\textfont}}
                   {\wt@helper{#1}{#2}{\textstyle}{\textfont}}
                   {\wt@helper{#1}{#2}{\scriptstyle}{\scriptfont}}
                   {\wt@helper{#1}{#2}{\scriptscriptstyle}{\scriptscriptfont}}%
        \endgroup
        #2%
}
\newcommand*\wt@helper[4]{%
        \def\currentfont{\the#41}%
        \def\currentskewchar{\char\the\skewchar\currentfont}%
        \setbox\tw@\hbox{\currentfont#2\currentskewchar}%
        \dimen@ii\wd\tw@
        \setbox\tw@\hbox{\currentfont#2{}\currentskewchar}%
        \advance\dimen@ii-\wd\tw@
        \rlap{\raisebox{-#1}{$\m@th#3\kern\dimen@ii\widetilde{\phantom{#2}}$}}%
}
\makeatother

\graphicspath{{figsPDF/}{figsPNG/}{figs/}}

\begin{document}
\title{Hybridisable discontinuous Galerkin solution of geometrically parametrised Stokes flows}

\author{
\renewcommand{\thefootnote}{\arabic{footnote}}
			  R. Sevilla\footnotemark[1]\textsuperscript{ \ ,}*, 
			  L. Borchini\footnotemark[1]\textsuperscript{ \ ,}\footnotemark[2] \ ,
			  M. Giacomini\footnotemark[2]\textsuperscript{ \ ,}\footnotemark[3]  \ and
             A. Huerta\footnotemark[2]\textsuperscript{ \ ,}\footnotemark[3]
}

\date{\today}
\maketitle

\renewcommand{\thefootnote}{\arabic{footnote}}

\footnotetext[1]{Zienkiewicz Centre for Computational Engineering, College of Engineering, Swansea University, Wales, UK}
\footnotetext[2]{Laboratori de C\`alcul Num\`eric (LaC\`aN), ETS de Ingenieros de Caminos, Canales y Puertos, Universitat Polit\`ecnica de Catalunya, Barcelona, Spain}
\footnotetext[3]{Centre Internacional de M\`etodes Num\`erics en Enginyeria (CIMNE), Barcelona, Spain.
\vspace{5pt}\\
* Corresponding author: Ruben Sevilla. \textit{E-mail:} \texttt{r.sevilla@swansea.ac.uk}
}

\begin{abstract}
This paper proposes a novel computational framework for the solution of geometrically parametrised flow problems governed by the Stokes equation. The proposed method uses a high-order hybridisable discontinuous Galerkin formulation and the proper generalised decomposition rationale to construct an off-line solution for a given set of geometric parameters. The generalised solution contains the information for all the geometric parameters in a user-defined range and it can be used to compute sensitivities. The proposed approach circumvents many of the weaknesses of other approaches based on the proper generalised decomposition for computing generalised solutions of geometrically parametrised problems. Four numerical examples show the optimal approximation properties of the proposed method and demonstrate its applicability in two and three dimensions.
\end{abstract}


\section{Introduction}

Reduced order models (ROMs) have become commonplace in many areas of computational sciences and engineering~\cite{quarteroni2014reduced}. Some popular ROMs used to reduce the complexity of high dimensional problems include the reduced basis method~\cite{rozza2008reduced}, the proper orthogonal decomposition (POD)~\cite{berkooz1993proper,lieu2006reduced,ballarin2019pod} and the proper generalised decomposition (PGD)~\cite{Chinesta-Keunings-Leygue,chinesta2011short,PGD-CCH:14}. 

One of the main attractive properties of the PGD is its ability to build reduced basis without prior knowledge of the solution~\cite{Chinesta-Keunings-Leygue,chinesta2011short,PGD-CCH:14}. However, the intrusive implementation and the difficulty in handling geometrically parametrised problems has often been considered a difficulty when considering its application to complex problems. In recent years, there have been an increase in non-intrusive implementations of the PGD~\cite{giraldi2015or,zou2018nonintrusive,tsiolakis2020nonintrusive}. In terms of geometrically parametrised problems, early work focused on solutions tailored to specific problems~\cite{FC-CLBACGAAH:13,leygue2010first,bognet2012advanced,heuze2015parametric} or strategies only applicable in a context of low order approximations~\cite{AH-AHCCL:14,SZ-ZDMH:15}. More recently, a general approach to deal with geometrically parametrised problems in a CAD environment was proposed~\cite{sevilla2020solution}. The PGD strategy presented in~\cite{sevilla2020solution} used a classical finite element (FE) discretisation of Stokes flow problems, leading to the need to use the so-called high-order PGD projection~\cite{DM-MZH:15} to separate some terms of the weak formulation. 

In this work a PGD strategy is proposed in the framework of the hybridisable discontinuous Galerkin (HDG) method~\cite{Jay-CGL:09,Cockburn-Encyclopedia,RS-SH:16,MG-GS:19,MG-GSH-20}. The use of a mixed formulation is shown to be beneficial as all the terms of the weak formulation can be written in a separated form, as required by the PGD, without invoking to the memory intensive high-order PGD projection. The use of the HDG method for the spatial discretisation also guarantees that equal order of approximation can be used for all the variables circumventing the so-called Ladyzhenskaya-Babu{\v s}ka-Brezzi  (LBB) condition. This is of special importance in this work, where geometrically parametrised domains are considered with curved boundaries. The use of the same degree of approximation for all the variables means that standard isoparametric elements can be used. In contrast, the work in~\cite{sevilla2020solution}, employing standard FEs, required the use of sub-parametric or super-parametric formulations in the presence of curved boundaries due to the different degree of approximation used for the velocity and pressure, as required to satisfy the LBB condition. Furthermore, the proposed HDG-PGD approach facilitates the imposition of the Dirichlet boundary conditions as in the HDG context all boundary conditions are weakly imposed.

The formulation is presented using Stokes flows as the model problem. However, it is worth mentioning that there has been a substantial effort in developing HDG methods for a variety of problems in different areas of science and engineering~\cite{JP-PNC:10-compressible,JP-NPC:10-incompressible,Nguyen-NPC:11-Maxwell,Nguyen-NPC:11-acoustics,GG-GFH:13,Cockburn-SCS:09-elasticity,Cockburn-KLC:15-NonLinearElasticity,RS-SGKH:18} and therefore, the proposed approach can be easily extended to a wide range of problems. It is also worth noting that the integration within a CAD environment proposed in~\cite{sevilla2020solution} is also feasible given the recent development of a coupled HDG-NEFEM formulation for fluid~\cite{HDG-NEFEM} and solid mechanics~\cite{HDG-NEFEM-Elas}.
	
The structure of the remainder of the paper is as follows. Section \ref{sc:problemStatement} presents the Stokes flow problem on a geometrically parametrised domain and the corresponding multi-dimensional parametric problem. The HDG formulation for the multi-dimensional parametric Stokes problem is described in section~\ref{sc:hdg}. The proposed PGD rationale is described in detail in section~\ref{sc:pgd}. Section~\ref{sc:examples} presents a series of numerical examples involving Stokes flow problems in two and three dimensions. Finally, section~\ref{sc:conclusions} presents the conclusions of the work that has been presented. 
	
\section{Problem statement} \label{sc:problemStatement}

\subsection{The Stokes problem on a parametrised domain} \label{sc:stokes}

Let us consider a parametrised domain $\Omega^{\bmu} \subset \mathbb{R}^{\nsd}$, where $\nsd$ is the number of spatial dimensions and $\bmu\in\bI\subset\mathbb{R}^{\npar}$ is a set of geometric parameters that controls the boundary representation of the domain, with $\npar$ being the number of geometric parameters.  It is worth noting that the set of geometric parameters can be written as $\bI := \I^1\times\I^2\times\dotsb\times\I^{\npar}$ with $\mu_j\in\I^j$ for $j=1,\dotsc , \npar$.
	
For any set of parameters $\bmu$, the goal is to find the parametric velocity, $\bu(\bx)$, and pressure, $p(\bx)$, fields that satisfy the Stokes problem given by
\begin{equation} \label{eq:stokesStrong}
\left\{\begin{aligned}
-\Div(\nu \Grad \bu - p \Insd) &= \bs  														&&\text{in $\Omega^{\bmu}$,}\\
\Div\bu &= 0  																						&&\text{in $\Omega^{\bmu}$,}\\
\bu &= \bu_D  																						&&\text{on $\Gamma_D^{\bmu}$,}\\
\bn^{\bmu} \cdot \bigl(\nu \Grad \bu -p\Insd \bigr) &= \bt        							&&\text{on $\Gamma_N^{\bmu}$,}\\
\bu \cdot \bD^{\bmu} + \bn^{\bmu} \cdot \bigl(\nu \Grad \bu -p\Insd \bigr)\bE^{\bmu} &= \bm{0}&&\text{on $\Gamma_S^{\bmu}$,}\\
\end{aligned}\right.
\end{equation}
where $\nu>0$ is the kinematic viscosity, $\bs$ is the volumetric source and $\bn^{\bmu}$ is the outward unit normal vector to $\partial\Omega^{\bmu}$. The boundary of the domain, $\partial\Omega^{\bmu}$, is partitioned into the non-overlapping Dirichlet, $\Gamma^{\bmu}_D$, Neumann, $\Gamma^{\bmu}_N$, and slip, $\Gamma^{\bmu}_S$, boundaries such that $\overline{\partial\Omega}^{\bmu}=\overline{\Gamma}^{\bmu}_D\cup\overline{\Gamma}^{\bmu}_N\cup\overline{\Gamma}^{\bmu}_S$. On the Dirichlet boundary the velocity is given by $\bu_D$. On the Neumann boundary the \textit{pseudo-traction} is given by $\bt$. Finally, on the slip boundary, the matrices $\bD^{\bmu}$ and $\bE^{\bmu}$ are given by $\bD^{\bmu} = [\bn^{\bmu},\bm{0}_{\nsd \times (\nsd-1)}]$ and $\bE^{\bmu} = [ \bm{0},\bm{t}_1^{\bmu},...,\bm{t}_{\nsd-1}^{\bmu}]$, as detailed in~\cite{MG-GSH-20}. The tangential vectors $\bm{t}_k^{\bmu}$, for $k=1,\ldots\nsd-1$ are such that $\{\bn^{\bmu},\bm{t}_1^{\bmu},...,\bm{t}_{\nsd-1}^{\bmu}\}$ form an orthonormal system of vectors.

The free divergence condition in equation~\eqref{eq:stokesStrong} induces the compatibility condition
\begin{equation} \label{eq:incompressibleConst}
\langle 1, \bu_D \cdot \bn^{\bmu} \rangle_{\Gamma^{\bmu}_D} + \langle 1, \bu \cdot \bn^{\bmu} \rangle_{\partial \Omega^{\bmu} \setminus \Gamma^{\bmu}_D}  = 0,
\end{equation}
where $\langle \cdot,\cdot \rangle_S$ denotes the standard $\eltwo$ scalar product in any domain $S \subset \partial \Omega^{\bmu}$. 

In addition, it is worth noting that, if $\Gamma^{\bmu}_N = \emptyset$, an additional constraint to avoid the indeterminacy of the pressure is required. One common option~\cite{Jay-CG:09,Cockburn-CS:14,Nguyen-NPC:10,giacomini2018superconvergent} that is considered here, consists of imposing the mean pressure on the boundary of the domain, namely
\begin{equation}
\Big{\langle}  \frac{1}{|\partial \Omega^{\bmu}|} p,1 \Big{\rangle}_{\partial \Omega^{\bmu}} = 0.
\end{equation}

\subsection{The multi-dimensional parametric Stokes problem} \label{sc:stokesMulti}

The classical strategy to solve the parametric Stokes problem is to solve equation~\eqref{eq:stokesStrong} for every set of parameters $\bmu\in\bI$. However, this strategy is not well suited when fast queries are required. 

Reduced order models have demonstrated to be a viable alternative to compute multi-dimensional parametric solutions in an offline phase. Once the offline solution is available, the computation of the solution for a given set of parameters has a very small computational cost, being very well suited for applications where fast queries are required.

The multi-dimensional parametric problem arises from interpreting $\bmu$ as additional parametric coordinates, rather than parameters of the problem. In the context of the Stokes problem considered here, the strategy is to consider the velocity and pressure fields as functions in a multidimensional space, namely $\bu(\bx,\bmu)$ and $p(\bx,\bmu)$. The multi-dimensional parametric Stokes problem can be written as
\begin{equation} 
\left\{\begin{aligned}
-\Div(\nu \Grad \bu - p \Insd) &= \bs  															&&\text{in $\Omega^{\bmu} \times \bI$,}\\
\Div\bu &= 0  																					&&\text{in $\Omega^{\bmu} \times \bI$,}\\
\bu &= \bu_D  																					&&\text{on $\Gamma_D^{\bmu} \times \bI$,}\\
\bn^{\bmu} \cdot \bigl(\nu \Grad \bu -p\Insd \bigr) &= \bt        								&&\text{on $\Gamma_N^{\bmu} \times \bI$,}\\
\bu \cdot \bD^{\bmu} + \bn^{\bmu} \cdot \bigl(\nu \Grad \bu -p\Insd \bigr)\bE^{\bmu} &= \bm{0}  &&\text{on $\Gamma_S^{\bmu} \times \bI$.}
\end{aligned}\right.
\end{equation}

For the multi-dimensional problem, the compatibility condition induced by the free divergence condition can be written as
\begin{equation} \label{eq:incompressibleConstMultidim}
\langle 1, \bu_D \cdot \bn^{\bmu} \rangle_{\Gamma_D^{\bmu} \times \bI} + \langle1, \bu \cdot \bn^{\bmu} \rangle_{ ( \partial \Omega^{\bmu} \setminus \Gamma_D^{\bmu} ) \times \bI}  = 0
\end{equation}
and the additional constraint to avoid the indeterminacy of the pressure, required when $\Gamma_N^{\bmu} = \emptyset$, becomes
\begin{equation}
\Big{\langle}  \frac{1}{|\partial \Omega^{\bmu}|} p,1 \Big{\rangle}_{\partial \Omega^{\bmu}  \times \bI} = 0.
\end{equation}

\section{Hybridisable discontinuous Galerkin formulation} \label{sc:hdg}

Let us consider a subdivision of the domain $\Omega^{\bmu}$ in $\numel$ disjoint subdomains $\Omega_e^{\bmu}$ such that
\begin{equation}
\overline{\Omega}^{\bmu} =  \bigcup_{e=1}^{\numel} \overline{\Omega}_e^{\bmu}.
\end{equation}
The interior boundaries of the subdomains define the so-called mesh skeleton or internal interface $\Gamma^{\bmu}$ as
\begin{equation}\label{eq:Gamma}
\Gamma^{\bmu} := \left[ \bigcup_{e=1}^{\numel} \partial\Omega_e^{\bmu} \right]\setminus\partial\Omega^{\bmu}.
\end{equation}

A partition of the parametric domains $\I^j$, for $j=1,\dotsc, \npar$, in $\numel^j$ disjoint subdomains $\I^j_e$ such that
\begin{equation}
\overline{\I}^j =  \bigcup_{e=1}^{\numel^j} \overline{\I}_e^j,
\end{equation} 
is also considered.

This section briefly presents the HDG formulation for the multi-dimensional parametric Stokes problem. The presentation is based on previous work on HDG methods found in~\cite{Jay-CG:09,Jay-CGNPS:11,Cockburn-CS:14,giacomini2018superconvergent}.

\subsection{Mixed formulation} \label{sc:mixed}

Introducing the so-called \textit{mixed variable} $\bL = -\nu \Grad \bu$, the Stokes problem can be written as a first-order system of equations in the broken computational domain, namely
\begin{equation} \label{eq:StokesMixed}
\left\{\begin{aligned}
\bL_e + \nu \Grad \bu_e & = \bm{0} 																	&&\text{in $\Omega_e^{\bmu} \times \bI$, and for $e=1,\dotsc ,\numel$,}\\
\Div \bigl(\bL_e + p_e\Insd \bigr) &= \bs      														&&\text{in $\Omega_e^{\bmu} \times \bI$, and for $e=1,\dotsc ,\numel$,}\\
\Div \bu_e &= 0  																					&&\text{in $\Omega_e^{\bmu} \times \bI$, and for $e=1,\dotsc ,\numel$,}\\
\bu_e &= \bu_D     																					&&\text{on $\left( \partial \Omega_e^{\bmu} \cap \Gamma_D^{\bmu} \right) \times \bI$,}\\
\bn^{\bmu} \cdot \bigl(\bL_e + p_e \Insd \bigr) &= -\bt    											&&\text{on $\left( \partial \Omega_e^{\bmu} \cap \Gamma_N^{\bmu}\right) \times \bI$,}\\
\bu_e \cdot \bm{D}^{\bmu} - \bn^{\bmu} \cdot \bigl(\bL_e + p_e\Insd \bigr)\bm{E}^{\bmu} &= \bm{0}  	&&\text{on $\left( \partial \Omega_e^{\bmu} \cap \Gamma_S^{\bmu}\right) \times \bI$,}\\
\jump{\bu \otimes \bn^{\bmu}} &=\bm{0}  															&&\text{on $\Gamma^{\bmu} \times \bI$,}\\
\jump{ \bn^{\bmu} \cdot \bigl(\bL + p \Insd \bigr) } &= \bm{0}  									&&\text{on $\Gamma^{\bmu} \times \bI$,}
\end{aligned} \right.
\end{equation}
where the last two equations, known as \textit{transmission conditions}, impose the continuity of the velocity and the normal flux on the mesh skeleton. Following~\cite{AdM-MFH:08}, the \emph{jump} operator $\jump{\cdot}$ is defined as the sum from the left, $\Omega_l$, and right, $\Omega_r$, elements of a given portion of the interface $\Gamma^{\mu} \times \bI$, that is
\begin{equation}
\jump{\odot} = \odot_l + \odot_r .
\end{equation}

\subsection{Strong form of the local and global problems} \label{sc:strongHDG}

The HDG method solves the mixed problem of equation~\eqref{eq:StokesMixed} in two steps. First, the so-called \textit{local problems} are considered 
\begin{equation} \label{eq:localStrong}
\left\{\begin{aligned}
\bL_e + \nu \Grad \bu_e & = \bm{0} 											&&\text{in $\Omega_e^{\bmu} \times \bI$, and for $e=1,\dotsc ,\numel$,}\\
\Div \bigl(\bL_e + p_e\Insd \bigr) &= \bs       							&&\text{in $\Omega_e^{\bmu} \times \bI$, and for $e=1,\dotsc ,\numel$,}\\
\Div \bu_e &= 0  															&&\text{in $\Omega_e^{\bmu} \times \bI$, and for $e=1,\dotsc ,\numel$,}\\
\bu_e &= \bu_D     															&&\text{on $\left( \partial \Omega_e^{\bmu} \cap \Gamma_D^{\bmu} \right) \times \bI$,}\\
\bu_e &= \bhu     															&&\text{on $\left( \partial \Omega_e^{\bmu} \setminus \Gamma_D^{\bmu}\right) \times \bI$,}\\
\Big{\langle} \frac{1}{|\partial \Omega_e^{\bmu}|} p_e,1 \Big{\rangle}_{\partial \Omega_e^{\bmu}  \times \bI} & = \rho_e,	&&\text{for $e=1,\dotsc ,\numel$,}
\end{aligned} \right.
\end{equation}
where $\bhu$ is the so-called hybrid variable, which is an independent variable representing the trace of the solution on the element faces, and $\rho_e$ is the mean value of the pressure on the boundary $\partial\Omega_e$. It is worth noting that the local problem is a pure Dirichlet problem and therefore, the last condition in equation~\eqref{eq:localStrong} is introduced to ensure the uniqueness of the pressure. The local problems can be solved independently, element by element, to write $\bL_e$, $\bu_e$ and $p_e$ in terms of $\bhu$ and $\rho_e$ along the interface $\Gamma^{\bmu} \cup \Gamma_N^{\bmu} \cup \Gamma_S^{\bmu}$.

Second, the so-called \textit{global problem} is defined to impose the continuity of the normal flux on the inter-element faces and the Neumann and slip boundary conditions, namely
\begin{equation} \label{eq:globalStrong}
\left\{\begin{aligned}
\jump{ \bn^{\bmu} \cdot \bigl(\bL + p \Insd \bigr) } &= \bm{0}  									&&\text{on $\Gamma^{\bmu} \times \bI$,}\\
\bn^{\bmu} \cdot \left(\bL_e + p_e \Insd \right) &= -\bt    										&&\text{on $\left( \partial \Omega_e^{\bmu} \cap \Gamma_N^{\bmu}\right) \times \bI$,}\\
\bu_e \cdot \bm{D}^{\bmu} - \bn^{\bmu} \cdot \left(\bL_e + p_e\Insd \right)\bm{E}^{\bmu} &= \bm{0}  &&\text{on $\left( \partial \Omega_e^{\bmu} \cap \Gamma_S^{\bmu}\right) \times \bI$}.
\end{aligned} \right.
\end{equation}

It is worth noting that, due to the unique definition of the hybrid variable on each face and the Dirichlet boundary condition in the local problems, there is no need to enforce the continuity of the solution in the global problem. 

The constraint of equation~\eqref{eq:incompressibleConstMultidim}, induced by the incompressibility condition, is also considered in the global problem and written in terms of the hybrid variable as
\begin{equation} \label{eq:incompressibleConstMultidimUhat}
\langle 1, \bu_D \cdot \bn^{\bmu} \rangle_{\Gamma_D^{\bmu} \times \bI} + \langle1, \bhu \cdot \bn^{\bmu} \rangle_{ ( \partial \Omega^{\bmu} \setminus \Gamma_D^{\bmu} ) \times \bI}  = 0.
\end{equation}

\subsection{Weak form of the local and global problems} \label{sc:weakHDG}

The following discrete functional spaces are introduced:
\begin{align*}
\sVh(\Omega^{\bmu}) &:= \{ v \in \eltwo(\Omega^{\bmu}) : v \vert_{\Omega_e^{\bmu}}\in \Pk(\Omega_e^{\bmu}) \;\forall\Omega_e^{\bmu} \, , \, e=1,\dotsc ,\numel \} , \\
\shVh(S)&:= \{ \hv \in [\eltwo(S)]^{\nsd} : \hv\vert_{\Gamma_i^{\bmu}}\in \Pk(\Gamma_i^{\bmu}) \;\forall\Gamma_i^{\bmu}\subset S\subseteq\Gamma^{\bmu}\cup\partial\Omega^{\bmu} \}, \\
\sLh(\I^j) &:= \{ v\in\eltwo(\I^j) : v\vert_{\I_e^j}\in\Pk(\I_e^j) \;\forall\I_e^j\, , \, e=1,\dotsc ,\numel^j \}, \\
\Lh(\bI) &:= \sLh(\I^1) \otimes \dotsb \otimes \sLh(\I^{\npar}), \\
\sVhMu  & := \sVh(\Omega^{\bmu}) \otimes \Lh(\bI), \\
\hVhMu  & :=  \left[ \shVh(\Gamma^{\bmu} \cup \Gamma_N^{\bmu} \cup \Gamma_S^{\bmu}) \otimes \Lh(\bI) \right]^{\nsd}, \\
\VhMu   & :=  \left[ \sVh(\Omega^{\bmu}) \otimes \Lh(\bI) \right]^{\nsd}, \\
\WhMu   & :=  \left[ \sVh(\Omega^{\bmu}) \otimes \Lh(\bI) \right]^{\nsd \times \nsd}, 
\end{align*}
where $\Pk(\Omega_e^{\bmu})$, $\Pk(\Gamma_i^{\bmu})$ and $\Pk(\I_e^j)$ stand for the spaces of polynomial functions of complete degree at most $k$ in $\Omega_e^{\bmu}$, on $\Gamma_i^{\bmu}$ and in $\I_e^j$ respectively.

The weak form of the local problems, for $e=1,\dotsc ,\numel$, reads: given $\bu_D$ on $\Gamma_D^{\bmu}$ and $\bhu^h$ on $\Gamma^{\bmu} \cup \Gamma_N^{\bmu} \cup \Gamma_S^{\bmu}$, find $(\bL_e^h, \bu_e^h, p_e^h) \in \WhMu \times \VhMu \times \sVhMu$ that satisfy
\begin{equation} \label{eq:localWeak}
	\begin{aligned}	
	A_{LL}(\bW,\bL_e^h) + A_{Lu}(\bW,\bu_e^h)  & = L_L(\bW) + A_{L\hu}(\bW,\bhu^h),  \\
	A_{uL}(\bv,\bL_e^h) + A_{uu}(\bv,\bu_e^h) + A_{up}(\bv,p_e^h) & = L_u(\bv) + A_{u\hu}(\bv,\bhu^h),  \\
	A_{pu}(v,\bu_e^h) & = L_p(v) + A_{p\hu}(v,\bhu^h),  \\
	A_{\rho p}(1,p_e^h) & = A_{\rho \rho}(1,\rho_e^h),
	\end{aligned}
\end{equation}
for all $(\bW,\bv,v) \in \WhMu \times \VhMu \times \sVhMu$, where the multi-dimensional bilinear and linear forms of the local problem are given by
\begin{equation} \label{eq:bilinearLocal}
	\begin{aligned}	
	A_{LL}(\bW,\bL)  	& \!:=\! -\intEMu{\bW}{\nu^{-1} \bL}, 	 		&
	A_{Lu}(\bW,\bu)  	& \!:=\! \intEMu{\Div \bW}{\bu}, 				\\
	A_{L\hu}(\bW,\bhu) 	& \!:=\! \intBNoDMu{\bn^{\bmu} \cdot \bW}{\bhu},&
	A_{uL}(\bv,\bL)  	& \!:=\! \intEMu{\bv}{\Div \bL},				\\
	A_{uu}(\bv,\bu)  	& \!:=\! \intBEMu{\bv}{\btau^{\bmu} \bu},		&
	A_{up}(\bv,p)    	& \!:=\! \intEMu{\bv}{\Grad p},					\\
	A_{u\hu}(\bv,\bhu) 	& \!:=\! \intBNoDMu{\bv}{\btau^{\bmu} \bhu}, 	&
	A_{pu}(v,\bu)    	& \!:=\! \intEMu{\Grad v}{\bu},					\\
	A_{p\hu}(v,\bhu) 	& \!:=\! \intBNoDMu{v}{\bhu \cdot \bn^{\bmu}},	&
	A_{\rho p}(w,p)  	& \!:=\! \intBEMu{w}{|\partial \Omega_e^{\bmu}|^{-1}p},		\\
	A_{\rho \rho}(w,\rho)  	& \!:=\! \intI{w}{\rho},				
	\end{aligned}	
\end{equation}
and
\begin{equation} \label{eq:linearLocal}
	\begin{aligned}	
	L_L(\bW) 	& \!:=\! \intBDMu{\bn^{\bmu} \cdot \bW}{\bu_D}, 				\\
	L_u(\bv) 	& \!:=\! \intEMu{\bv}{\bs} + \intBDMu{\bv}{\btau^{\bmu} \bu_D}, 	\\
	L_p(v)   	& \!:=\! \intBDMu{v}{\bu_D \cdot \bn^{\bmu}},	
	\end{aligned}	
\end{equation}
respectively, where $\left(\cdot,\cdot\right)_{D}$ denotes the standard $\eltwo$ scalar product in a generic subdomain $D$ and $\btau^{\bmu}$ is the stabilisation tensor, whose selection has an important influence on the accuracy, stability and convergence properties of the resulting HDG method~\cite{Jay-CGL:09,Nguyen-NPC:09,Nguyen-NPC:09b,Nguyen-NPC:10}. The choice of the stabilisation tensor for geometrically parametrised problems will be discussed in the next section.

Similarly, the weak form of the global problem is: find $\bhu^h \in \hVhMu$ and $\rho^h \in \mathbb{R}^{\numel} \otimes \Lh(\bI) $ that satisfies
\begin{equation} \label{eq:globalWeak}
	\begin{aligned}	
	\sum_{e=1}^{\numel} \left\{ A_{\hu L}(\bhv,\bL_e^h) + A_{\hu u}(\bhv,\bu_e^h) + A_{\hu p}(\bhv,p_e^h)  + A_{\hu \hu}(\bhv,\bhu^h) \right\} & = \sum_{e=1}^{\numel} \left\{ L_{\hu}(\bhv) \right\},  \\
	A_{p \hu}(1,\bhu^h)  & = -L_p(1),
	\end{aligned}
\end{equation}
for all $\bhv \in \hVhMu$, where the multi-dimensional bilinear and linear forms of the global problem are given by
\begin{equation} 
	\begin{aligned}
	A_{\hu L}(\bhv,\bL) 		 := & \intBNoDSMu{\bhv}{\bn^{\bmu} \cdot \bL}  - \intBSMu{\bhv}{\bn^{\bmu} \cdot \bL \bE^{\bmu}} 					\\
	A_{\hu u}(\bhv,\bu)			 := & \intBNoDSMu{\bhv}{\btau^{\bmu} \bu} - \intBSMu{\bhv}{(\btau^{\bmu} \bu) \!\cdot\! \bE^{\bmu}} 				\\
	A_{\hu p}(\bhv,p)			 := & \intBNoDSMu{\bhv}{p \bn^{\bmu}} 																				\\
	A_{\hu \hu}(\bhv,\bhu)		 := & - \intBNoDSMu{\bhv}{\btau^{\bmu} \bhu}  																		\\
								    & + \intBSMu{\bhv}{\bhu \!\cdot\! \bD^{\bmu} + (\btau^{\bmu} \bhu) \!\cdot\! \bE^{\bmu}}  
	\end{aligned}	
\end{equation}
and 
\begin{equation}
L_{\hu}(\bhv) \!:=\! - \intBNMu{\bhv}{\bt} ,
\end{equation}
respectively.

\section{The proper generalised decomposition strategy} \label{sc:pgd}

The solution of the parametric problem of dimension $\nsd + \npar$, presented in the previous section, with the standard HDG approach is usually not affordable, even for a relatively small number of parameters. To circumvent the \textit{curse of dimensionality}, this section proposes the use of the PGD framework. As it will be shown in this section, the use of an HDG formulation has important advantages compared to other formulations such as standard finite elements~\cite{sevilla2020solution}.

To simplify the presentation, the subindex $_e$ and the superindex  $^h$ used in the previous section to specify the element and the discrete approximations will be omitted here, unless they are needed to follow the development.

\subsection{Separated spatial mapping to obtain generalised solutions} \label{sc:separatedMapping}

As discussed in detail in~\cite{Patera-Rozza:07,Rozza:14,sevilla2020solution}, the solution of the parametric problem described in section~\ref{sc:hdg} requires that the bilinear and linear forms in the weak form can be expressed, or well approximated, by a sum of products of parametric functions and operators that are parameter-independent. To enforce the affine parameter dependence, the integrals appearing in the weak form must involve domains that are not dependent upon the parameters. Following the work of~\cite{AH-AHCCL:14,SZ-ZDMH:15,sevilla2020solution}, a mapping between a parameter-independent reference domain, $\Omega$, and the geometrically parametrised domain is considered, namely
\begin{equation} \label{eq:mappingDomain}
\begin{aligned}
\mapping\, :
\Omega\times\bI & \longrightarrow \Omega^{\bmu} \\
(\bX,\bmu) 		& \longmapsto \bx = \mapping(\bX,\bmu).
\end{aligned}
\end{equation}
The coordinates of the reference, or \textit{undeformed}, domain are denoted by $\bX$ whereas the coordinates of the parametric, or \textit{deformed}, domain are denoted by $\bx$. To ensure the affine parameter dependence, the mapping is assumed to be given in separated form as
\begin{equation} \label{eq:displacementSep}
\mapping(\bX,\bmu) = \sum_{k=1}^{\nmap} \bM^k(\bX) \phi^k(\bmu).
\end{equation}

\begin{remark} \label{rk:analyticalMapping}
To simplify the presentation here, it is assumed that the separated representation of the mapping is given analytically. As mentioned earlier, a general strategy to construct a separable mapping was described in~\cite{sevilla2020solution} using an exact boundary description of the computational domain by means of NURBS.
\end{remark}

The separated representation of the mapping leads to the following separated representation of its Jacobian
\begin{equation}\label{eq:JacoSep}
\Jaco(\bX,\bmu) = \frac{\partial \bx}{\partial\bX}(\bX,\bmu) =  \sum_{k=1}^{\nmap} \Jk(\bX) \phi^k(\bmu).
\end{equation}
In addition, the separated description of the mapping and its Jacobian can be used to obtain a separated expression of the determinant and the adjoint of the Jacobian using the Leibniz formula and the Leverrier's algorithm as explained in detail in~\cite{sevilla2020solution}. 
The separated expression of the determinant of the Jacobian and its adjoint are written in compact form as
\begin{equation} \label{eq:DeterminantSep}
\det(\Jaco)(\bX,\bmu) = \sum_{k=1}^{\ndet} D^k(\bX) \theta^k(\bmu)
\end{equation}
and
\begin{equation} \label{eq:DeterminantAdj}
\adj(\Jaco)(\bX,\bmu) = \sum_{k=1}^{\nadj} \bA^k(\bX) \vartheta^k(\bmu) ,
\end{equation}
respectively.

\subsection{Affine parameter dependence of the HDG bilinear and linear forms} \label{sc:HDG-PGD}

Introducing the mapping $\mapping$ of equation~\eqref{eq:mappingDomain} into the weak form of the local and global problems, it is possible to write the integrals over the reference domain, $\Omega$, and its boundary, $\partial\Omega$, not dependent on the parameters $\bm{\mu}$. The bilinear and linear forms for the local problems can be written as
\begin{equation} \label{eq:bilinearLocalRef}
	\begin{aligned}	
	A_{LL}(\bW,\bL)  	& = -\intExI{\bW}{\nu^{-1} \det{(\Jaco) \bL } }, 	\\
	A_{Lu}(\bW,\bu)  	& = \intExI{ \adj{(\Jaco)}\DivX \bW}{\bu}, 			\\
	A_{L\hu}(\bW,\bhu) 	& = \intBNoDxI{\adj{(\Jaco)} \bn \cdot \bW}{\bhu}, 	\\
	A_{uL}(\bv,\bL)  	& = \intExI{\bv}{\adj{(\Jaco)}\DivX \bL},			\\
	A_{uu}(\bv,\bu)  	& = \intBExI{\bv}{\btau \bu},						\\
	A_{up}(\bv,p)    	& = \intExI{\bv}{\adj{(\Jaco)} \GradX p},			\\
	A_{u\hu}(\bv,\bhu) 	& = \intBNoDxI{\bv}{\btau \bhu}, 					\\
	A_{pu}(v,\bu)    	& = \intExI{\adj{(\Jaco)} \GradX v}{\bu},			\\
	A_{p\hu}(v,\bhu) 	& = \intBNoDxI{v}{\bhu \cdot \adj{(\Jaco)} \bn},	\\
	A_{\rho p}(w,p)  	& = \intBExI{w}{|\partial \Omega_e|^{-1}p},			\\
	A_{\rho \rho}(w,\rho)  	& = \intI{w}{\rho},						
	\end{aligned}	
\end{equation}
and
\begin{equation} \label{eq:linearLocalRef}
	\begin{aligned}	
	L_L(\bW) 	& = \intBDxI{\adj{(\Jaco)} \bn \cdot \bW}{\bu_D}, 									\\
	L_u(\bv) 	& = \intExI{\bv}{\det{(\Jaco)} \bs} + \intBDxI{\bv}{\btau \bu_D}, 	\\
	L_p(v)   	& = \intBDxI{v}{\bu_D \cdot \adj{(\Jaco)} \bn},				
	\end{aligned}	
\end{equation}
respectively, where the adjoint operator is defined as $\adj(\bm{A}) =\det(\bm{A})\,\bm{A}^{-1}$ and the stabilisation parameter in the deformed domain is chosen as
\begin{equation} \label{eq:tauDeformed}
\btau^\bmu := \frac{1}{\| \adj(\Jaco) \bn \| }\btau.
\end{equation}
The scaling factor $\|\! \adj(\Jaco) \bn \|$ in equation~\eqref{eq:tauDeformed} accounts for the increased or decreased area of the deformed face, $\partial\Omega_e^\bmu$, with respect to the reference one, $\partial\Omega_e$. This definition, inspired by the expression of the penalty coefficient in classical interior penalty DG methods~\cite{Arnold1982}, ensures that the larger the deformation of the face, the  smaller the value of $\btau^\bmu$ is. This ensures that a weaker continuity is imposed for large deformations and it is justified by the expected loss of accuracy in the hybrid variable when the mapping introduces a large deformation.

Following previous work on HDG methods for Stokes problems~\cite{giacomini2018superconvergent}, the stabilisation parameter in the reference domain is selected as $\btau = (\tau \nu/\ell) \Insd$, where $\tau$ is a numerical parameter, selected as $\tau=10$ in this work, and $\ell$ is a characteristic length of the domain.

\begin{remark} \label{rk:rhoConstraint}
As mentioned above, it holds that $ \| \adj(\Jaco) \bn \|=|\partial \Omega_e^{\bmu}|/|\partial \Omega_e|$. Hence, no parametric dependence appears in the arguments of the bilinear form $A_{\rho p}$.
\end{remark}

Analogously, the bilinear and linear forms for the global problem can be written as
\begin{equation} \label{eq:bilinearGlobalRef}
	\begin{aligned}
	A_{\hu L}(\bhv,\bL) 		 = & \intBNoDSxI{\bhv}{\adj{(\Jaco)} \bn \cdot \bL} 	\\
	 &- \intBSxI{\bhv}{\adj{(\Jaco)} \bn \cdot \bL \bE} 												\\
	A_{\hu u}(\bhv,\bu)		 = & \intBNoDSxI{\bhv}{\btau \bu} - \intBSxI{\bhv}{(\btau \bu) \!\cdot\! \bE} 	\\
	A_{\hu p}(\bhv,p)			 = & \intBNoDSxI{\bhv}{p \adj{(\Jaco)} \bn}  \\
	A_{\hu \hu}(\bhv,\bhu) = & - \intBNoDSxI{\bhv}{\btau \bhu} 											\\
	 &+ \intBSxI{\bhv}{\bhu \!\cdot\! \adj{(\Jaco)} \bD + (\btau \bhu) \!\cdot\! \bE}
	\end{aligned}	
\end{equation}
and 
\begin{equation} \label{eq:linearGlboalRef}
 L_{\hu}(\bhv)	 = - \intBNxI{\bhv}{\bt} ,
\end{equation}
respectively.

\begin{remark} \label{rk:GammaS}
The derivation of the terms on the slip boundary in~\eqref{eq:bilinearGlobalRef} follows from the relationship $\intBSMu{\bhv}{\bn^{\bmu} \cdot \bF} = \intBSxI{\bhv}{\adj{(\Jaco)} \bn \cdot \bF}$ and the definition~\eqref{eq:tauDeformed}. The slip boundary condition is used here to enforce a symmetry condition and therefore, it is assumed that the orientation of the vectors $\{\bn^{\bmu},\bm{t}_1^{\bmu},...,\bm{t}_{\nsd-1}^{\bmu}\}$ is preserved by the mapping $\mapping$. It is worth noting that this does not imply that $\Gamma_S^\bmu=\Gamma_S$ as it will be shown with numerical examples. 
\end{remark}

\begin{remark} \label{rk:GammaN}
As classical in the context of shape optimisation~\cite{Allaire-book}, in~\eqref{eq:linearGlboalRef} it is assumed that Neumann boundaries, where a traction (or pseudo-traction) is imposed, are fixed, that is, $\Gamma_N^\bmu=\Gamma_N$. On the contrary, deformable Neumann boundaries, also known as \emph{free boundaries}, are traction-free, whence $\bt$  is null.
\end{remark}

\subsection{Separated representation of the data} \label{sc:separatedData}

As usual in a PGD context, the data is assumed to be given in separated form. For the Stokes problem under consideration, this means that the Dirichlet and Neumann data and the source term can be written as
\begin{equation} \label{eq:separatedData}
	\begin{aligned}	
	\bu_D& \!=\! \sum_{l=1}^{\nDir} \bgD^l(\bX) \lambda^l_D(\bmu),
	\\ 
	\bt	 & \!=\! \sum_{l=1}^{\nNeu} \bgN^l(\bX) \lambda^l_N(\bmu),
	\\ 
	\bs	 & \!=\! \sum_{l=1}^{\nSou} \bgS^l(\bX) \lambda^l_S(\bmu).
	\end{aligned}	
\end{equation}
Even if the data is not directly given in this form, it is possible to obtain a good approximation in a separated form, see~\cite{Chinesta-Keunings-Leygue}.

\subsection{Separated representation of the primal, mixed and hybrid variables} \label{sc:separatedSols}

Following the \emph{predictor-corrector} PGD rationale, see~\cite{tsiolakis2020nonintrusive}, each variable of the HDG formulation, presented in section~\ref{sc:hdg}, is written as a rank-$m$ separable approximation, that is 
\begin{equation}\label{eq:PGDapprox}
	\begin{aligned}
	\Lpgd^m	 (\bX,\bmu) &= \ampL^m  [\fL^m (\bX) \, \psi^m(\bmu)  + \dLpgd^m(\bX, \bmu)] + \Lpgd^{m-1}  (\bX, \bmu) , \\
	\bupgd^m (\bX,\bmu) &= \ampU^m [\fU^m (\bX) \, \psi^m(\bmu) + \dbupgd^m(\bX, \bmu)] + \bupgd^{m-1} (\bX, \bmu) , \\
	\ppgd^m  (\bX,\bmu) &= \ampP^m  [\fP^m (\bX) \, \psi^m(\bmu) + \dppgd^m(\bX, \bmu)] + \ppgd^{m-1}  (\bX, \bmu) , \\
	\bhupgd^m(\bX,\bmu) &= \ampHU^m [\fHU^m(\bX) \, \psi^m(\bmu) + \dbhupgd^m(\bX, \bmu)] + \bhupgd^{m-1}(\bX, \bmu) , \\
	\rpgd^m  (\bX,\bmu) &= \ampR^m  [\fR^m (\bX) \, \psi^m(\bmu) + \drpgd^m(\bX, \bmu)] + \rpgd^{m-1}  (\bX, \bmu) ,
	\end{aligned}
\end{equation}
where $\ampL^m  \fL^m \psi^m$, $\ampU^m \fU^m \psi^m$, $\ampP^m  \fP^m \psi^m$, $\ampHU^m \fHU^m \psi^m$ and $\ampR^m  \fR^m \psi^m$ are the predictors of the $m$-th mode in the PGD expansion, whereas $\ampL^m  \dLpgd^m$, $\ampU^m \dbupgd^m$, $\ampP^m  \dppgd^m$, $\ampHU^m \dbhupgd^m$ and $\ampR^m  \drpgd^m$ are the corresponding correction terms. Introducing the variation $\De$, the correctors are defined as
\begin{equation}\label{eq:PGDcorrections}
	\begin{aligned}
	\dLpgd^m	 (\bX,\bmu) &:= \De\fL (\bX) \, \psi^m(\bmu)  + \fL^m (\bX) \, \De\psi(\bmu) + \De\fL (\bX) \, \De\psi(\bmu) , \\
	\dbupgd^m (\bX,\bmu) &:= \De\fU (\bX) \, \psi^m(\bmu)  + \fU^m (\bX) \, \De\psi(\bmu) + \De\fU (\bX) \, \De\psi(\bmu) , \\
	\dppgd^m  (\bX,\bmu) &:= \De\fP (\bX) \, \psi^m(\bmu)  + \fP^m (\bX) \, \De\psi(\bmu) + \De\fP (\bX) \, \De\psi(\bmu) , \\
	\dbhupgd^m(\bX,\bmu) &:= \De\fHU (\bX) \, \psi^m(\bmu)  + \fHU^m (\bX) \, \De\psi(\bmu) + \De\fHU (\bX) \, \De\psi(\bmu) , \\
	\drpgd^m  (\bX,\bmu) &:= \De\fR (\bX) \, \psi^m(\bmu)  + \fR^m (\bX) \, \De\psi(\bmu) + \De\fR (\bX) \, \De\psi(\bmu) ,
	\end{aligned}
\end{equation}
where the least term denotes a high-order variation and it is henceforth neglected.

Each term, or \emph{mode}, of the PGD approximation is the product of a function that depends upon the spatial coordinates and a function that depends upon the parameters. In addition, the parametric functions are assumed to be the product of functions that depend upon a single parameter, namely
\begin{equation} \label{eq:phiM}
\psi^m(\bmu) = \prod_{j=1}^{\npar} \psi^m_j(\mu_j).
\end{equation}
As usual in a PGD context, the number of terms is a priori unknown and it is determined using a greedy algorithm. Assuming that $m-1$ modes are known, the computation of the $m$-th mode requires the solution of a nonlinear problem as described in the next section.

\begin{remark} \label{rk:singleParameterPGD}
This work considers the so-called \emph{single-parameter} approach, where the parametric function of the $m$-th mode, $\psi^m$, is the same for all the variables. Other approaches, including a different parametric function for each variable or even the use of vector-valued parametric functions in the approximation of vector fields are discussed in~\cite{diez2017generalized}.	
\end{remark}


The tangent manifold for $\bL$ is characterised by choosing $\bW$ as variations of $\fL$ and $\psi$, that is
\begin{equation}\label{eq:tangentL}
\bW = \de \fL \psi^m + \ampL^m \fL^m  \de \psi ,
\end{equation}
for $\de \fL \in \Wh \!:=  \left[ \sVh(\Omega)  \right]^{\nsd \times \nsd} $ and $\de \psi \in \Lh(\bI)$. Similarly, the tangent manifolds for $\bu$, $p$, $\bhu$ and $\rho$ are characterised by choosing 
\begin{equation}\label{eq:tangentU}
\begin{aligned}
\bv  &= \de \fU \psi^m + \ampU^m \fU^m  \de \psi ,			&
v    &= \de \fP \psi^m + \ampP^m \fP^m  \de \psi ,			\\
\bhv &= \de \fHU \psi^m + \ampHU^m \fHU^m  \de \psi ,		&
w    &= \de \fR \psi^m + \ampR^m \fR^m \de \psi, 			
\end{aligned}
\end{equation}
for $\de \fU \in \Vh \!:=  \left[ \sVh(\Omega)  \right]^{\nsd}$, $\de \fP \in \sVh$, $\de\fHU \in \hVh \!:=  \left[ \shVh(\Gamma \cup \Gamma_N \cup \Gamma_S)  \right]^{\nsd}$ and $\de \fR \in \mathbb{R}^{\numel}$.

\subsection{Alternating direction scheme} \label{sc:AD}

With the separated structure of the PGD approximations, the weighting functions and the bilinear and linear HDG forms described in the previous sections, it is possible to drastically reduce the complexity of the problem by projecting the high-dimensional problem on the tangent manifold and applying an alternating direction strategy. First, in the so-called \textit{spatial iteration}, the parametric function of the $m$-th mode is assumed known and the spatial functions are determined. As it will be shown, this step requires to solve a system of equations with a very similar structure to the non-parametric HDG problem. Second, in the so-called \textit{parametric iteration}, the parametric function is computed using the spatial functions determined in the first step. This process is repeated until convergence is achieved. It is worth noting that the order of the spatial and parametric iterations can be swapped without affecting the alternating direction algorithm.

Let us assume that we have computed the first $m-1$ modes and it is of interest to compute the $m$-th mode. In the next two sections, the alternating direction strategy to compute the spatial and parametric modes is detailed. 

\subsubsection{The spatial iteration} \label{sc:spatialIteration}

In the spatial iteration, it is assumed that the parametric function $\psi^m$ and the spatial predictions $\ampL^m \fL^m$, $\ampU^m \fU^m$, $\ampP^m \fP^m$, $\ampHU^m \fHU^m$ and $\ampR^m \fR^m$ are known and the goal is to compute the corresponding corrections $\ampL^m \De \fL$, $\ampU^m \De \fU$, $\ampP^m \De \fP$, $\ampHU^m \De \fHU$ and $\ampR^m \De \fR$. 

Taking into account that $\de\psi=0$ when $\psi^m$ is known and introducing the expression of the PGD approximations and the weighting functions in the weak form of the HDG local problems, the following weak form of the local problem for the spatial iteration is obtained: find $(\ampL^m \De \fL, \ampU^m \De \fU, \ampP^m \De \fP) \in \Wh \times \Vh \times \sVh$ that satisfy
\begin{equation}\label{eq:localWeakPGDSpatial}
	\begin{aligned}
	\sum_{k=1}^{\ndet} \beta_{\theta}^k \mathcal{A}^k_{LL}(\de \fL, \ampL^m \De \fL) & + \sum_{k=1}^{\nadj} \beta^k_{\vartheta} \mathcal{A}^k_{Lu}(\de \fL, \ampU^m \De \fU)	\\
	= & \mathcal{R}_L^m (\de \fL \psi^m) 
	+ \sum_{k=1}^{\nadj} \beta^k_{\vartheta} \mathcal{A}^k_{L\hu}(\de \fL, \ampHU^m \De \fHU),  
	\\	
	\sum_{k=1}^{\nadj} \beta^k_{\vartheta} \mathcal{A}^k_{uL}(\de \fU, \ampL^m \De \fL)  + \beta & \mathcal{A}_{uu} (\de \fU, \ampU^m \De \fU) \\
	+ \sum_{k=1}^{\nadj} \beta^k_{\vartheta} \mathcal{A}^k_{up}(\de \fU, \ampP^m \De \fP)  = & \mathcal{R}_u^m (\de \fU \psi^m) + \beta \mathcal{A}_{u\hu}(\de \fU, \ampHU^m \De \fHU),
	\\
	\sum_{k=1}^{\nadj} \beta^k_{\vartheta} \mathcal{A}^k_{pu}(\de \fP, \ampU^m \De \fU)  = & \mathcal{R}_p^m (\de \fP \psi^m) + \sum_{k=1}^{\nadj} \beta^k_{\vartheta} \mathcal{A}^k_{p\hu}(\de \fP, \ampHU^m \De \fHU) 
	\\
	\beta \mathcal{A}^k_{\rho p}(1, \ampP^m \De \fP)  = & \mathcal{R}_{\overline{p}}^m (\psi^m) + \beta \mathcal{A}_{\rho \rho}(1, \ampR^m \De \fR),
	\end{aligned}
\end{equation}
for all $(\de \fL,\de \fU,\de \fP) \in \Wh \times \Vh \times \sVh$.

The bilinear and linear forms of the local problem are detailed in equation~\eqref{eq:bilinearPGDSpatialLocal}, in~\ref{sc:bilinearPGD}, and equation~\eqref{eq:linearLocalPGDSpatial}, in~\ref{sc:linearPGD}, respectively. The constants in equation~\eqref{eq:localWeakPGDSpatial} are given by 
\begin{equation} \label{eq:cttLocalPGDSpatial}
\beta^k_{\theta}    :=  \mathcal{A}^k_{\theta}(\psi^m,\psi^m) 		\qquad
\beta^k_{\vartheta} :=  \mathcal{A}^k_{\vartheta}(\psi^m,\psi^m),  	\qquad
\beta    			:=  \mathcal{A}(\psi^m,\psi^m) , 				
\end{equation}
where the bilinear forms involved in the definitions of these constants are introduced in equation~\eqref{eq:bilinearPGDParam}, in~\ref{sc:bilinearPGD}.

As mentioned earlier, in remark~\ref{rk:singleParameterPGD}, this work considers the same parametric function for all the variables. It is worth noting that this choice reduces the number of different constants in equation~\eqref{eq:localWeakPGDSpatial}.

Similarly, the weak form of the global problem is: find $\ampHU^m \De \fHU \in \hVh$ and $\ampR^m \De \fR \in \mathbb{R}^{\numel}$ that satisfy
\begin{subequations}
\begin{equation}\label{eq:globalWeakPGDSpatialTransmission}
	\begin{aligned}	
		\sum_{e=1}^{\numel} \Biggl\{ \sum_{k=1}^{\nadj}\beta^k_{\vartheta} \mathcal{A}^k_{\hu L}  (\de \fHU, \ampL^m \De \fL) 	
	 + \beta \mathcal{A}_{\hu u}(\de \fHU, \ampU^m \De \fU) & \\
	 + \sum_{k=1}^{\nadj} \beta^k_{\vartheta} \mathcal{A}^k_{\hu p}(\de \fHU, \ampP^m \De \fP)  
     + \beta \mathcal{A}_{\hu \hu}(\de \fHU, & \ampHU^m \De \fHU) \\
	  + \sum_{k=1}^{\nadj} \beta^k_{\vartheta}\mathcal{A}^k_{\hu \hu} (\de \fHU, \ampHU^m \De \fHU) \Biggr\}
	& = \sum_{e=1}^{\numel} \mathcal{R}_{\hu}^m (\de \fHU \psi^m),  
	\end{aligned}
\end{equation}
for all $\de \fHU \in \hVh$, with the incompressibility constraint 
\begin{equation}\label{eq:globalWeakPGDSpatialConstraint}
	\sum_{k=1}^{\nadj} \beta^k_{\vartheta} \mathcal{A}^k_{p \hu} (1,\ampHU^m \De \fHU)  = \mathcal{R}_{\rho}^m (\psi^m), \quad e=1,\dots,\numel .
\end{equation}
\end{subequations}
The bilinear and linear forms of the global problem are detailed in equation~\eqref{eq:bilinearPGDSpatialGlobal}, in~\ref{sc:bilinearPGD}, and equation~\eqref{eq:linearGlobalPGDSpatial}, in~\ref{sc:linearPGD}, respectively.

\subsubsection{The parametric iteration} \label{sc:parametricIteration}

After computing the spatial corrections following the procedure described in the previous section, the spatial modes are updated, namely
\begin{equation}\label{eq:newSpatialModes}
\begin{aligned}
\ampL^m \fL^m &\gets \ampL^m \fL^m + \ampL^m \De \fL ,\\
\ampU^m \fU^m &\gets \ampU^m \fU^m + \ampU^m  \De \fU , \\
\ampP^m \fP^m &\gets \ampP^m \fP^m + \ampP^m \De \fP , \\
\ampHU^m \fHU^m &\gets \ampHU^m \fHU^m + \ampHU^m \De \fHU , \\
\ampR^m \fR^m &\gets \ampR^m \fR^m + \ampR^m  \De \fR ,
\end{aligned}
\end{equation}
where the constant $\sigma_{\diamond}^m$ on the left hand side denotes the amplitude of the newly computed $m$-th mode of the function $\diamond$, e.g. $\ampP^m \gets \| \ampP^m \fP^m + \ampP^m \De \fP \|$.

In the parametric iteration, the goal is to compute the parametric correction $\De\psi$ given the prediction $\psi^m$ and the known spatial functions in~\eqref{eq:newSpatialModes}. Following the assumption that such functions are known, it holds that $\de\fL = \de\fU = \de\fP = \de\fHU = \de\fR = 0$. Introducing the expression of the PGD approximations and the weighting functions in the weak form of the HDG local problems, the following weak form of the local problem for the spatial iteration is obtained: find $\De\psi \in \Lh(\bI)$ such that 
\begin{equation}\label{eq:localWeakPGDParam}
	\begin{aligned}	
	\sum_{k=1}^{\ndet} \gamma_{LL}^k \mathcal{A}^k_{\theta}(\de\psi,\De\psi)  & + \sum_{k=1}^{\nadj} \gamma_{Lu}^k \mathcal{A}^k_{\vartheta}(\de\psi,\De\psi) \\
	 = & \mathcal{R}_L^m (\ampL^m \fL^m \de \psi)  
	 + \sum_{k=1}^{\nadj} \gamma_{L\hu}^k \mathcal{A}^k_{\vartheta}(\de\psi,\De\psi),  
	\\
	\sum_{k=1}^{\nadj} \gamma_{uL}^k \mathcal{A}^k_{\vartheta}(\de\psi,\De\psi)  + 
	\gamma_{uu} & \mathcal{A}(\de\psi,\De\psi) \\
	 + 	\sum_{k=1}^{\nadj} \gamma_{up}^k \mathcal{A}^k_{\vartheta}(\de\psi,\De\psi) = & \mathcal{R}_u^m (\ampU^m \fU^m \de \psi)  
	+  \gamma_{u\hu} \mathcal{A}(\de\psi,\psi^m),  
	\\
	\sum_{k=1}^{\nadj} \gamma_{pu}^k \mathcal{A}^k_{\vartheta}(\de\psi,\De\psi)  = & \mathcal{R}_p^m (\ampP^m \fP^m \de \psi)  + \sum_{k=1}^{\nadj} \gamma_{p\hu}^k \mathcal{A}^k_{\vartheta}(\de\psi,\De\psi),  \\
	\gamma_{\rho p} \mathcal{A}(\de \psi,\De\psi)  = & \mathcal{R}_{\overline{p}}^m (\de \psi) + \gamma_{\rho \rho} \mathcal{A}(\de \psi,\De\psi),
	\end{aligned}
\end{equation}
for all $\de \psi \in \Lh(\bI)$,

Similarly, the weak form of the global problem is: find $\De\psi \in \Lh(\bI)$ that satisfies
\begin{equation}\label{eq:globalWeakPGDParam}
\begin{aligned}	
\sum_{e=1}^{\numel} \left\{ \sum_{k=1}^{\nadj} \right. \gamma^k_{\hu L} \mathcal{A}_{\vartheta}^k(\de\psi,\De\psi)
 + \gamma_{\hu u} \mathcal{A}(\de\psi,\De\psi)  + \sum_{k=1}^{\nadj} & \gamma^k_{\hu p} \mathcal{A}_{\vartheta}^k(\de\psi,\De\psi) \\
 \left. + \gamma_{\hu \hu} \mathcal{A}(\de\psi,\De\psi)  + \sum_{k=1}^{\nadj} \gamma^k_{\hu \hu} \mathcal{A}_{\vartheta}^k(\de\psi,\De\psi) \right\}
= & \sum_{e=1}^{\numel} \mathcal{R}_{\hu}^m (\ampHU^m \fHU^m \de \psi),  
\\
\sum_{k=1}^{\nadj} \gamma_{\rho \hu}^k \mathcal{A}_{\vartheta}^k (\de\psi,\psi^m)  = & \mathcal{R}_{\rho}^m (\de\psi),
\end{aligned}
\end{equation}
for all $\de \psi \in \Lh(\bI)$.

The constants in equations~\eqref{eq:localWeakPGDParam} and~\eqref{eq:globalWeakPGDParam} are defined as
\begin{equation} \label{eq:cttPGDParam}
\begin{aligned}	
\gamma_{LL}^k 	& := \mathcal{A}^k_{LL}		(\ampL^m \fL^m, \ampL^m \fL^m),		&
\gamma_{Lu}^k 	& := \mathcal{A}^k_{Lu}		(\ampL^m \fL^m, \ampU^m \fU^m), 	\\
\gamma_{L\hu}^k & := \mathcal{A}^k_{L\hu}	(\ampL^m \fL^m, \ampHU^m \fHU^m), 	&
\gamma_{uL}^k 	& := \mathcal{A}^k_{uL}		(\ampU^m\fU^m, \ampL^m \fL^m),		\\
\gamma_{uu} 	& := \mathcal{A}_{uu}		(\ampU^m\fU^m, \ampU^m\fU^m),		&
\gamma_{up}^k	& := \mathcal{A}^k_{up}		(\ampU^m\fU^m, \ampP^m\fP^m),		\\
\gamma_{u\hu}	& := \mathcal{A}_{u\hu}		(\ampU^m\fU^m, \ampHU^m\fHU^m),	&
\gamma_{pu}^k	& := \mathcal{A}^k_{pu}		(\ampP^m \fP^m, \ampU^m\fU^m),		\\
\gamma_{p\hu}^k	& := \mathcal{A}^k_{p\hu}	(\ampP^m\fP^m, \ampHU^m\fHU^m),  	&
\gamma_{\rho p}	& := \mathcal{A}_{\rho p}	(1, \ampP^m\fP^m), 	\\
\gamma_{\rho \rho}& := \mathcal{A}_{\rho \rho}(1, \ampR^m\fR^m), & & \\
\gamma_{\hu L}^k 	& := \mathcal{A}^k_{\hu L}		(\ampHU^m\fHU^m, \ampL^m \fL^m),		&
\gamma_{\hu u} 	& := \mathcal{A}_{\hu u}		(\ampHU^m\fHU^m, \ampU^m\fU^m),		\\
\gamma_{\hu p}^k	& := \mathcal{A}^k_{\hu p}		(\ampHU^m\fHU^m, \ampP^m\fP^m),		&
\gamma_{\hu \hu}	& := \mathcal{A}_{\hu \hu}		(\ampHU^m\fHU^m, \ampHU^m\fHU^m),	\\
\gamma_{\hu \hu}^k	& := \mathcal{A}_{\hu\hu}^k		(\ampHU^m\fHU^m, \ampHU^m\fHU^m),	&
\gamma_{\rho\hu}^k	& := \mathcal{A}_{p\hu}^k		(1, \ampHU^m\fHU^m).
\end{aligned}	
\end{equation}

The choice of a single parameter approximation implies that we can combine equations~\eqref{eq:localWeakPGDParam} and~\eqref{eq:globalWeakPGDParam} to obtain the following parametric problem: find $\De\psi \in \Lh(\bI)$ that satisfies
\begin{equation}\label{eq:WeakPGDParam}
\sum_{k=1}^{\ndet} \gamma_{LL}^k \mathcal{A}^k_{\theta}(\de\psi,\De\psi) +  \sum_{k=1}^{\nadj} \gamma_{\vartheta}^k \mathcal{A}_{\vartheta}^k(\de\psi,\De\psi) + \gamma \mathcal{A}(\de \psi,\De\psi) =
\mathcal{R}^m(\de\psi),
\end{equation}
for all $\de \psi \in \Lh(\bI)$, where
\begin{equation}\label{eq:WeakPGDParamCtt}
\begin{aligned}
\gamma_{\vartheta}^k 	 :=& \gamma_{Lu}^k - \gamma_{L\hu}^k + \gamma_{uL}^k + \gamma_{up}^k + \gamma_{pu}^k - \gamma_{p \hu}^k + \gamma^k_{\hu L} + \gamma^k_{\hu p} + \gamma^k_{\hu \hu} + \gamma_{\rho \hu}^k , \\
\gamma   				 :=& \gamma_{uu} - \gamma_{u\hu} + \gamma_{\rho p} - \gamma_{\rho \rho} + \gamma_{\hu u} + \gamma_{\hu \hu}, \\
\mathcal{R}^m(\de\psi) 			  :=& \mathcal{R}_L^m (\ampL^m \fL^m \de \psi) + \mathcal{R}_u^m (\ampU^m \fU^m \de \psi) + \mathcal{R}_p^m (\ampP^m \fP^m \de \psi)\\
&+ \mathcal{R}_{\overline{p}}^m (\de \psi) + \mathcal{R}_{\hu}^m (\ampHU^m \fHU^m \de \psi) + \mathcal{R}_{\rho}^m (\de \psi).
\end{aligned}
\end{equation}

\begin{remark}\label{rk:parProb}
Alternative formulations of the parametric problem may be devised, e.g. by considering only equation~\eqref{eq:localWeakPGDParam} or \eqref{eq:globalWeakPGDParam}. In this work, equation~\eqref{eq:WeakPGDParam} has been considered in the parametric iteration in order to account for the information of both the local and the global HDG problems.
\end{remark}

As detailed in equation~\eqref{eq:phiM}, the parametric iteration involves $\npar$ geometric parameters. To reduce the size of the problem of the parametric iteration, $\npar$ one-dimensional problems are solved sequentially, as commonly done in a PGD framework~\cite{PGD-CCH:14}.

\subsection{The HDG-PGD algorithm} \label{sc:algorithm}

The HDG solver for geometrically parametrised Stokes equation is described in algorithm~\ref{alg:HDG-PGD}. Differently from traditional PGD strategies relying on continuous Galerkin approximations, Dirichlet boundary conditions do not require a special treatment in the context of HDG-PGD. More precisely, Dirichlet conditions are imposed in a weak sense and appear in the linear forms~\eqref{eq:linearLocal} of the HDG local problem.

\begin{algorithm}
\caption{The HDG-PGD implementation}\label{alg:HDG-PGD}
\begin{algorithmic}[1]
\REQUIRE{For the greedy enrichment loop, the value $\eta^\star$ of the tolerance. For the alternating direction iterations, the values $\eta_{\hu}$ and $\eta_\circ^r$ of the tolerances for the mode amplitude $\ampHU$ and the residuals $r_{\circ}$ obtained from the linear forms in~\ref{sc:linearPGD}, respectively. For the spatial and parametric problems, the typical values $\text{typ}_\circ$ of the residuals. $\circ=\hu,\psi$.}
\STATE{Set $m \gets 1$ and initialise the amplitude of the spatial mode $\ampHU^1 \gets 1$.}

\WHILE{$\ampHU^m > \eta^\star\,\ampHU^1$}
\STATE{Set $q \gets 1$ and initialise the parametric predictor $\psi^m {\gets} 1$.}
\STATE{Compute the spatial constants \eqref{eq:cttLocalPGDSpatial}.}
\STATE{Solve the HDG global \eqref{eq:globalWeakPGDSpatialTransmission}-\eqref{eq:globalWeakPGDSpatialConstraint} and local problems \eqref{eq:localWeakPGDSpatial}.}

\STATE{Initialise $\varepsilon_{\hu} \gets 1$, $\varepsilon_\circ^r \gets \text{typ}_\circ$.}
             
\WHILE{$\varepsilon_{\hu} > \eta_{\hu}$ or $\varepsilon_\circ^r > \eta_\circ^r$}

\STATE{Compute the parametric constants \eqref{eq:cttPGDParam}.}
\STATE{Solve the parametric linear system \eqref{eq:WeakPGDParam}.}
\STATE{Update the parametric predictor $\psi^m {\gets} (\psi^m + \De\psi)/\norm{\psi^m+\De\psi}$.}

\STATE{Compute the spatial constants \eqref{eq:cttLocalPGDSpatial}.}
\STATE{Solve the HDG global \eqref{eq:globalWeakPGDSpatialTransmission}-\eqref{eq:globalWeakPGDSpatialConstraint} and local problems \eqref{eq:localWeakPGDSpatial}.}
\STATE{Normalise the spatial predictor $\ampHU^m {\gets} \norm{\ampHU^m \fHU^m  + \ampHU^m \De\fHU}$.}
\STATE{Update the spatial predictor $\ampHU^m \fHU^m  {\gets} \ampHU^m \fHU^m  + \ampHU^m \De \fHU$.}
             
\STATE{Update the stopping criteria for the mode amplitude $\varepsilon_{\hu} {\gets} \norm{\ampHU^m \De\fHU}/\ampHU^m$ and the residuals $\varepsilon_\circ^r {\gets} \norm{r_\circ}$.}

\STATE{Increase the counter of the alternating direction iterations $q \gets q+1$.}

\ENDWHILE

\STATE{Increase the mode counter $m \gets m+1$.}
\ENDWHILE
\end{algorithmic}
\end{algorithm}

In the greedy enrichment loop, first a predictor of the spatial mode is computed as the solution of the HDG global and local problems using a guess for the parametric mode (Algorithm \ref{alg:HDG-PGD} - Steps 3-5). Then, the alternating direction scheme computes the corrections of the parametric (Algorithm \ref{alg:HDG-PGD} - Steps 8-10) and spatial mode (Algorithm \ref{alg:HDG-PGD} - Steps 11-14) solving a parametric linear system and the HDG global and local problems, respectively. The nonlinear iterations of the alternating direction scheme stop when the amplitude $\ampHU^m \De \fHU$ of the correction is negligible with respect to the amplitude $\ampHU^m$ of the current mode and the residuals of the spatial and parametric problems are below a given tolerance (Algorithm \ref{alg:HDG-PGD} - Steps 7 and 15). The stopping criterion for the greedy enrichment algorithm relies on the relative amplitude $\ampHU^m$ of the current mode being negligible with respect to the first mode $\ampHU^1$ (Algorithm \ref{alg:HDG-PGD} - Step 2). Alternative stopping criteria based on normalising the amplitude of the current mode with respect to the cumulative amplitudes of the previous modes have also been considered in the literature, see e.g.~\cite{tsiolakis2020nonintrusive}. Note that for the purpose of normalisation (Algorithm \ref{alg:HDG-PGD} - Step 14), an appropriate norm needs to be defined and the $\elinf$ norm has been utilised for the simulations in section~\ref{sc:examples}.

\subsubsection{Discretisation of the spatial and parametric problems} \label{sc:discretisation}

The discretisation of the local problems of the spatial iteration using an isoparametric formulation with equal interpolation for all the variables~\cite{RS-SH:16,HDG-NEFEM,RS-SGKH:18}, leads to a system of equations for each element with the following structure:
\begin{equation} \label{eq:localProblemSystem}
\begin{bmatrix}
\mat{A}_{LL} & \mat{A}_{Lu} & \bm{0} & \bm{0} \\
\mat{A}_{Lu}^T & \mat{A}_{uu} & \mat{A}_{up} & \bm{0} \\
\bm{0} & \mat{A}_{up}^T & \bm{0} & \bm{a}_{\rho p}^T \\
\bm{0} & \bm{0} & \bm{a}_{\rho p} & 0
\end{bmatrix}_{\!\! e}
\begin{Bmatrix}
\vect{F}_{\!\! L} \\
\vect{F}_{\!\! u} \\
\vect{F}_{\!\! p} \\
\text{F}_{\!\! \zeta}
\end{Bmatrix} _{\!\! e}
= 
\begin{Bmatrix}
\vect{f}_L \\
\vect{f}_u \\
\vect{f}_p \\
0
\end{Bmatrix}_{\!\! e}
+
\begin{bmatrix}
\mat{A}_{L\hat{u}} \\
\mat{A}_{u\hat{u}} \\
\mat{A}_{p\hat{u}} \\
\bm{0}
\end{bmatrix}_{\!\! e}
\vect{F}_{\!\! \hu}
+
\begin{Bmatrix}
\bm{0} \\
\bm{0} \\
\bm{0} \\
1
\end{Bmatrix}_{\!\! e}
\text{F}_{\!\! \rho}
,
\end{equation}
where $\vect{F}_{\!\! L}$, $\vect{F}_{\!\! u}$, $\vect{F}_{\!\! p}$ and $\vect{F}_{\!\! \hu}$ denote the nodal values of the unknown spatial functions $\ampL^m \De \fL$, $\ampU^m \De \fU$, $\ampP^m \De \fP$ and $\ampHU^m \De \fHU$ respectively and the constraint on the mean value $\text{F}_{\!\! \rho}$ of the pressure on the element boundaries is enforced using the Lagrange multiplier $\text{F}_{\!\! \zeta}$.

The only difference between the local system obtained in the spatial iteration of the proposed HDG-PGD approach and the local system of a standard HDG method~\cite{HDG-NEFEM,RS-SGKH:18} lies in the construction of the blocks forming the matrices $\mat{A}_{\odot \circledcirc}$ and vectors $\vect{f}_{\odot}$. As an example, let us consider the matrix $\mat{A}_{LL}$. In the proposed HDG-PGD framework, this matrix is defined as
\begin{equation} \label{eq:ALL}
\left( \mat{A}_{LL} \right)_{IJ} = - \sum_{k=1}^{\ndet} \beta_{\theta}^k \intE{N_I}{\nu^{-1} D^k N_J} \InsdTwo
\end{equation}
whereas in a standard HDG approach, the corresponding matrix is defined as
\begin{equation} \label{eq:ALL-HDG}
\left( \mat{A}_{LL} \right)_{IJ} = -\intE{N_I}{\nu^{-1}  N_J} \InsdTwo.
\end{equation}
In the above expressions $\{N_I\}$ denotes the set of shape functions used to define the spatial approximation of the mixed variable.

Similarly, the discretisation of the global problem of the spatial iteration leads to a system of equations for the trace of the velocity on the element boundaries and the mean value of the pressure in each element, namely
\begin{equation}\label{eq:globalProblemSystem}
\begin{aligned}
\sum_{e=1}^{\numel}\Big\{
\begin{bmatrix} 
\mat{A}_{\hu L} & \mat{A}_{\hu u} & \mat{A}_{\hu p} 
\end{bmatrix}_{\! e}
\begin{Bmatrix} 
\vect{F}_{\!\! L} \\
\vect{F}_{\!\! u} \\
\vect{F}_{\!\! p} \\
\end{Bmatrix}_{\!\! e}
+
[\mat{A}_{\hu\hu}]_{e} \, \vect{F}_{\!\! \hu} \Big\}
&= 
\sum_{e=1}^{\numel} [\vect{f}_{\hu}]_{e}, 
\\
\vect{1}^T \, [\mat{A}_{p\hu}]_{e} \vect{F}_{\!\! \hu} &= - \vect{1}^T \, [\vect{f}_p]_{e}.
\end{aligned}
\end{equation}

As usual in an HDG context, the local problem of equation~\eqref{eq:localProblemSystem} is used to express the spatial part of the gradient of the velocity, the velocity and the pressure in terms of the spatial part of the trace of the velocity and the mean pressure. Introducing these expressions into the global problem, leads to the global system
\begin{equation}\label{eq:globalProblemSystemFinal}
\begin{bmatrix}
\widehat{\mat{K}} & \mat{G} \\
\mat{G}^T         & \mat{0}  
\end{bmatrix}
\begin{Bmatrix}
\vect{F}_{\!\! \hu}\\ 
\vect{F}_{\!\! \rho} 
\end{Bmatrix} 
=
\begin{Bmatrix}
\vect{\hat{f}}_{\hu} \\ 
\vect{\hat{f}}_{\rho} 
\end{Bmatrix},
\end{equation}
where the only unknowns are the spatial parts of the trace of the velocity and the mean pressure.

In a similar fashion, the discretisation of the parametric problem \eqref{eq:WeakPGDParam} using Lagrange shape functions leads to an algebraic system of equations whose unknowns are the nodal values of the parametric modes.

\subsubsection{A remark for a computationally efficient implementation} \label{sc:computational}

The evaluation of the right hand sides of the PGD spatial and parametric iterations tends to become computationally expensive when approximations with a large number of modes are considered. Indeed, the number of terms involved in such computation experiences a geometric growth rate during the iterations of the greedy algorithm. 

In order to ease the computational burden of the overall algorithm, the number of terms in the modal approximations $\bupgd^m, \ppgd^m, \Lpgd^m, \bhupgd^m \text{ and } \rpgd^m$ is reduced. It is well known that the terms in the PGD reduced basis are not orthogonal to each other and repeated information may appear. Hence, orthogonal separable approximations featuring $\widetilde{m} < m$ modes are constructed via the PGD compression~\cite{DM-MZH:15,PD-DZGH-20}, that is, a least-squares higher-order projection minimising the $\eltwo$ norm of the difference between target and test functions, namely
\begin{equation*}
	\begin{aligned}
	\Lpgd^{\widetilde{m}} &= \argmin_{\bW \in \Wh} \| \bW - \Lpgd^m \|_{\eltwo(\Omega \times \bI)} , \\
	\bupgd^{\widetilde{m}} &= \argmin_{\bv \in \Vh} \| \bv - \bupgd^m \|_{\eltwo(\Omega \times \bI)} , \\
	\ppgd^{\widetilde{m}} &= \argmin_{v \in \sVh} \| v - \ppgd^m \|_{\eltwo(\Omega \times \bI)} , \\
	\bhupgd^{\widetilde{m}} &= \argmin_{\bhv \in \hVh} \| \bhv - \bhupgd^m \|_{\eltwo(\Gamma \cup \Gamma_N \cup \Gamma_S \times \bI)} , \\
	\rpgd^{\widetilde{m}} &= \argmin_{q \in \mathbb{R}^{\numel} \otimes \Lh(\bI)} \| q - \rpgd^m \|_{\eltwo(\mathbb{R}^{\numel} \times \bI)} .
	\end{aligned}
\end{equation*}

From a practical point of view, the PGD compression is applied during the enrichment strategy described in algorithm \ref{alg:HDG-PGD}. A trade-off between the cost of performing the greedy iterations with a larger number of modes and the extra cost required by the PGD compression needs to be achieved. For the simulations in section~\ref{sc:examples}, PGD compression is applied every ten new computed modes for the analytical examples and every five for the microfluidics test cases.

\section{Numerical examples} \label{sc:examples}

This section presents four numerical examples. The first two examples are used to validate the implementation of the proposed approach as well as to study a number of properties of the proposed ROM. The last two examples consider two applications taken from the biomechanics community and involve the Stokes flow around a micro-swimmer formed by two spheres and the flow around a sphere in a corrugated channel. All the examples consider geometric parameters as extra coordinates within the proposed PGD approach.

\subsection{Coaxial Couette flow} \label{sc:Couette}

The first example considers the well known coaxial Couette flow problem~\cite{childs2010rotating}, consisting of the flow confined within two infinite coaxial circular cylinders with radius $\Rin$ and $\Rout$ respectively, with $\Rin < \Rout$. The boundary conditions introduce the known angular velocities, $\Vin$ and $\Vout$, at $\Rin$ and $\Rout$, respectively. The problem has analytical solution, given by the 
azimuthal component of the velocity as
\begin{equation} \label{eq:couetteAnalytical}
v_{\phi} = \frac{ \Rout^2 \Vout -\Rin^2 \Vin }{ \Rout^2 -\Rin^2 } r + \frac{ ( \Vin - \Vout) \Rout^2 \Rin^2 }{ \Rout^2 -\Rin^2 } \frac{1}{r}
\end{equation}
where $r$ is the distance to the axis of the cylinders.

To demonstrate the applicability of the proposed ROM the problem is considered in two dimensions, with $\Omega^{\mu} = \{\bx \in \mathbb{R}^2 \; | \; \mu_1 \leq r^{\mu} \leq \Rout \}$, with $\Rout=5$ and $\mu_1 \in \bI = \I^1 = [1,3]$ and where $r^{\mu} = \sqrt{(x_1^{\mu})^2 + (x_2^{\mu})^2}$. The reference domain is chosen to be $\Omega = \{ \bX \in \mathbb{R}^2 \; | \; 1 \leq r \leq \Rout \}$ and the mapping between the reference and the geometrically parametrised domains is defined by the general separable expression of equation~\eqref{eq:mappingDomain} with the mapping of equation~\eqref{eq:displacementSep} given by
\begin{equation} \label{eq:couetteMappingD}
\begin{aligned}
\bM^1(\bX) & = \frac{1}{r} \bX
\qquad &
\psi^1(\mu) & = \dfrac{ \Rout (\mu - 1)}{\Rout - 1}, 
\\
\bM^2(\bX) & = \bX
\qquad &
\psi^2(\mu) & = \dfrac{\Rout - \mu}{\Rout - 1},
\end{aligned}
\end{equation}
where $r = \sqrt{x_1^2 + x_2^2}$.
The Jacobian of the mapping is also written in the general separated form of equation~\eqref{eq:JacoSep}, with
\begin{equation}\label{eq:couetteMappingJ}
\bJ^1(\bX) = \frac{1}{r^3}
\begin{bmatrix}
x_2^2 & -x_1 x_2 \\
-x_1 x_2 & x_1^2
\end{bmatrix},
\qquad
\bJ^2(\bX) = \mat{I}_2.
\end{equation}

\begin{figure}[!tb]
\centering
\subfigure[Mesh 1]{\includegraphics[width=0.24\textwidth]{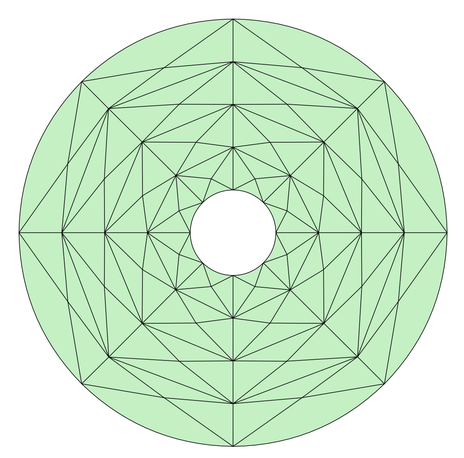}}
\subfigure[Mesh 2]{\includegraphics[width=0.24\textwidth]{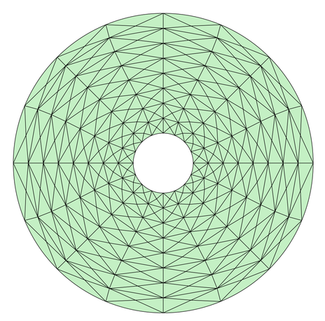}}
\subfigure[Mesh 3]{\includegraphics[width=0.24\textwidth]{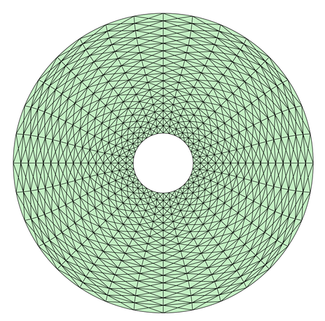}}
\subfigure[Mesh 4]{\includegraphics[width=0.24\textwidth]{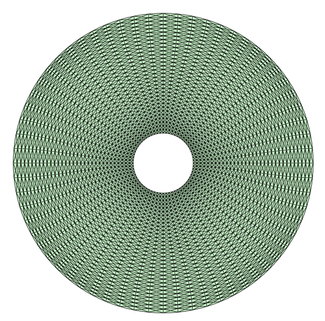}}
\caption{Coaxial Couette flow: Four triangular meshes of the reference domain.}
\label{fig:couetteRefMesh}
\end{figure}

For the numerical experiments in this section, four triangular meshes of the reference domain are generated, as shown in Figure~\ref{fig:couetteRefMesh}.
The meshes have 128, 512, 2,048 and 8,192 elements respectively. 

\begin{figure}[!tb]
	\centering
	\subfigure[$m=1$]{\includegraphics[width=0.24\textwidth]{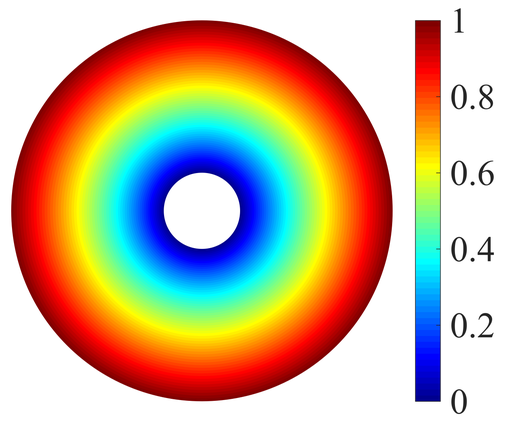}}
	\subfigure[$m=2$]{\includegraphics[width=0.24\textwidth]{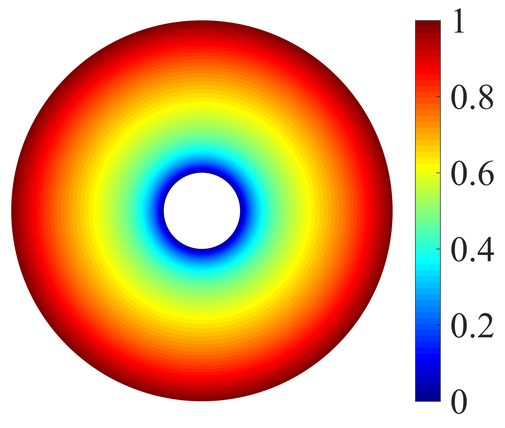}}
	\subfigure[$m=3$]{\includegraphics[width=0.24\textwidth]{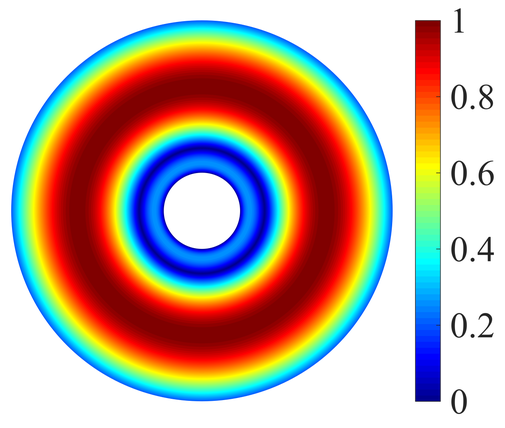}}
	\subfigure[$m=4$]{\includegraphics[width=0.24\textwidth]{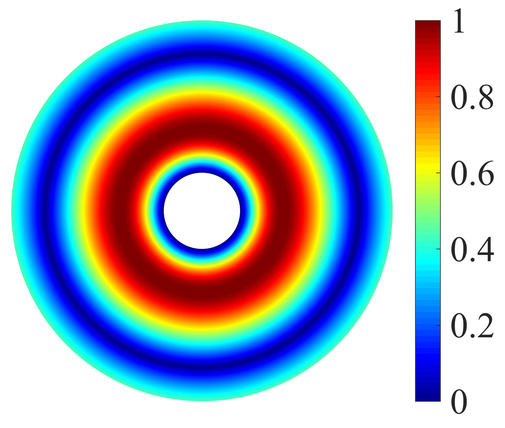}}
	\caption{Coaxial Couette flow: First four normalised spatial modes of the velocity field.}
	\label{fig:couetteVModes}
\end{figure}

The proposed HDG-PGD framework is used to obtain the generalised solution of the parametric Stokes problem. The first four normalised modes of the magnitude of the velocity field are displayed in figure~\ref{fig:couetteVModes}. 
The computation was performed using the second mesh shown in figure~\ref{fig:couetteRefMesh} with a degree of approximation $k=4$ for all the variables and with a mesh of 1,000 elements in the parametric dimension with also $k=4$. As usual in a the context of ROMs, the first modes capture the most relevant and global features of the solution whereas the features captured for the next modes only introduce localised features.

\begin{figure}[!tb]
	\centering
	\includegraphics[width=0.49\textwidth]{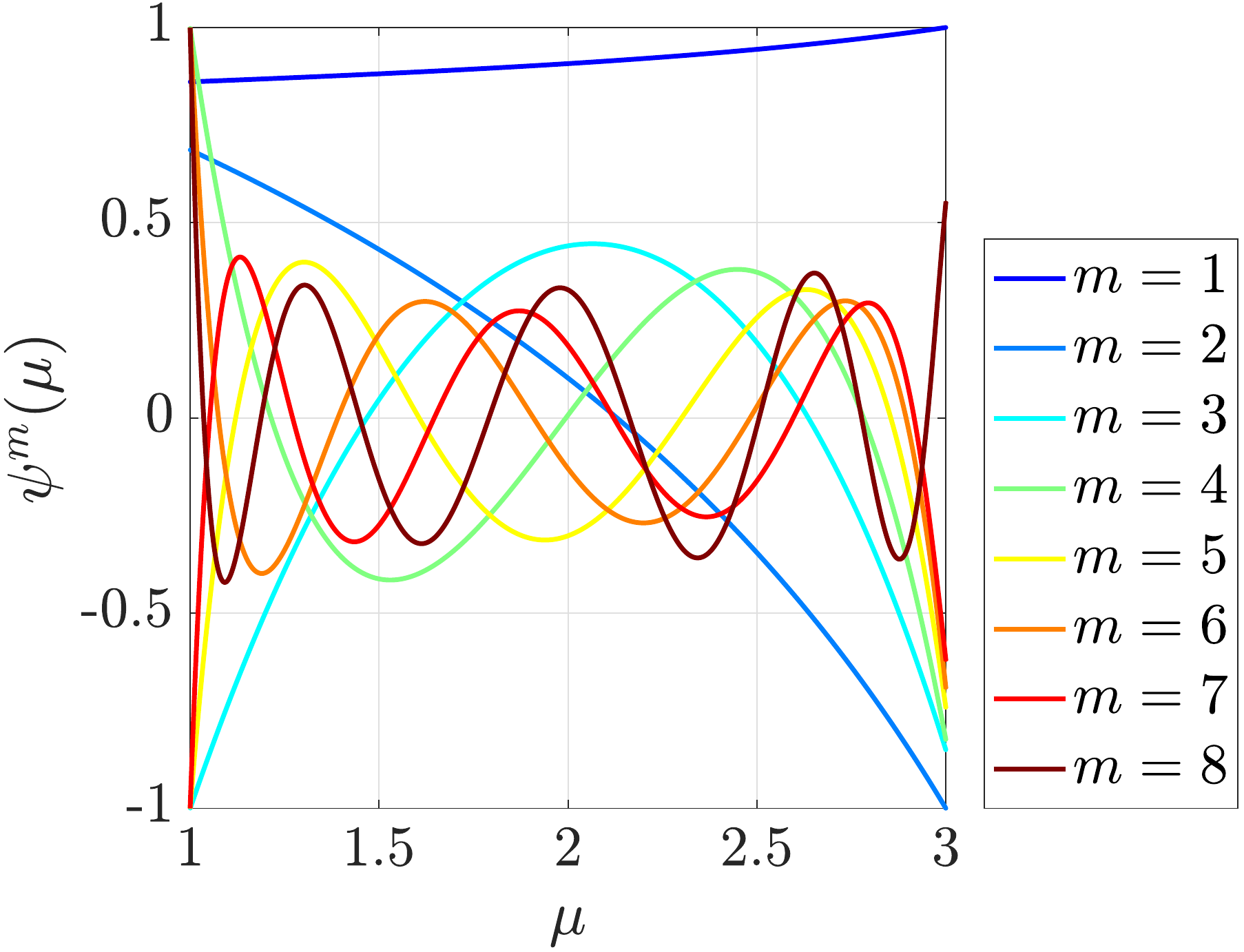}
	\caption{Coaxial Couette flow: First eight normalised parametric modes.}
	\label{fig:couetteParamModes}
\end{figure}

Figure~\ref{fig:couetteParamModes} shows the first eight normalised parametric modes computed.
It can be observed that the first three modes are smooth, whereas the next modes, that have a less relevant contribution to the generalised solution, show a more oscillatory character. 

\begin{figure}[!tb]
	\centering
	\includegraphics[width=0.45\textwidth]{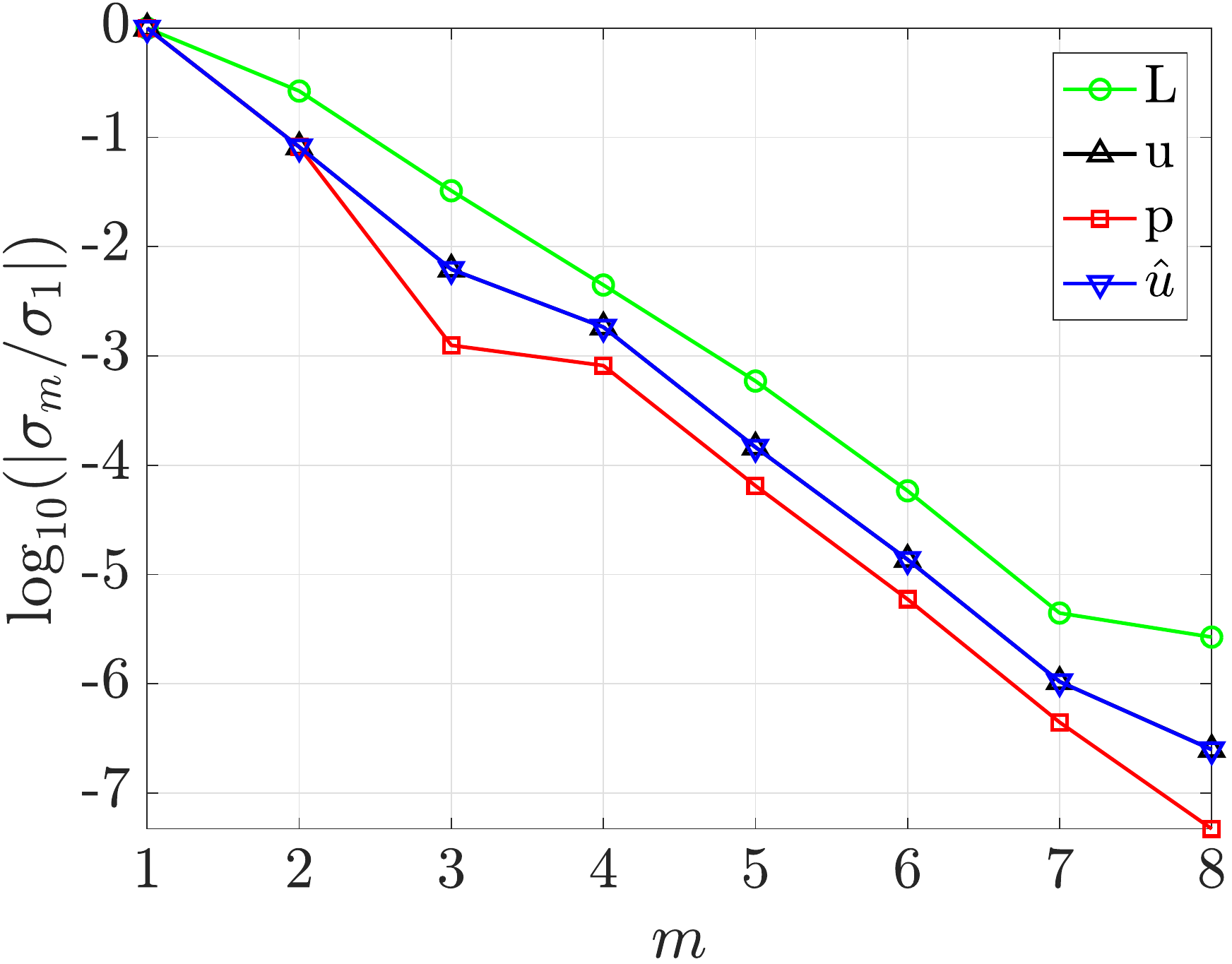}
	\caption{Coaxial Couette flow: Convergence of the mode amplitudes.}
	\label{fig:couetteAmplitudes}
\end{figure}

To quantify the importance of the modes on the generalised solution, figure~\ref{fig:couetteAmplitudes} shows the relative amplitudes of the modes with respect to the amplitude of the first mode for all the variables.
It can be clearly observed that the fourth mode has an amplitude that is already more than 100 times smaller than the amplitude of the first mode. After computing only nine modes the relative amplitude is already of the order of $10^{-6}$. It is worth noting that in practice it is not required to add modes with such a lower relative amplitude with respect to the first mode, but in this first example nine modes are computed to show the rapid decrease in their amplitudes.

\begin{figure}[!tb]
	\centering
	\subfigure[$\mu_1 = 1$, $m=1$] {\includegraphics[width=0.24\textwidth]{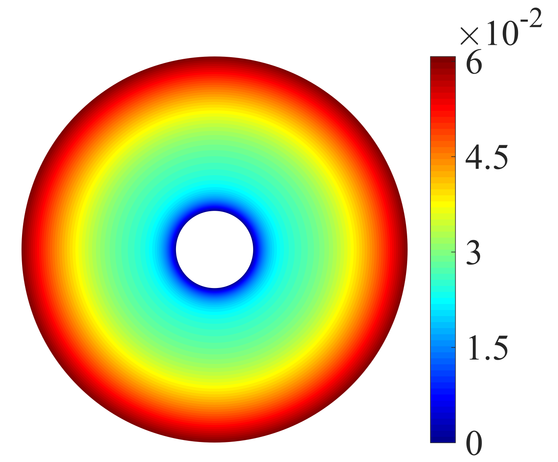}}
	\subfigure[$\mu_1 = 1$, $m=2$] {\includegraphics[width=0.24\textwidth]{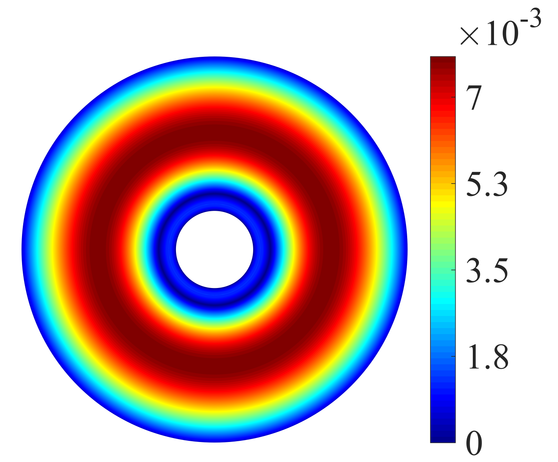}}
	\subfigure[$\mu_1 = 1$, $m=3$] {\includegraphics[width=0.24\textwidth]{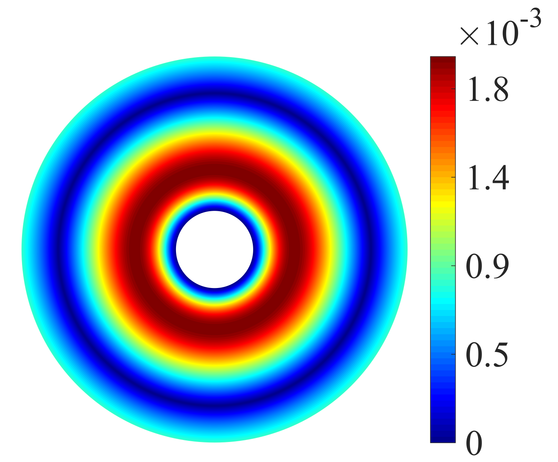}}
	\subfigure[$\mu_1 = 1$, $m=4$] {\includegraphics[width=0.24\textwidth]{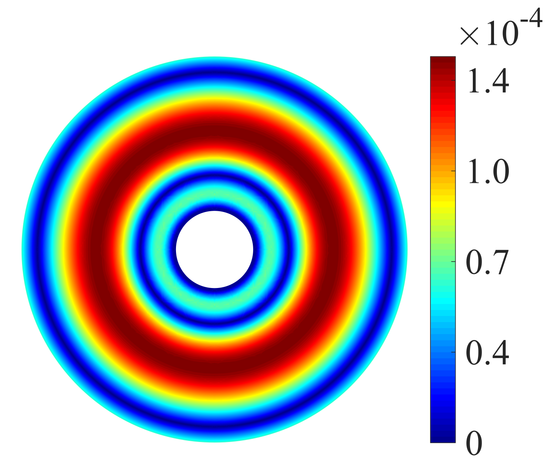}}
	\subfigure[$\mu_1 = 2$, $m=1$] {\includegraphics[width=0.24\textwidth]{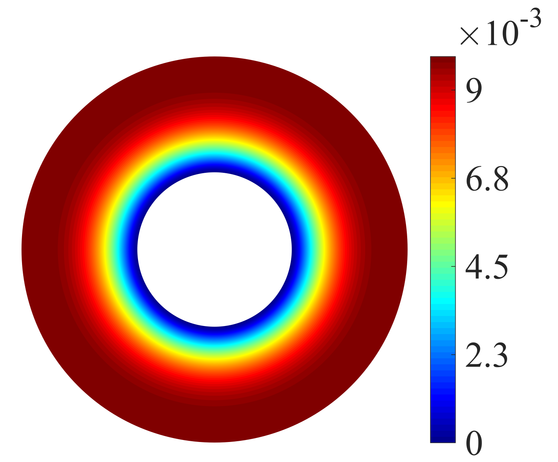}}
	\subfigure[$\mu_1 = 2$, $m=2$] {\includegraphics[width=0.24\textwidth]{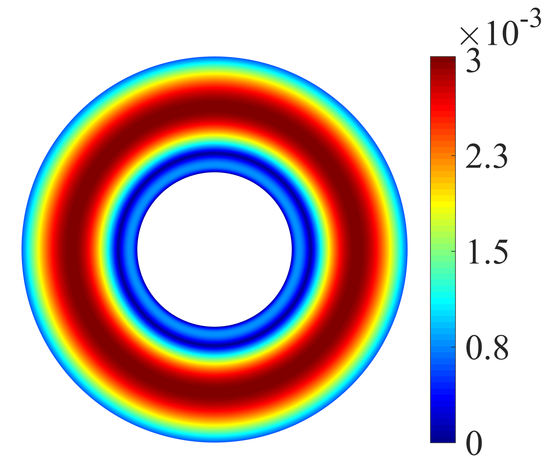}}
	\subfigure[$\mu_1 = 2$, $m=3$] {\includegraphics[width=0.24\textwidth]{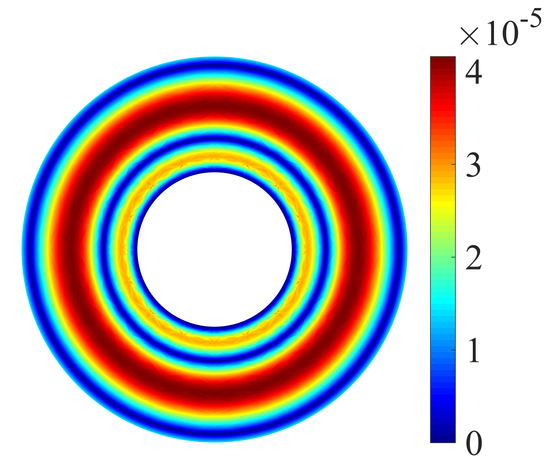}}
	\subfigure[$\mu_1 = 2$, $m=4$] {\includegraphics[width=0.24\textwidth]{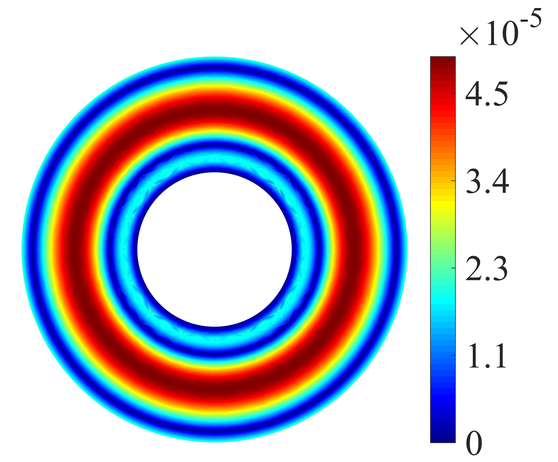}}
	\subfigure[$\mu_1 = 3$, $m=1$] {\includegraphics[width=0.24\textwidth]{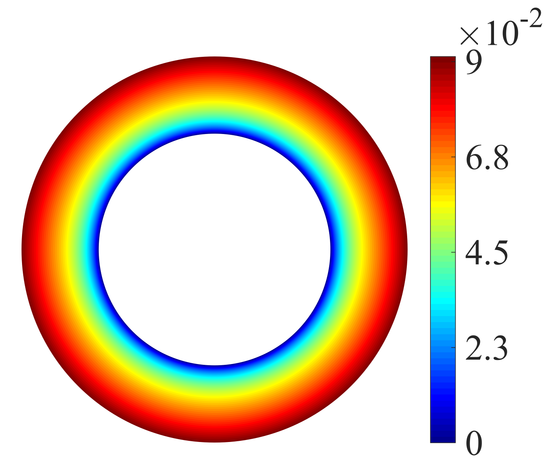}}
	\subfigure[$\mu_1 = 3$, $m=2$] {\includegraphics[width=0.24\textwidth]{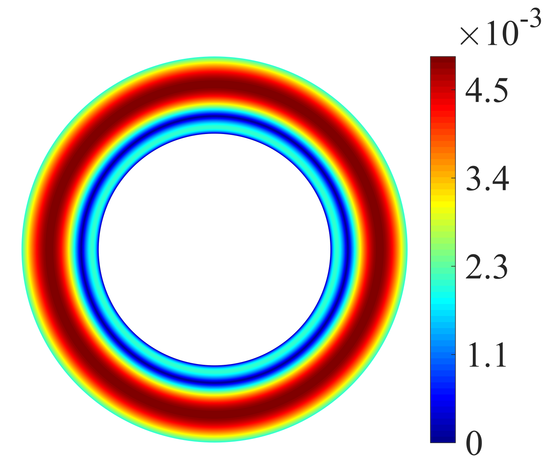}}
	\subfigure[$\mu_1 = 3$, $m=3$] {\includegraphics[width=0.24\textwidth]{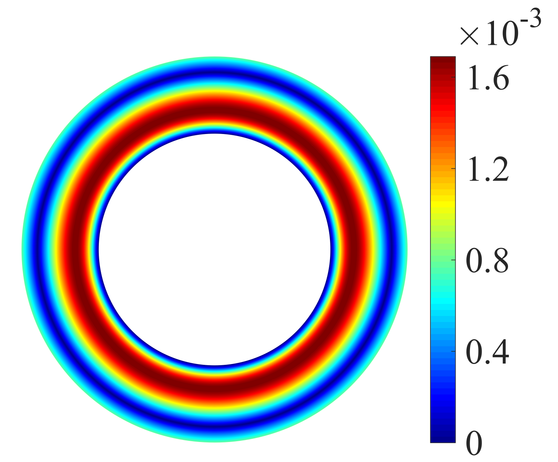}}
	\subfigure[$\mu_1 = 3$, $m=4$] {\includegraphics[width=0.24\textwidth]{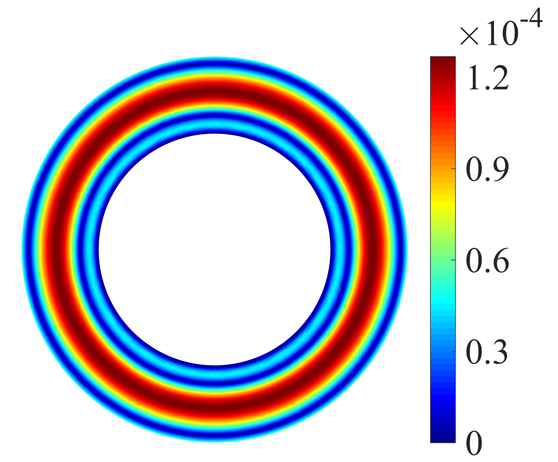}}
	\caption{Coaxial Couette flow: Absolute value of the error of the velocity magnitude using $n$ PGD modes and for different values of the geometric parameter $\mu_1$. A quartic approximation is used for all variables in the second mesh of figure~\ref{fig:couetteRefMesh}.}
	\label{fig:couetteModeErr}
\end{figure}

Once the generalised solution is computed, it is of interest to quantify its accuracy. Figure~\ref{fig:couetteModeErr} shows the absolute value of the error of the velocity magnitude using as the number of modes is increased for three relevant configurations corresponding to the parameter $\mu_1 = 1$, $\mu_1 = 2$ and $\mu_1 = 3$.
The results show that with only one PGD mode an absolute error below $10^{-1}$ is already obtained for all three configurations, with more accurate results for the case with $\mu_1 = 2$. With two PGD modes the error drops substantially, being less than $7 \times 10^{-3}$ in all cases, and with only PGD modes the error is below $2 \times 10^{-4}$ for the three configurations considered.

\begin{figure}[!tb]
	\centering
	\subfigure[$\bL$] {\includegraphics[width=0.49\textwidth]{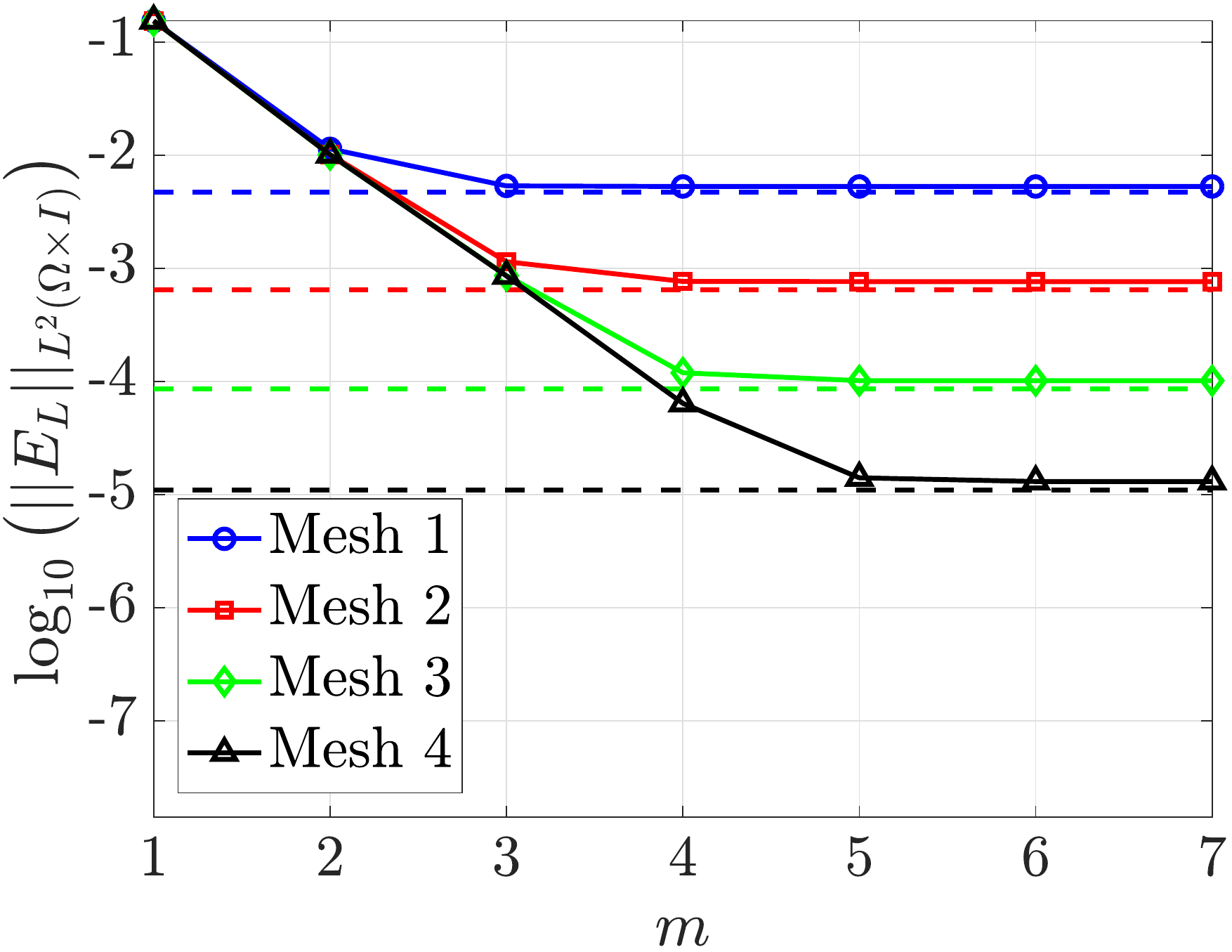}}
	\subfigure[$\bu$] {\includegraphics[width=0.49\textwidth]{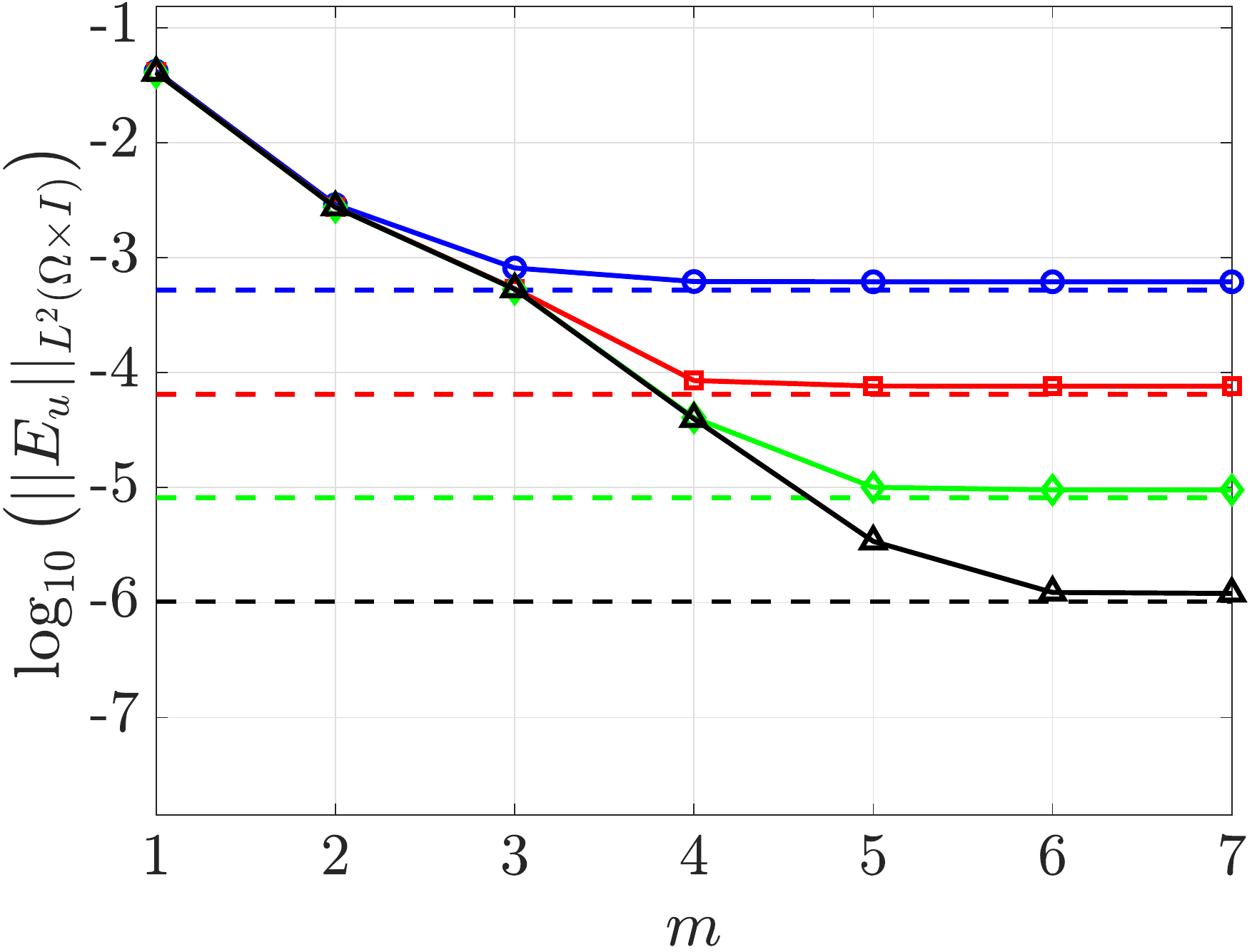}}
	\subfigure[$p$]     {\includegraphics[width=0.49\textwidth]{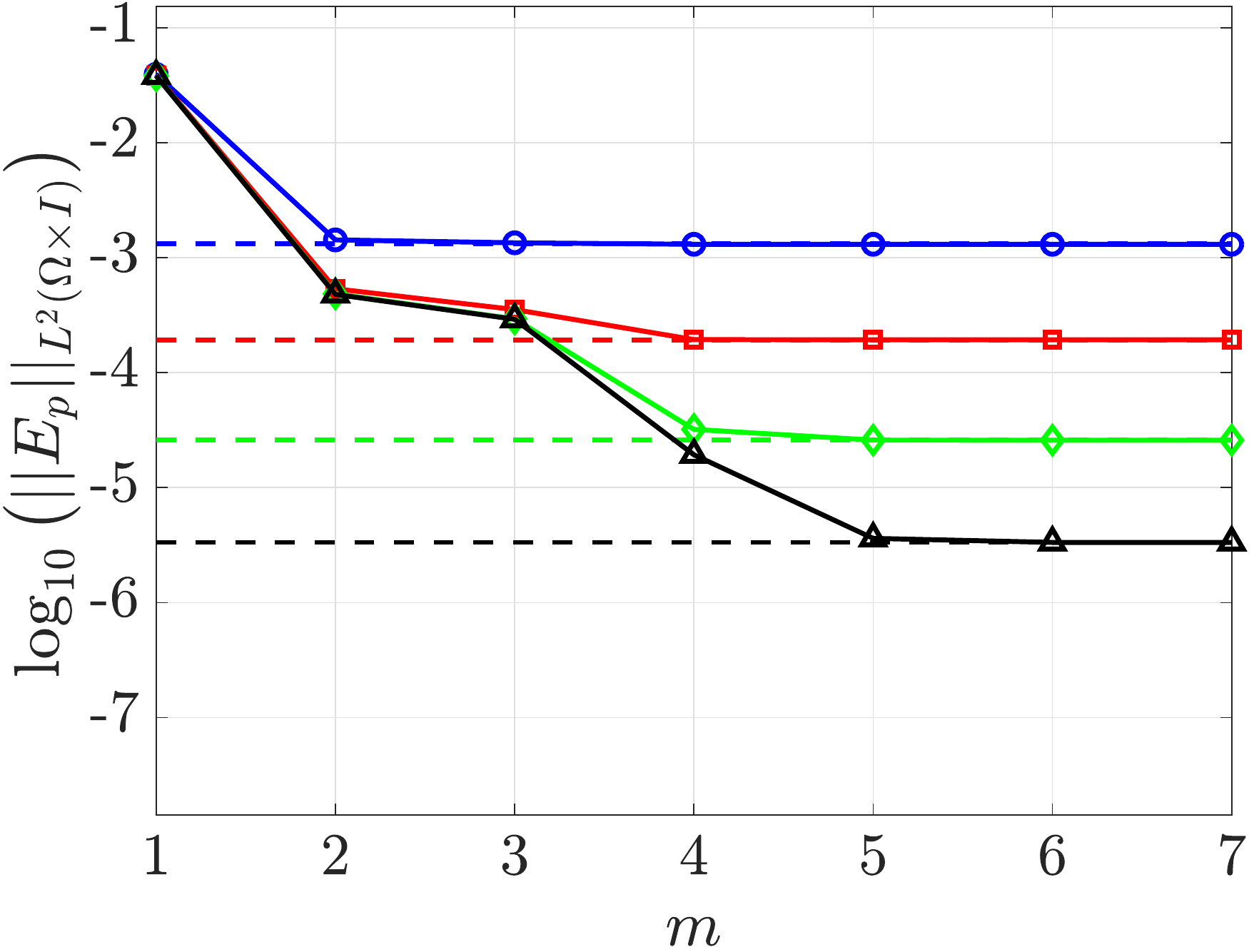}}
	\subfigure[$\bhu$]{\includegraphics[width=0.49\textwidth]{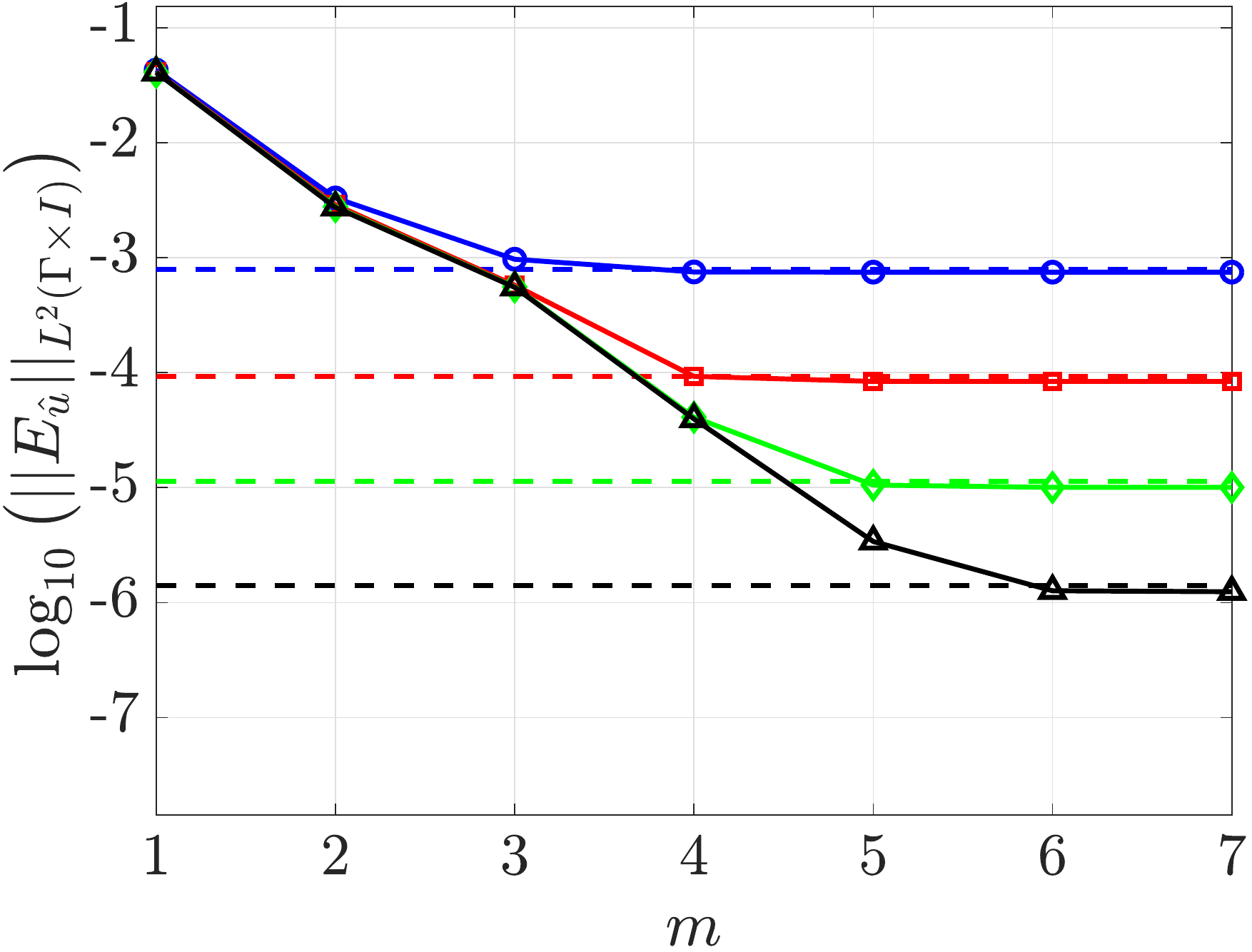}}
	\caption{Coaxial Couette flow: convergence of the $\eltwo$ norm of the error for $\bL$, $\bu$, $p$ and $\bhu$ as the number of PGD modes is increased. A quadratic approximation is used for all the variables.}
	\label{fig:couetteModeConv}
\end{figure}

To further illustrate the accuracy of the proposed HDG-PGD approach, the relative error in the $\eltwo(\Omega \times \bI)$ norm, defined as
\begin{equation} \label{eq:PGDerrorAnalytical}
\varepsilon_{\texttt{PGD}} 
= \left( \frac{ \displaystyle \int_{\I_1} \int_{\Omega} ( \bm{u}_{\texttt{PGD}} - \bm{u} ) \cdot ( \bm{u}_{\texttt{PGD}} - \bm{u} ) d\Omega \,  d\mu }
{ \displaystyle \int_{\I_1} \int_{\Omega}  \bm{u} \cdot \bm{u}\, d\Omega \, d\mu} \right)^{1/2} ,
\end{equation}
is studied and compared to the error of the full order HDG approach. Figure~\ref{fig:couetteModeConv} shows the evolution of $\varepsilon_{\texttt{PGD}}$, for all the variables, as the number of PGD modes is increased, for different meshes using a quadratic degree of approximation.
The discontinuous lines in Figure~\ref{fig:couetteModeConv} show the relative error of the full order HDG method, measured in the $\eltwo(\Omega \times \bI)$ norm. It is worth noting that the computation of the error for the full order approach requires the computation of a large number of solutions. More precisely, the number of HDG solutions required is equal to the number of elements in the parametric space multiplied by number of integrations points in each element. 

The results show that the error of the proposed ROM converges monotonically to the error of the full order approach with as the number of modes is increased. In all cases the number of PGD modes required to reach the maximum accuracy on a given mesh is lower than six. Furthermore, the results in figure~\ref{fig:couetteModeConv} illustrate the increased level of accuracy obtained as the spatial and parametric discretisations are refined. Analogous results, not reported here for brevity, are obtained for lower and higher orders of approximation.

\begin{figure}[!tb]
	\centering
	\subfigure[$\bL$] {\includegraphics[width=0.49\textwidth]{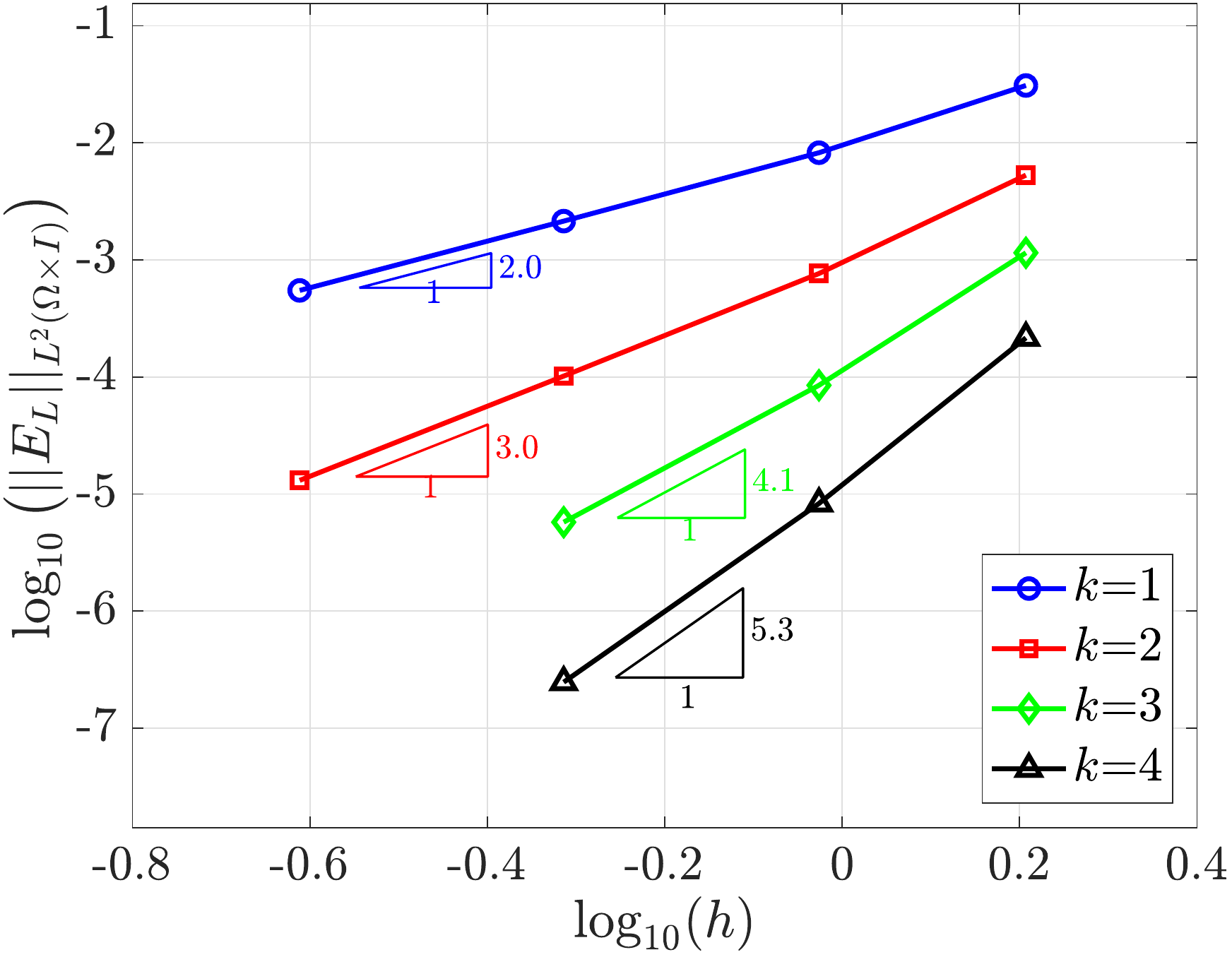}}
	\subfigure[$\bu$] {\includegraphics[width=0.49\textwidth]{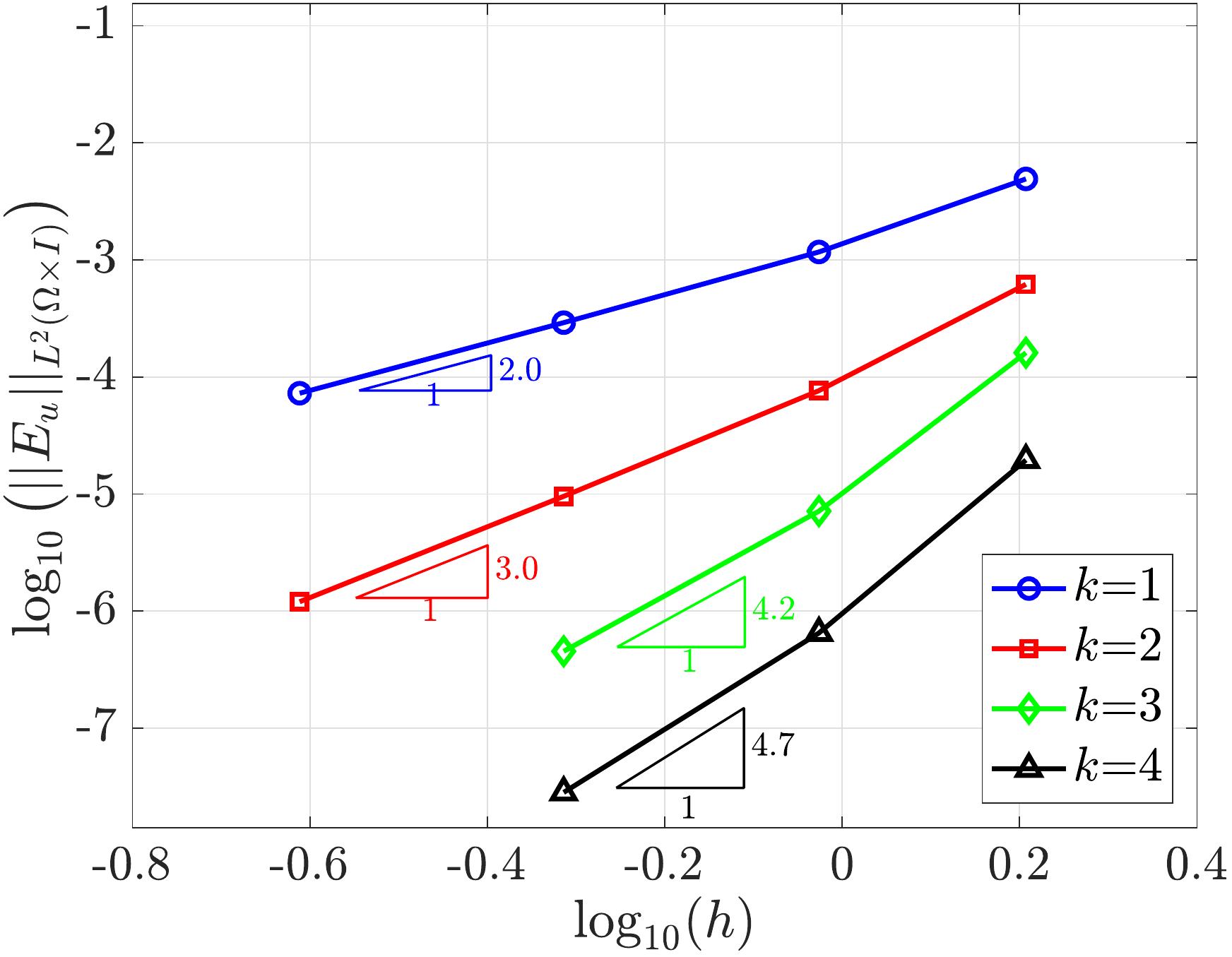}}
	\subfigure[$p$]   {\includegraphics[width=0.49\textwidth]{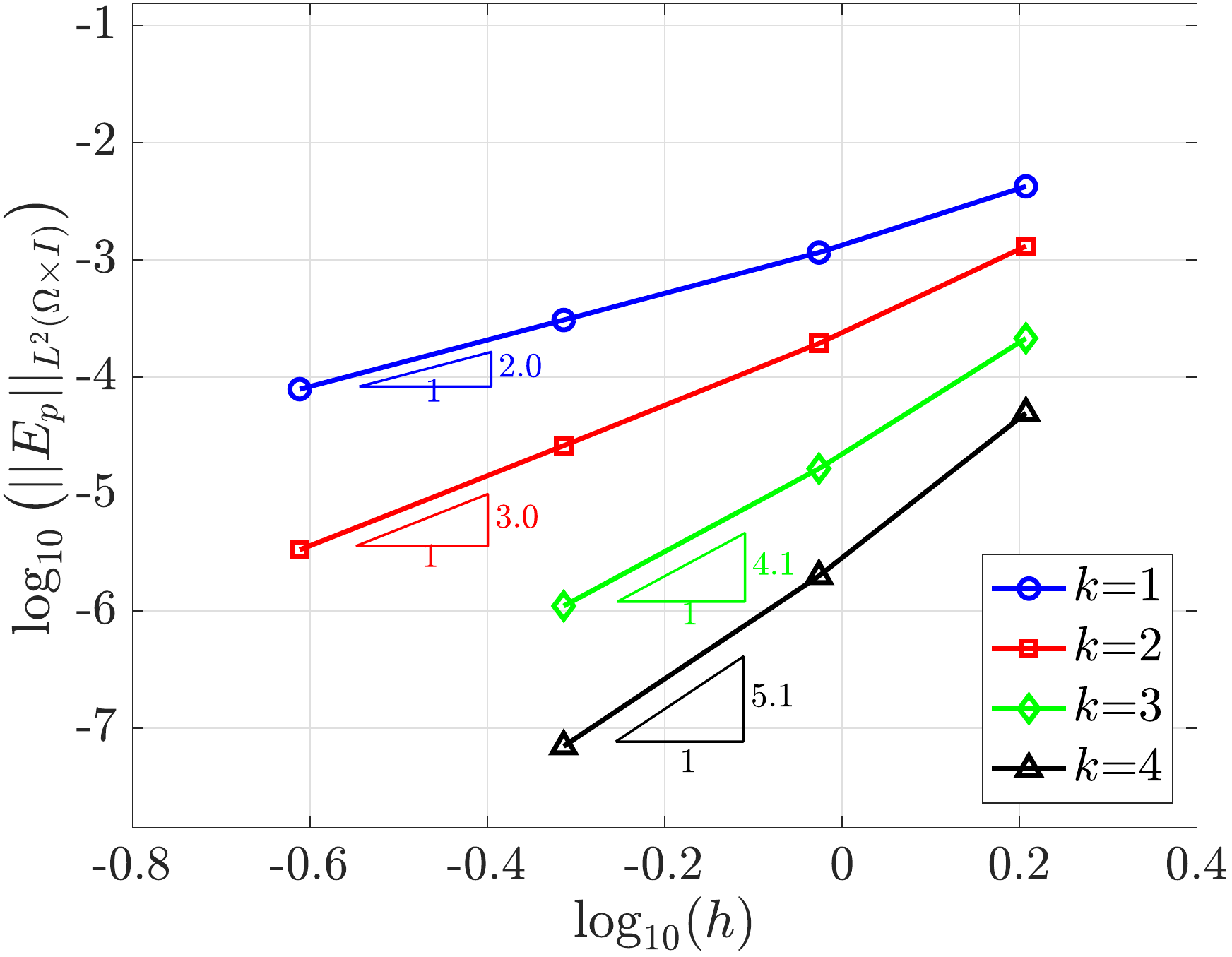}}
	\subfigure[$\bhu$]{\includegraphics[width=0.49\textwidth]{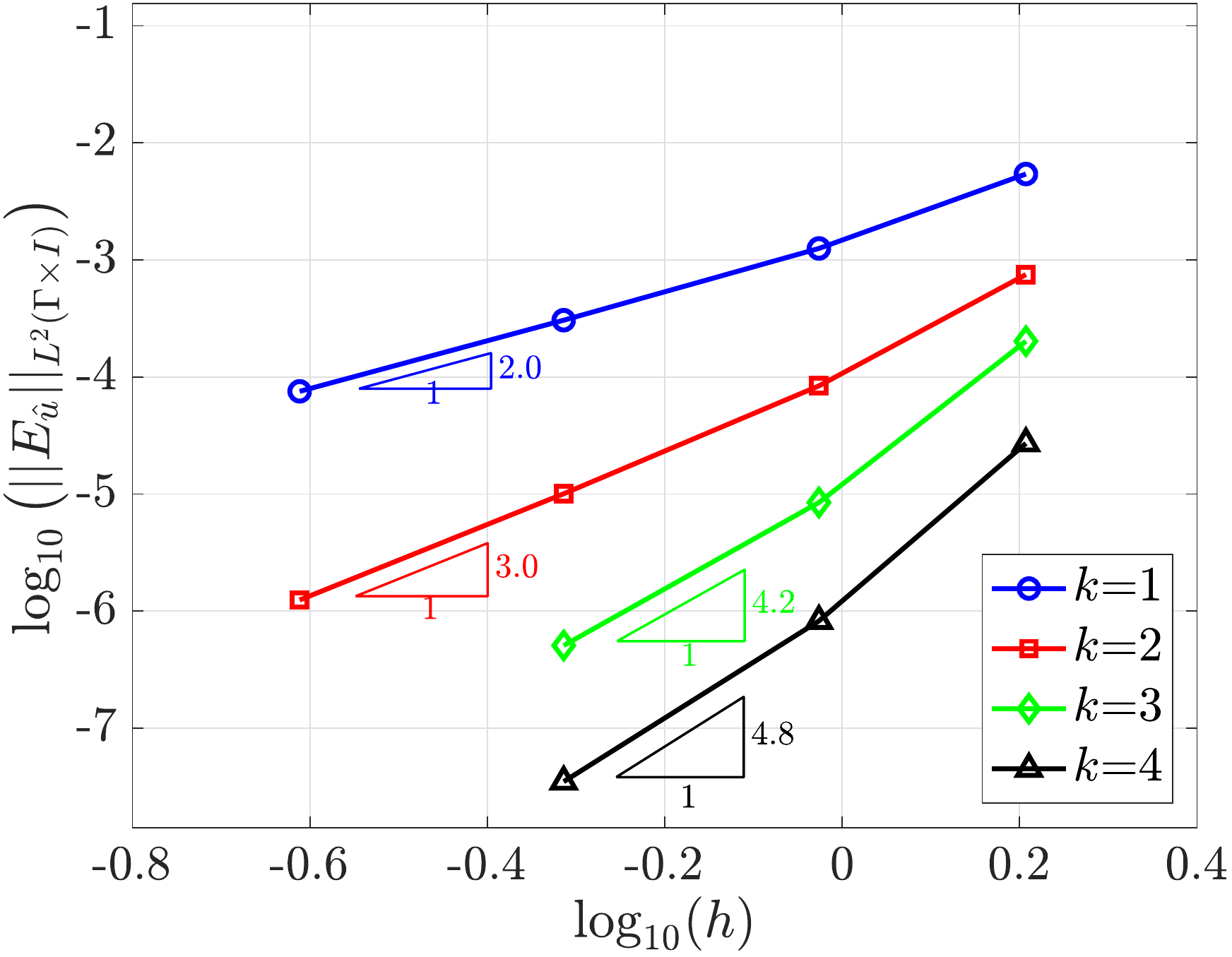}}
	\caption{Coaxial Couette flow: mesh convergence of the $\eltwo$ norm of the error for $\bL$, $\bu$, $p$ and $\bhu$.}
	\label{fig:couetteHConv}
\end{figure}

Next, the optimal approximation properties of the proposed HDG-PDG method are studied by performing a mesh convergent study. Figure~\ref{fig:couetteHConv} shows the evolution of the relative error in the $\eltwo(\Omega \times \bI)$ norm as a function of the characteristic element size, $h$, for different orders of approximation and for all the variables of the HDG formulation.
The optimal rate of convergence, equal to $h^{k+1}$, is approximately observed for all the variables. In each case, the minimum number of PGD modes required to achieve the accuracy of the full order method is selected, as previously discussed when presenting the results of figure~\ref{fig:couetteModeConv}. 

Finally, it is worth mentioning the differences between the proposed HDG-PGD approach presented here and the recently proposed PGD approach for geometrically parametrised domains in~\cite{sevilla2020solution} using standard finite elements for the spatial discretisation. First, the current approach does not require the higher order PGD projection to separate the inverse of the determinant of the Jacobian, given the first-order character of the problem solved with HDG. Second, the current approach enables the use of the same degree of approximation for velocity and pressure, contrary to the standard FE approach where specific choices are required to satisfy the LBB condition. In the context of geometrically parametrised domains with curved boundaries this implies that the current approach enables the use of isoparametric elements whereas super-parametric or sub-parametric elements are required in the FE context. Third, the weak imposition of the Dirichlet boundary conditions, as usually done in a DG context, facilitates the construction of the generalised solution without the need for specific choices for the modes that satisfy the Dirichlet boundary conditions, as required by approaches. Finally, the results in figure~\ref{fig:couetteHConv} can be compared to the results in~\cite{sevilla2020solution}.

\subsection{Axisymmetric  Stokes flow past a sphere} \label{sc:sphere}

The second example considers the Stokes flow past a sphere, a typical test case for axisymmetric Stokes flow solvers. The domain of interest is selected as the region confined by two concentric spheres with radius $\Rin$ and $\Rout$ respectively, with $\Rin < \Rout$. This problem also has analytical solution, given, in polar coordinates, by the following velocity and pressure fields
\begin{equation} \label{eq:sphereAnalytical}
\begin{aligned}
u_r &= \frac{v_\infty}{2r^3} \left(  2r^3 - 3\Rin r^2 + \Rin^3 \right) \cos\theta, \\
u_\theta &= -\frac{v_\infty}{4r^3} \left(  4r^3 - 3\Rin r^2 + \Rin^3 \right) \sin\theta, \\
p &= p_{\infty} - \dfrac{3}{2r^2} \nu  v_\infty \Rin \cos\theta,
\end{aligned}
\end{equation} 
where $v_\infty$ and $p_{\infty}$ are the magnitude of the velocity and the pressure of the undisturbed flow, far away from the obstacle. A typical quantity of interest in this problem is the drag force, whose exact value is given by $F_D = 6 \pi \nu v_\infty \Rin$

Similar to the previous example, the geometric parameter considered here is the radius of the inner sphere. The parametric domain considers the axial symmetry of the problem is defined as 
$\Omega^{\mu} = \{\bx \in \mathbb{R}^2 \; | \; x_2^{\mu}\geq 0 \text{ and } \mu_1 \leq r^{\mu} \leq \Rout \}$, with $\Rout=5$ and $\mu_1 \in \bI = \I^1 = [1,3]$. The reference domain is chosen to be $\Omega = \{ \bX \in \mathbb{R}^2 \; | \; x_2 \geq 0 \text{ and }  1 \leq r \leq \Rout \}$. The mapping between the reference and the geometrically parametrised domains is exactly the same mapping utilised in the previous example, given by the two terms in equation~\eqref{eq:couetteMappingD}.

A no-slip boundary condition is imposed on the inner sphere, a Dirichlet boundary condition corresponding to the exact solution on the outer boundary and axial symmetry is imposed on the  rest of the boundary. The axial symmetry is imposed by selecting $\alpha=\beta=0$ in the matrices $\bD^{\bmu}$ and $\bE^{\bmu}$ in equation~\eqref{eq:stokesStrong}. As mentioned earlier, in Remark~\ref{rk:GammaS}, the portion of the boundary where the axial symmetry is imposed depends on the geometric parameter, but the normal and tangent to the boundary are independent on the geometric changes. Therefore, the matrices $\bD$ and $\bE$ do not depend upon the geometric parameters.

\begin{figure}[!tb]
	\centering
	\subfigure[$m=1$]{\includegraphics[width=0.24\textwidth]{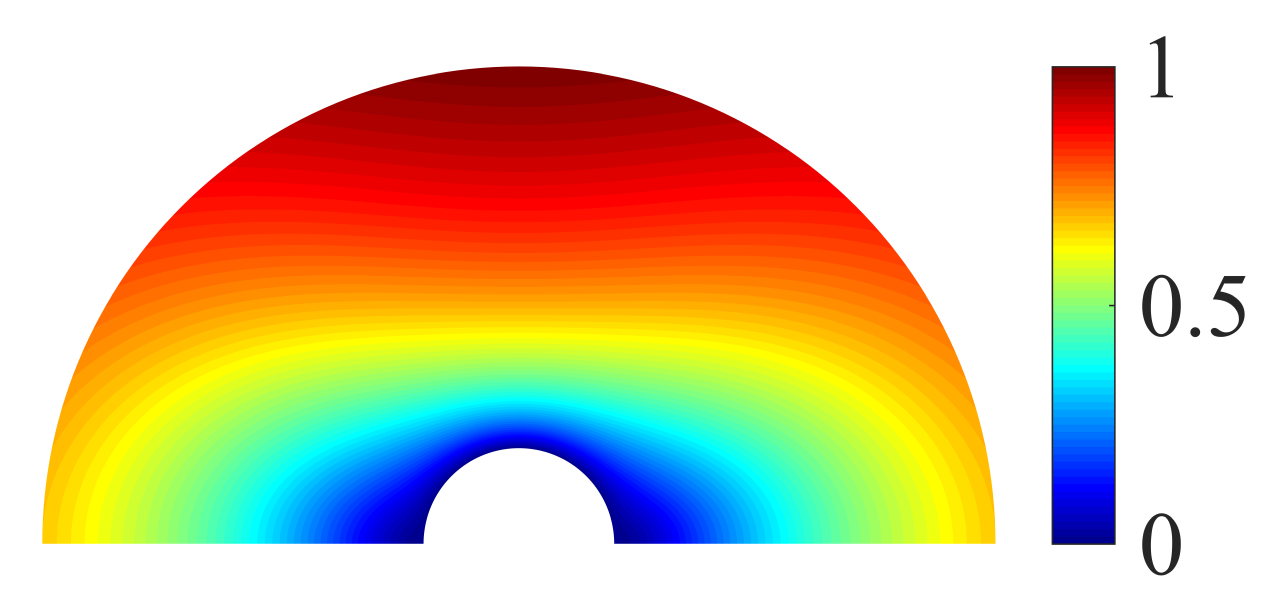}}
	\subfigure[$m=2$]{\includegraphics[width=0.24\textwidth]{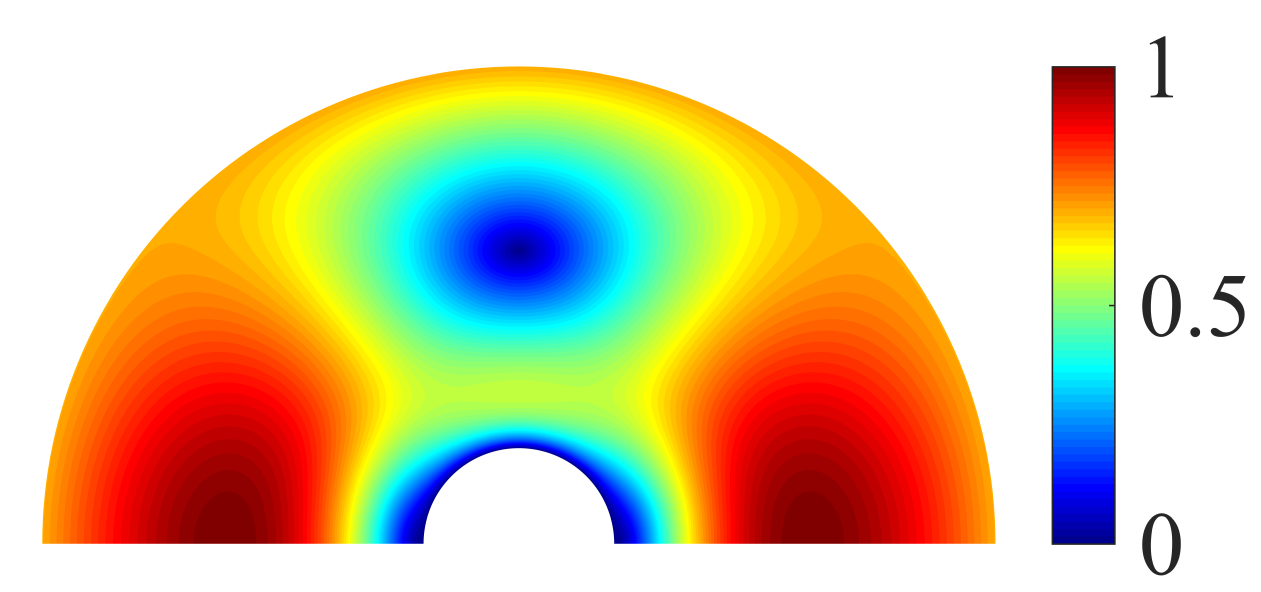}}
	\subfigure[$m=3$]{\includegraphics[width=0.24\textwidth]{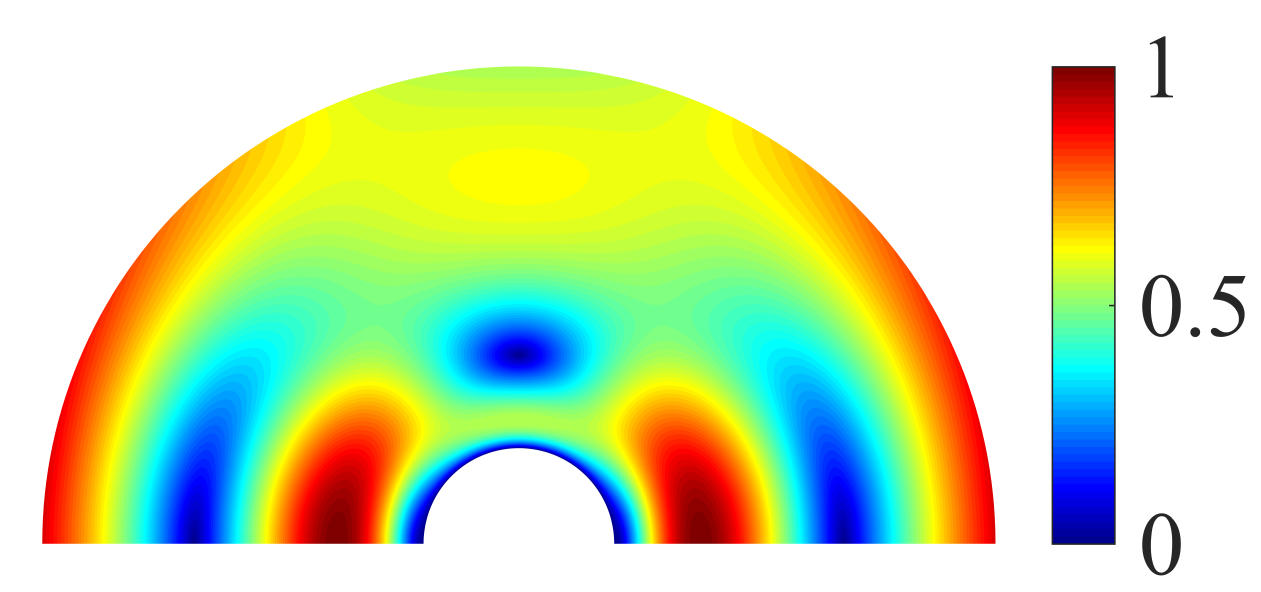}}
	\subfigure[$m=4$]{\includegraphics[width=0.24\textwidth]{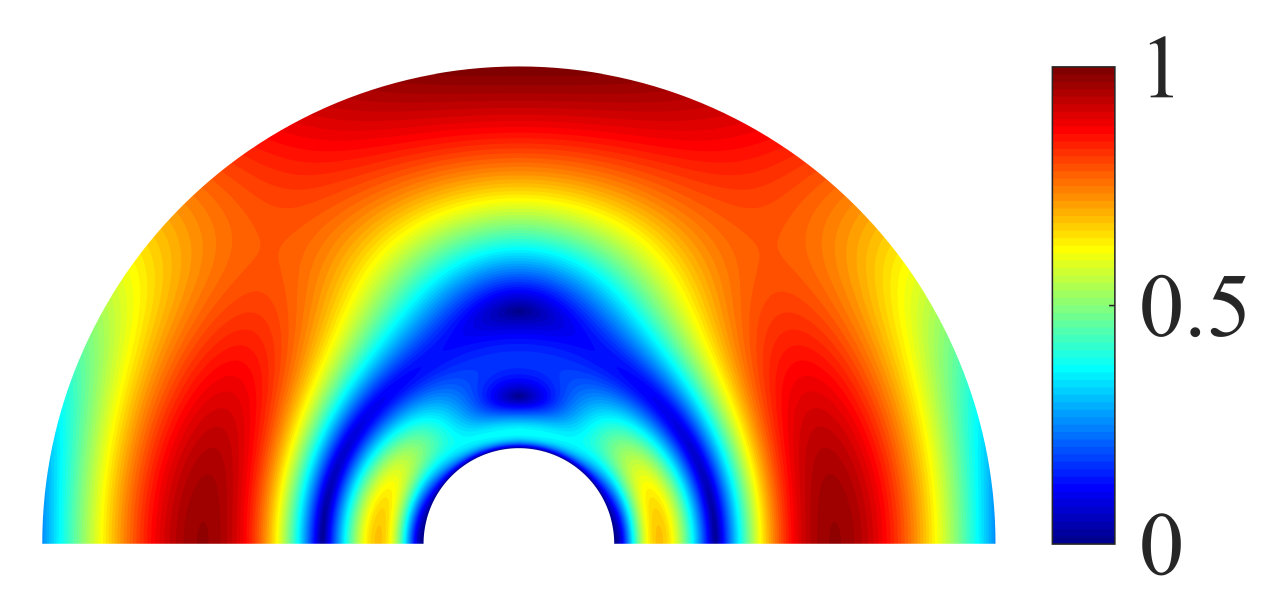}}
	\caption{Axisymmetric flow past a sphere: First four normalised spatial modes of the norm of the velocity field.}
	\label{fig:sphereVModes}
\end{figure}

\begin{figure}[!tb]
	\centering
	\subfigure[$m=1$]{\includegraphics[width=0.24\textwidth]{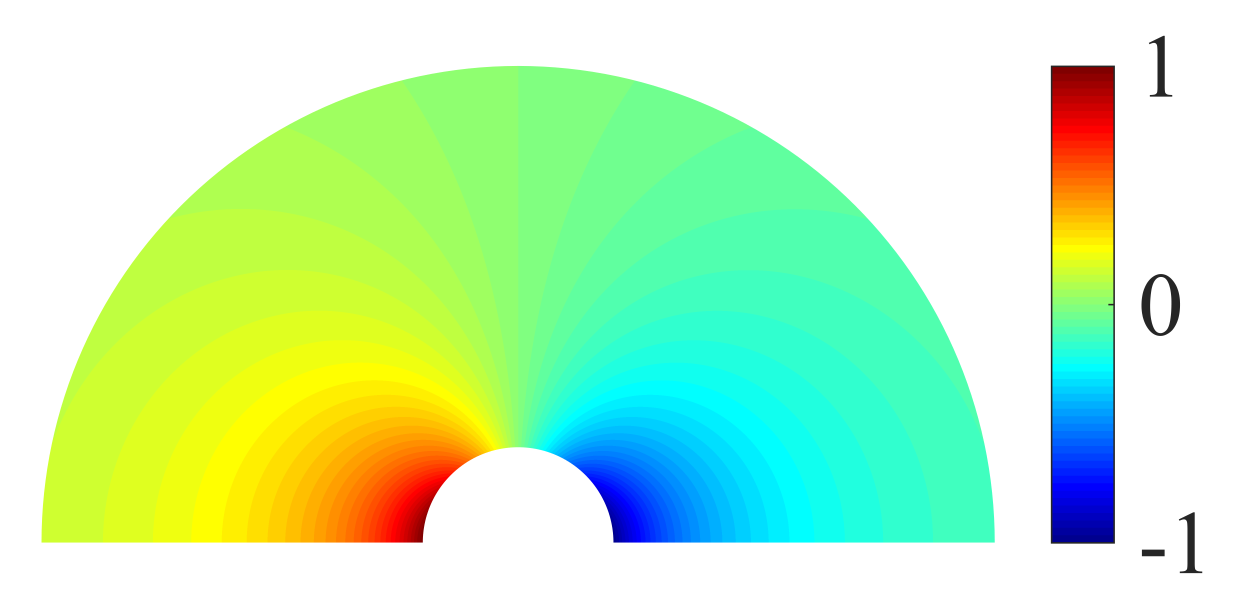}}
	\subfigure[$m=2$]{\includegraphics[width=0.24\textwidth]{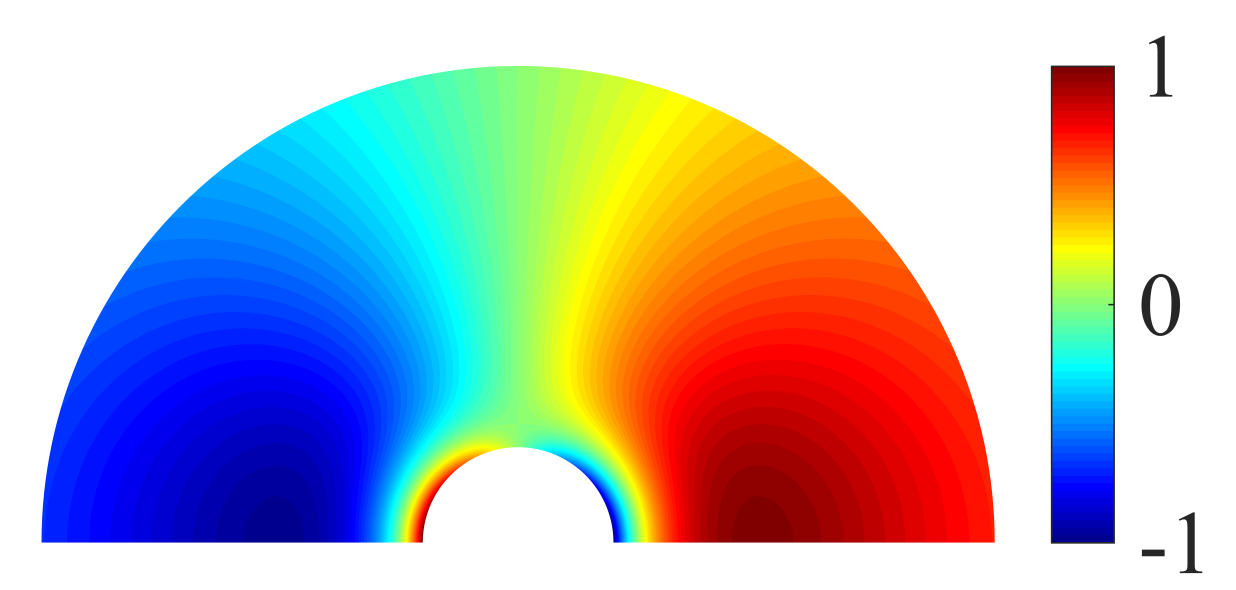}}
	\subfigure[$m=3$]{\includegraphics[width=0.24\textwidth]{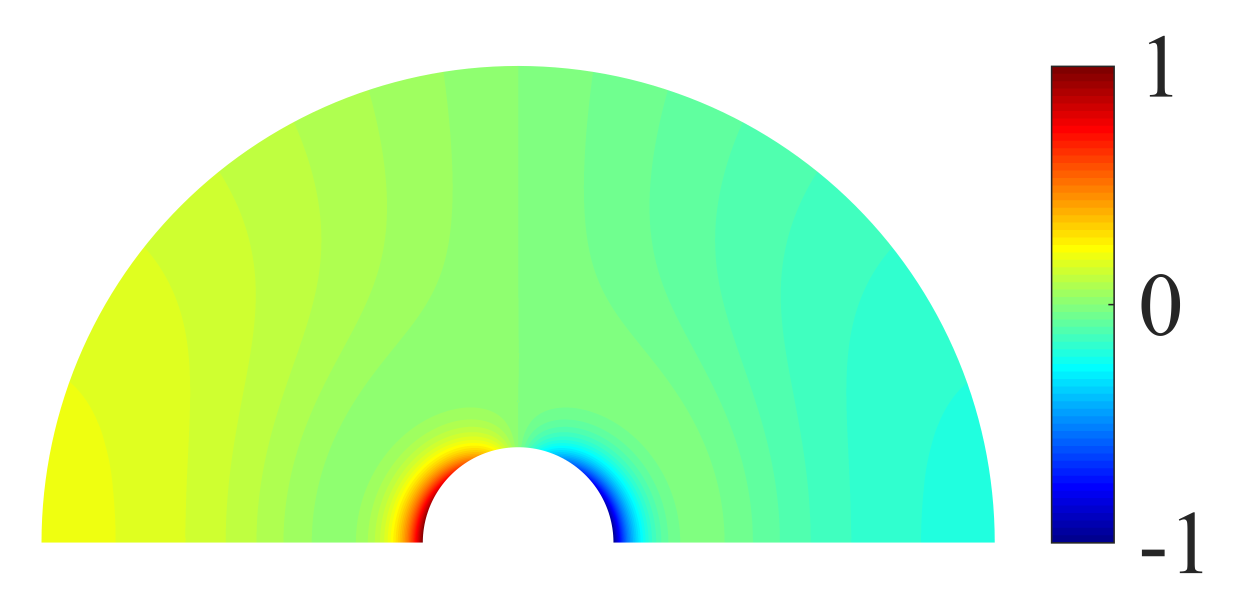}}
	\subfigure[$m=4$]{\includegraphics[width=0.24\textwidth]{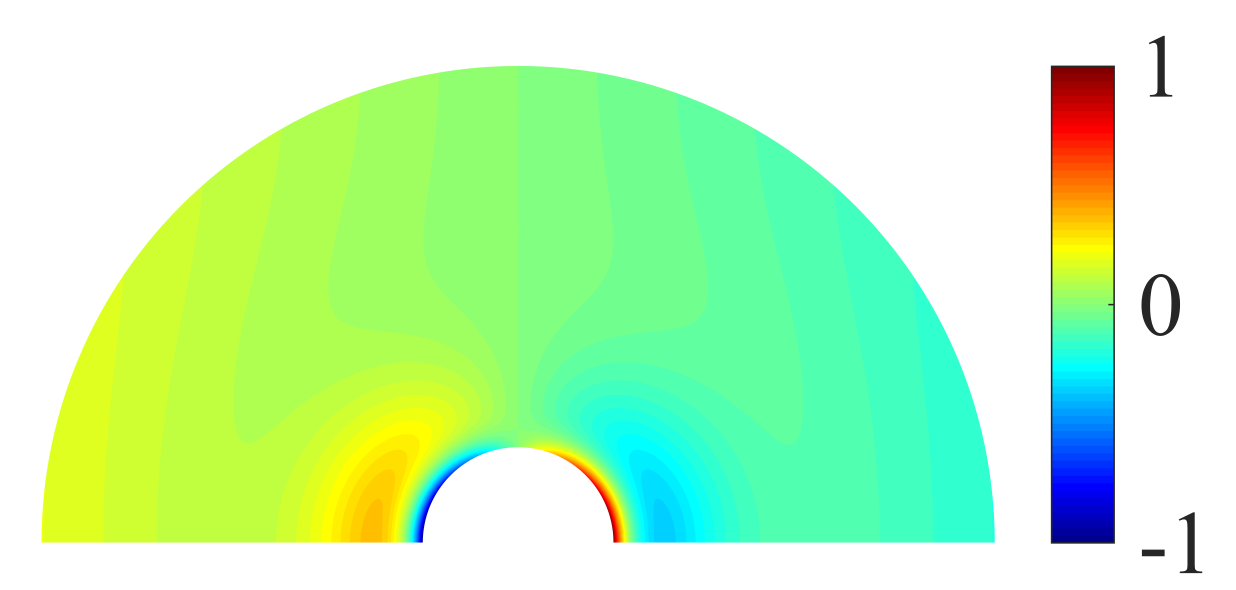}}
	\caption{Axisymmetric flow past a sphere: First four normalised spatial modes of the pressure field.}
	\label{fig:spherePModes}
\end{figure}

\begin{figure}[!tb]
	\centering
	\includegraphics[width=0.49\textwidth]{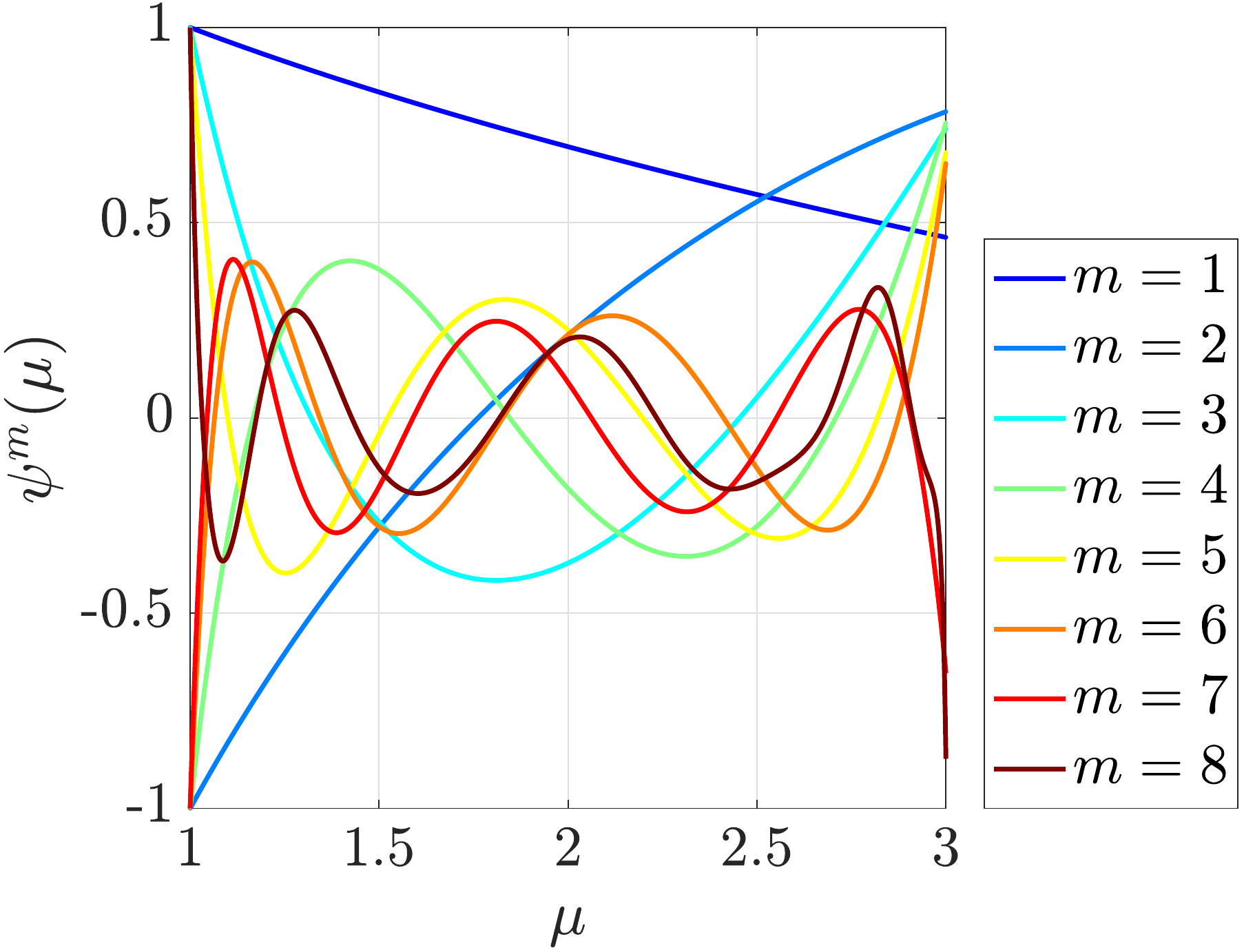}
	\caption{Axisymmetric flow past a sphere: First eight normalised parametric modes.}
	\label{fig:sphereParamModes}
\end{figure}

The proposed ROM is used to obtain the generalised solution of the parametric axisymmetric Stokes problem. The first four normalised modes of the magnitude of the velocity field and the pressure are shown in figures~\ref{fig:sphereVModes} and~\ref{fig:spherePModes}. 
The computation was performed using the second mesh with a degree of approximation $k=4$ for all the variables and with a mesh of 1,000 elements in the parametric dimension with also $k=4$. Figure~\ref{fig:sphereParamModes} shows the first eight normalised parametric modes computed.
It is worth noting that despite the different nature of the flow and the axisymmetric boundary condition, the parametric modes have a similar behaviour when compared to the modes obtained in the previous example. This is mainly attributed to the geometric parameter describing an analogous variation of the computational domain.

\begin{figure}[!tb]
	\centering
	\includegraphics[width=0.45\textwidth]{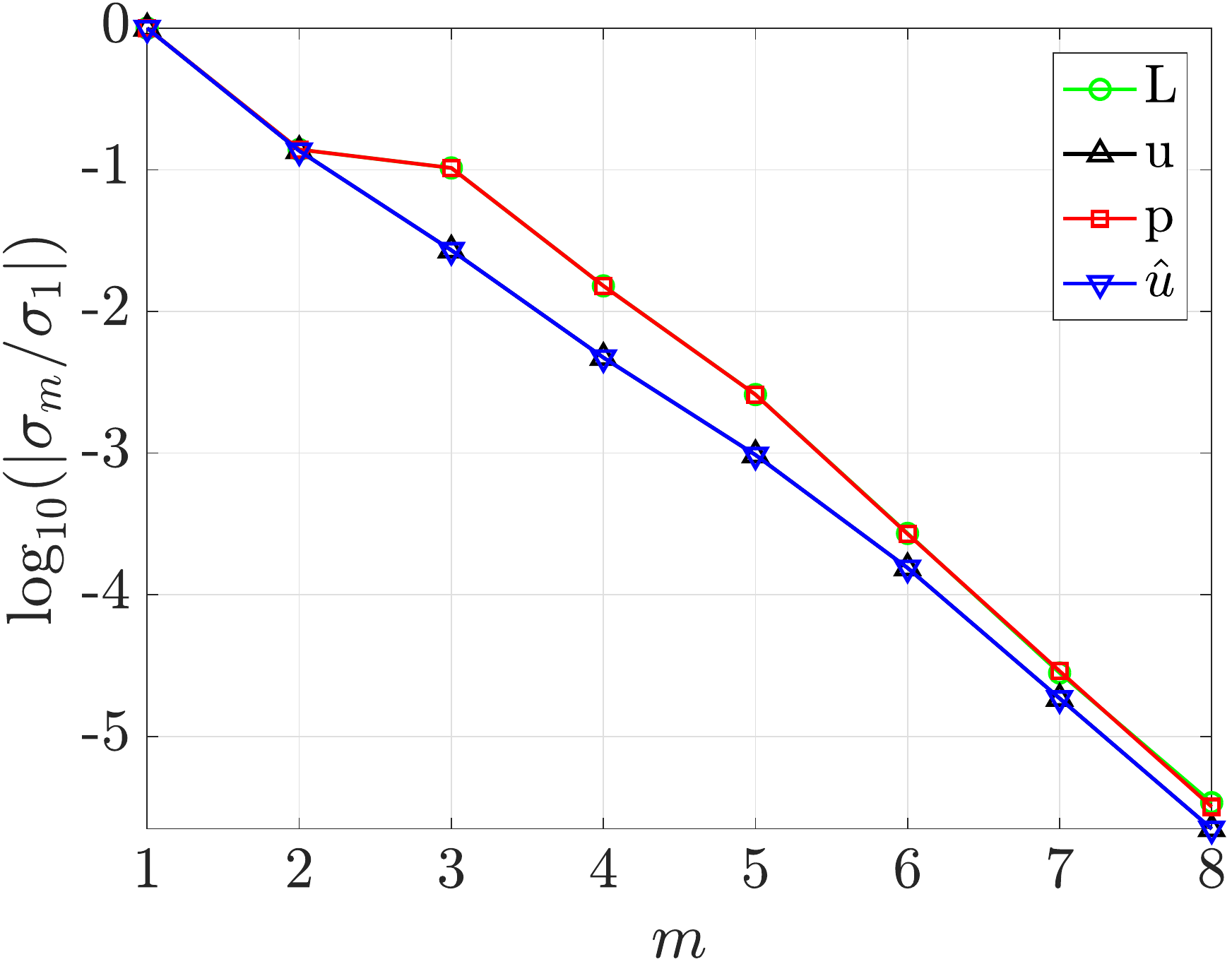}
	\caption{Axisymmetric flow past a sphere: Convergence of the mode amplitudes.}
	\label{fig:sphereAmplitudes}
\end{figure}

As in the previous example, the evolution of the relative amplitude of the modes is shown in~\ref{fig:sphereAmplitudes}. 
The rapid decrease shows that it is possible to compute a generalised solution to this problem with a very small number of modes. With eight modes the relative amplitude is already below $10^{-5}$.

\begin{figure}[!tb]
	\centering
	\subfigure[$\bL$] {\includegraphics[width=0.49\textwidth]{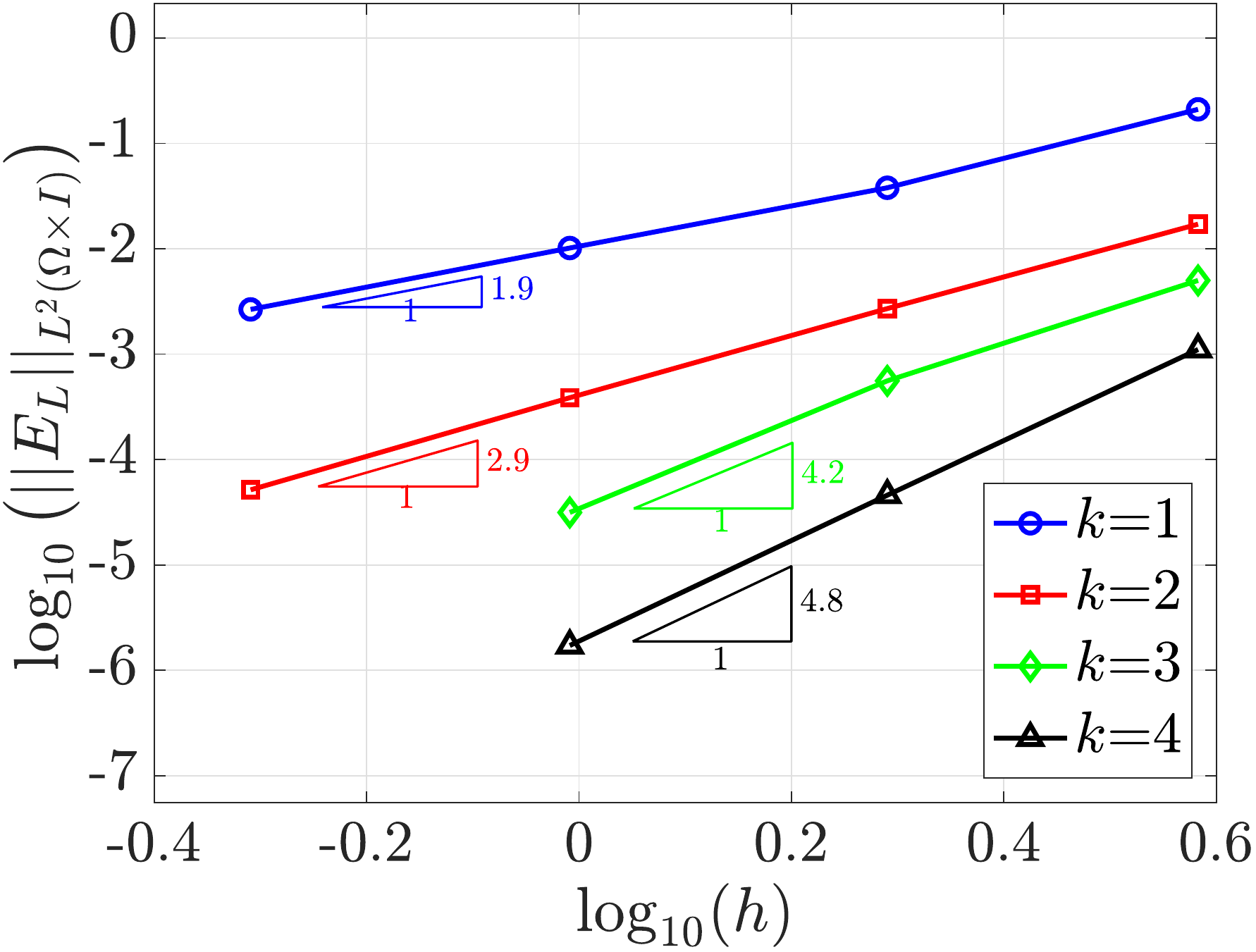}}
	\subfigure[$\bu$] {\includegraphics[width=0.49\textwidth]{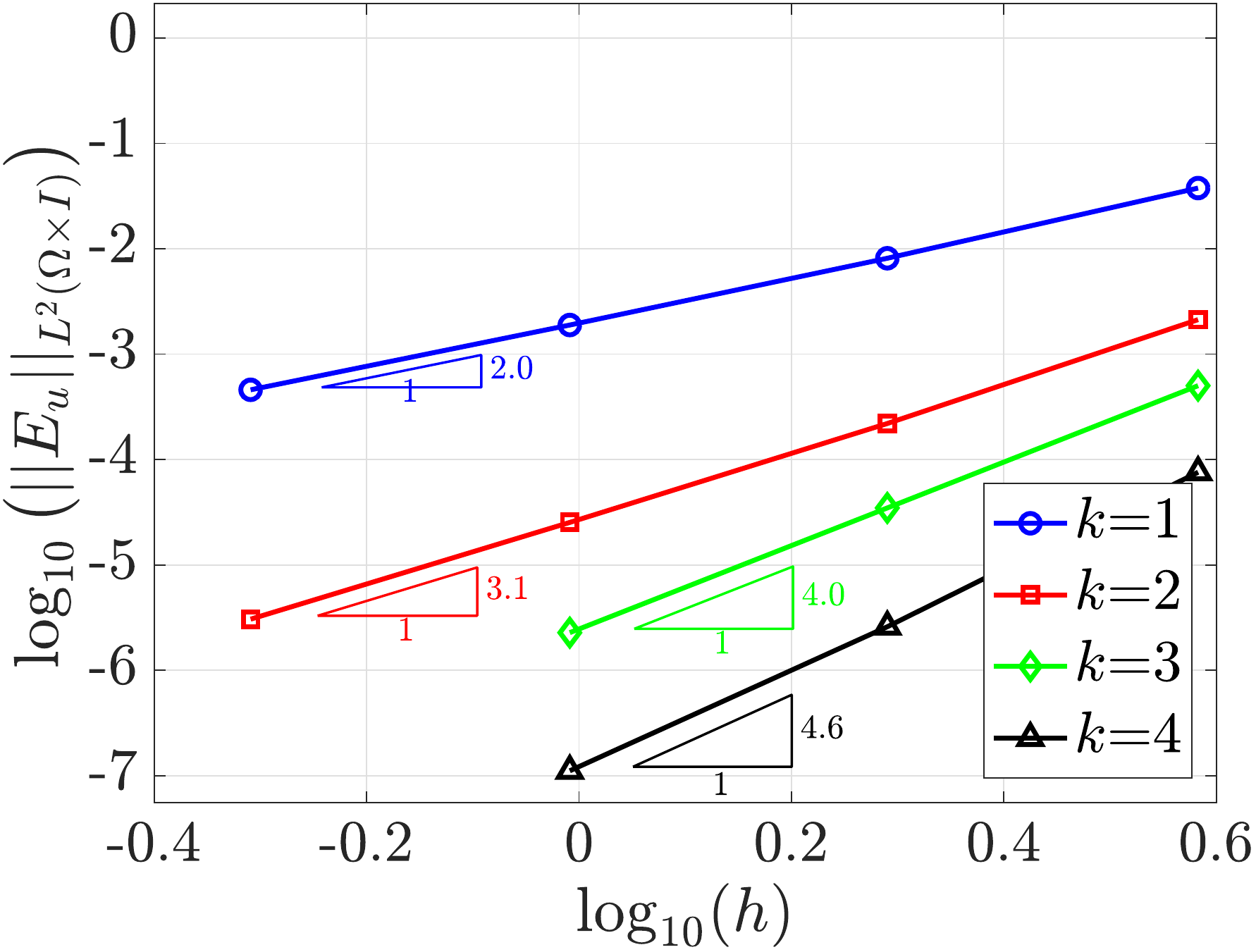}}
	\subfigure[$p$]   {\includegraphics[width=0.49\textwidth]{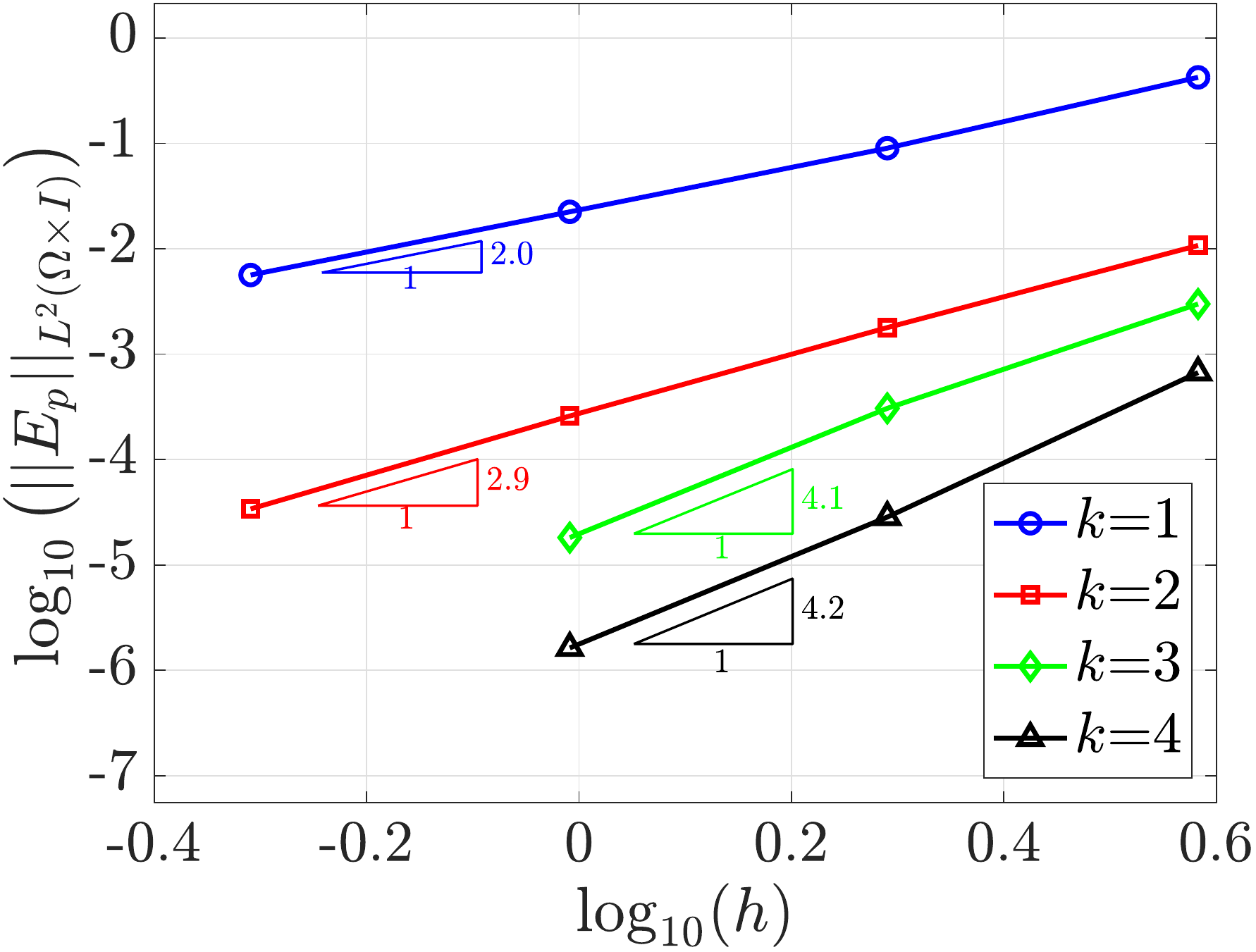}}
	\subfigure[$\bhu$]{\includegraphics[width=0.49\textwidth]{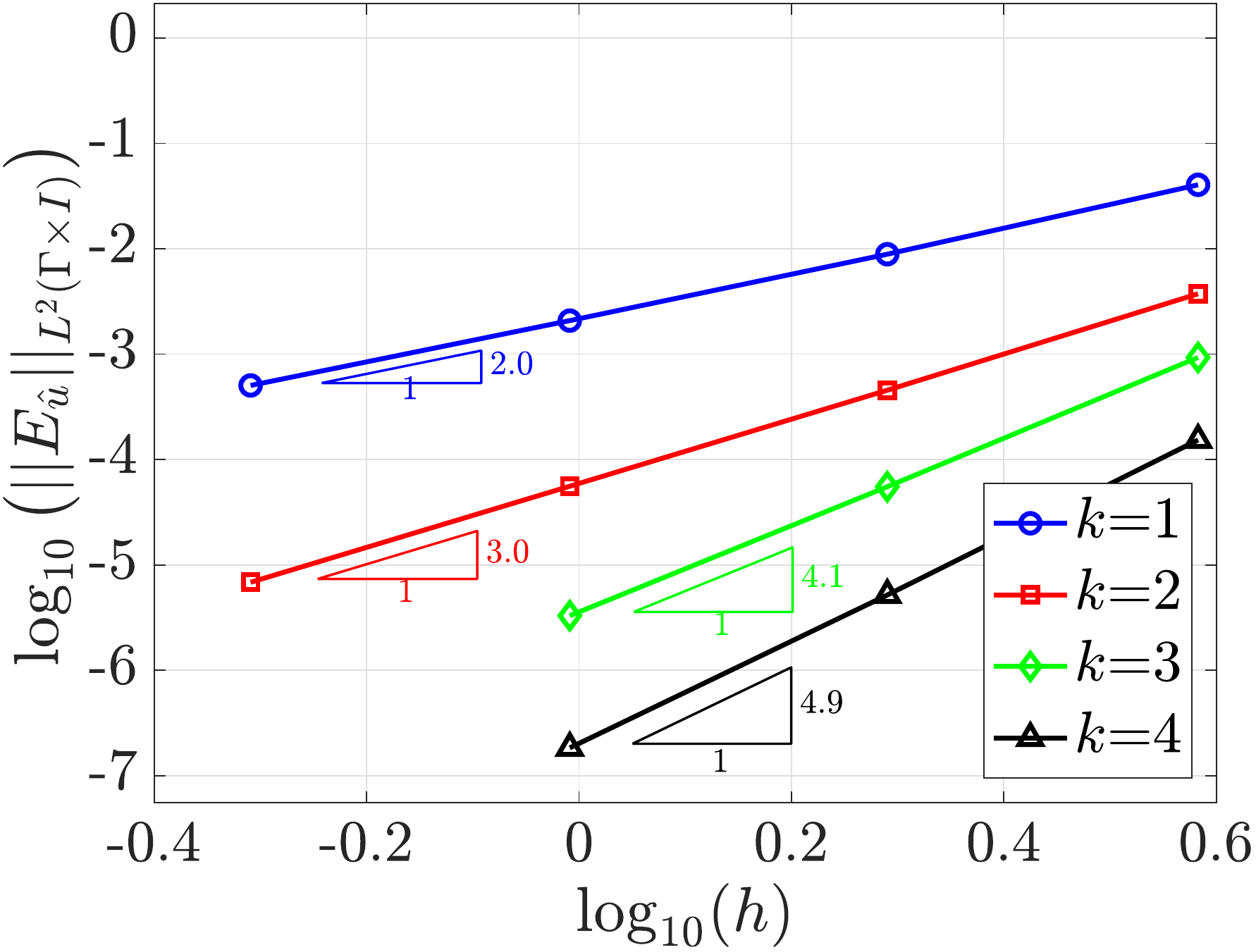}}
	\caption{Axisymmetric flow past a sphere: mesh convergence of the $\eltwo$ norm of the error for $\bL$, $\bu$, $p$ and $\bhu$.}
	\label{fig:sphereHConv}
\end{figure}

Next, the optimal approximation properties of the proposed HDG-PGD method are studied by performing a mesh convergent study. Figure~\ref{fig:sphereHConv} shows the evolution of the relative error in the $\eltwo(\Omega \times \bI)$ norm as a function of the characteristic element size, $h$, for different orders of approximation and for all the variables of the HDG formulation.
The optimal rate of convergence, equal to $h^{k+1}$, is approximately observed for all the variables. 

\begin{figure}[!tb]
	\centering
	\subfigure[$\mu_1=1$ \label{fig:sphereDragA}]{\includegraphics[width=0.32\textwidth]{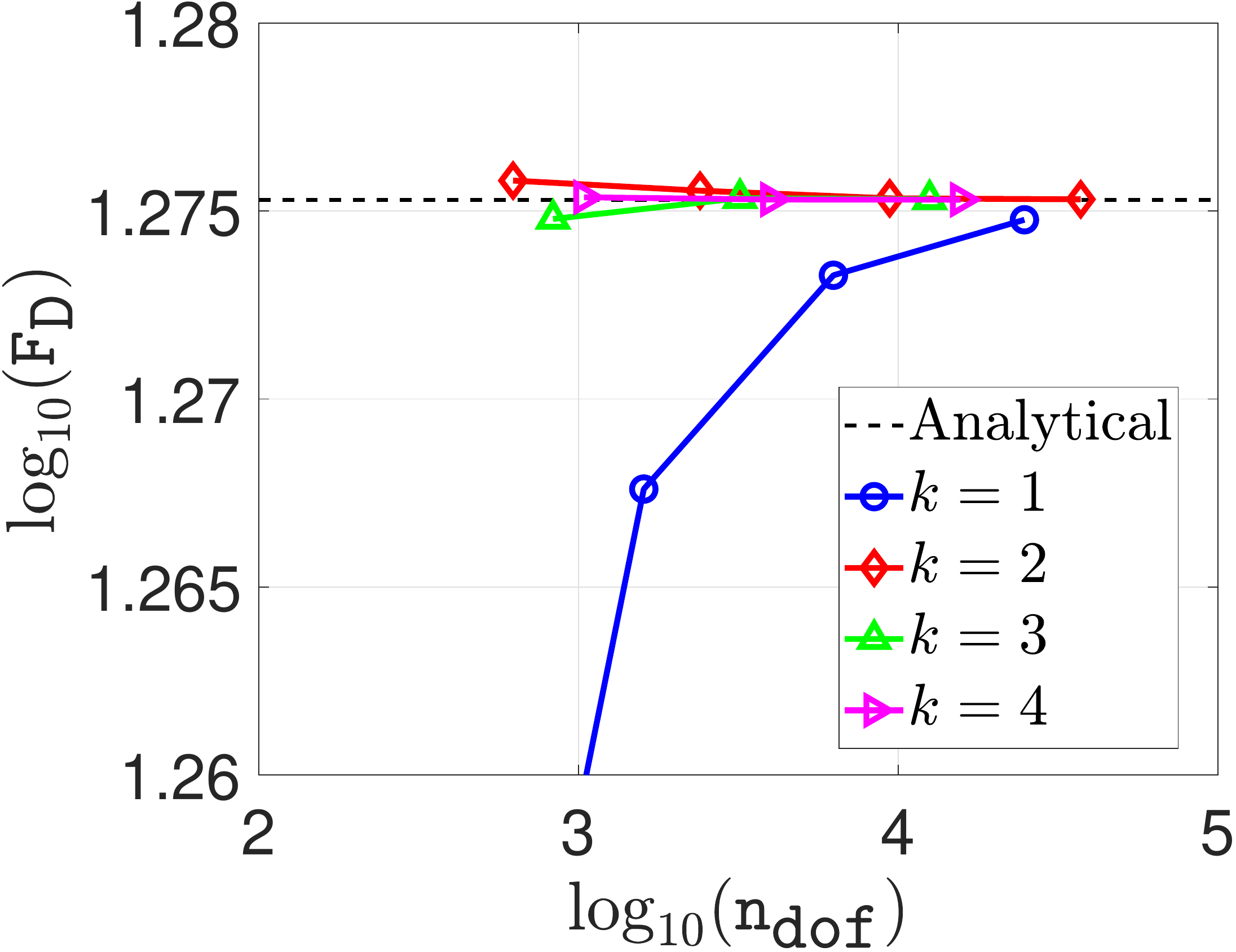}} 
	\subfigure[$\mu_1=2$ \label{fig:sphereDragB}]{\includegraphics[width=0.32\textwidth]{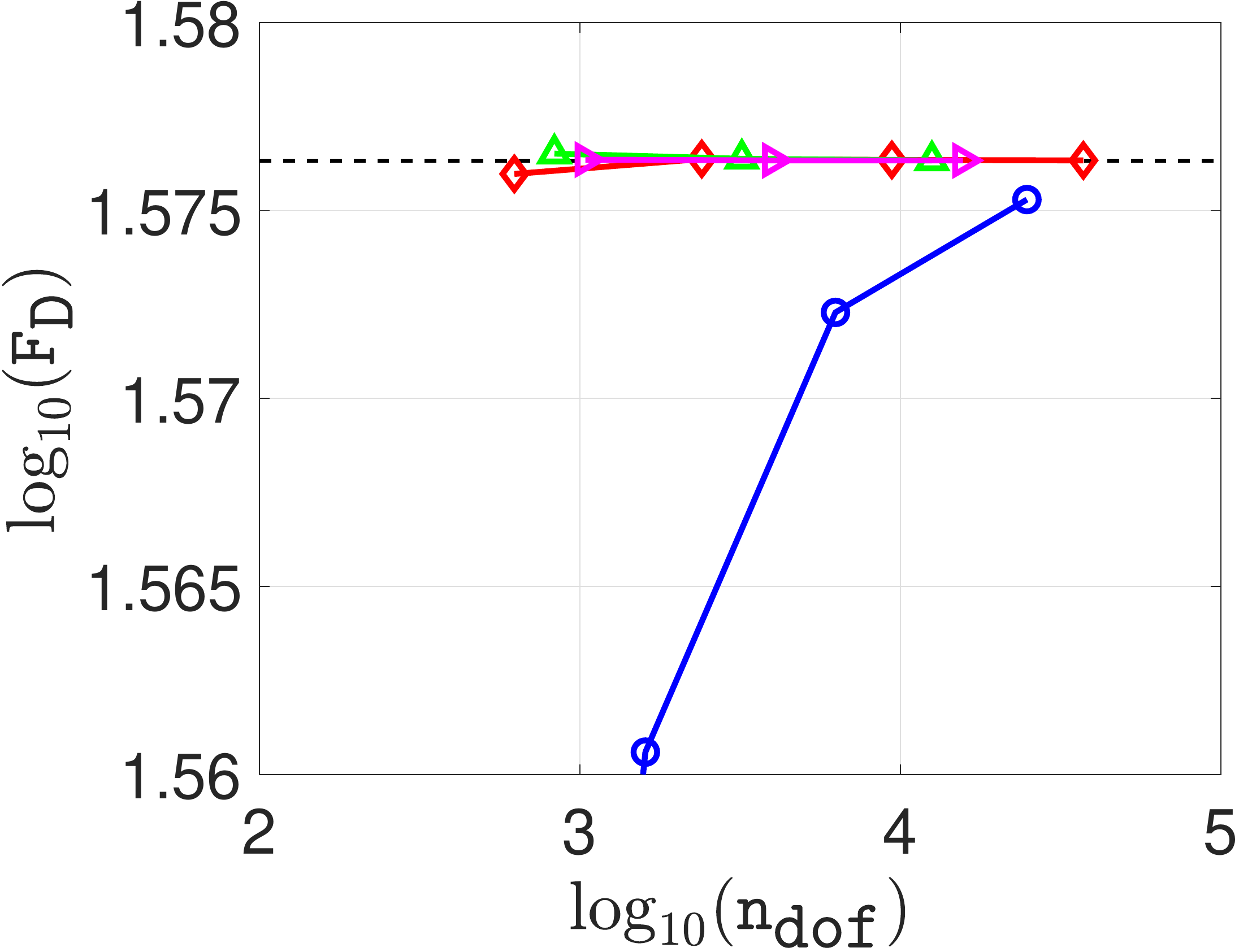}}
	\subfigure[$\mu_1=3$ \label{fig:sphereDragC}]{\includegraphics[width=0.32\textwidth]{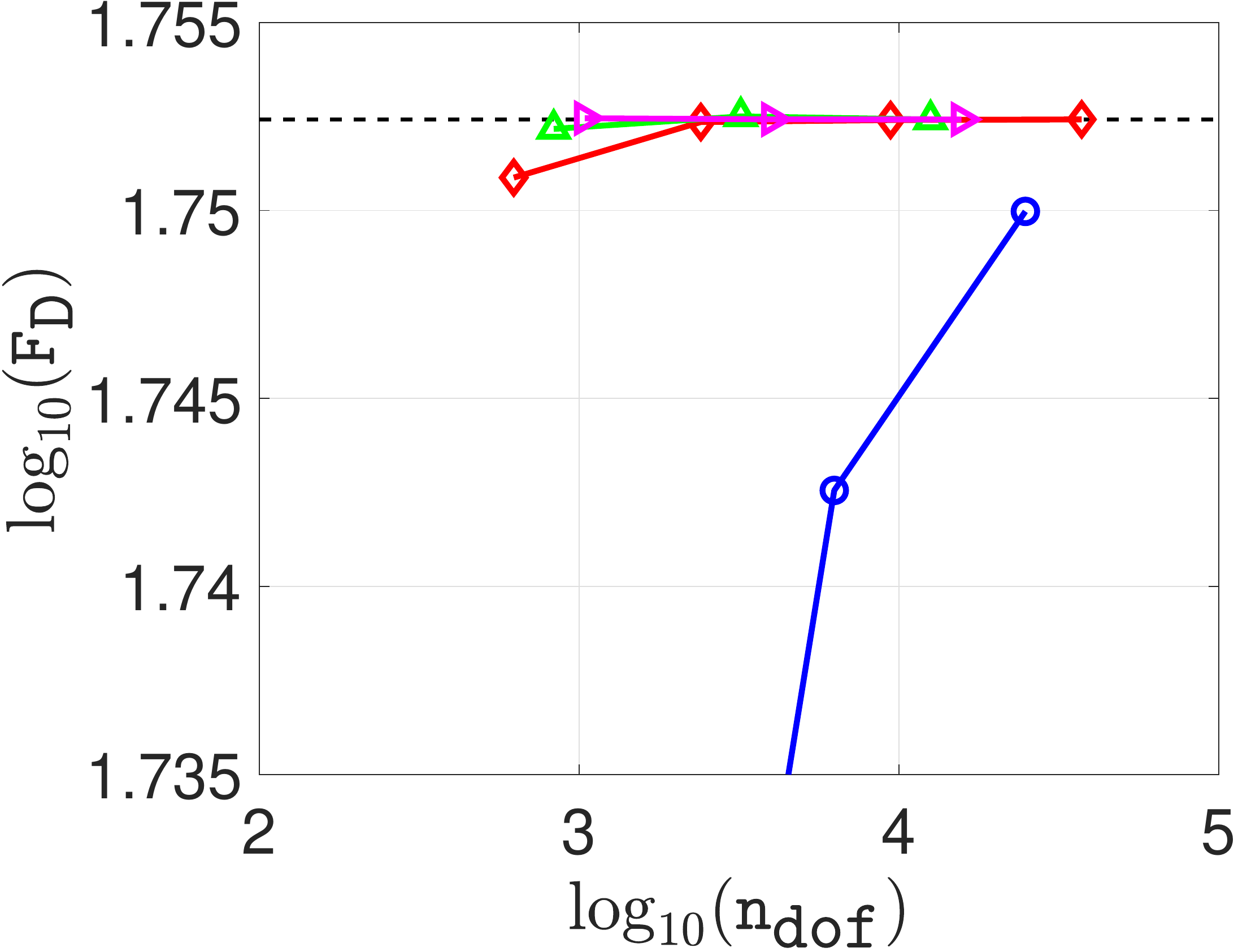}}
	\caption{Axisymmetric flow past a sphere: evolution of the error in the drag force as the number of modes is increased for three different geometric configurations.}
	\label{fig:sphereDrag}
\end{figure}

Finally, the accuracy of the HDG-PGD approach on the drag force is studied for three different configurations corresponding to $\mu_1=1$, $\mu_1=2$ and $\mu_1=3$. Figure~\ref{fig:sphereDrag} shows evolution of the error in the drag force as the number of of degrees of freedom is increased for the three different geometric configurations and for different orders of approximation.
The number of degrees of freedom refers to the size of the HDG global problem as this is the most time consuming part of the spatial iteration.

The results show the variation of the drag force induced by the variation of the geometric parameter and how the generalised solution produces accurate results for any value of the geometric parameter. In all cases, convergence to exact value is observed, and the superiority of using high-order approximations is clearly appreciated. For the first configuration, the results in figure~\ref{fig:sphereDragA} show that with a linear approximation requires the solution of a global problem with 24,832 degrees of freedom to obtain relative error in the drag force of 0.0181. In contrast , using a quartic approximation, the error in the first mesh is 0.0021, solving a global problem with only 416 degrees of freedom, that is an error one order of magnitude lower with almost 20 times less degrees of freedom.

The results also show that for higher values of the geometric parameter the solution is slightly more difficult to capture and the number of degrees of freedom required is slightly higher. In fact, the advantages of high-order approximations are more noticeable for the case of $\mu_1=3$. 

\subsection{Axisymmetric Stokes flow around two micro-swimmers} \label{sc:swimmers}

The next example considers the Stokes flow around the so-called \textit{push-me-push-you microswimmer}, proposed in~\cite{avron2005pushmepullyou}. This swimmer consists of two spherical bladders that have the ability to change their mutual distance and individual volume, whilst maintaining the total volume of the two spheres. The swimmer is placed in a cylindrical channel of length $L$ and diameter $D$. 

Two geometric parameters are considered in this example. The first one, $\mu_1 \in \I_1 = [-1,1]$, controls the radius of the two spheres in such a way that the total volume of the two spheres is maintained. The second parameter, $\mu_2 \in \I_2 = [-3,2]$, controls the distance between the centre of the two spheres. The value of $\mu_1=-1$ corresponds to the configuration where the radius of the first sphere is $R_1=0.3096$ and the radius of the second sphere is $R_2=0.116$, whereas the value of $\mu=1$ corresponds to the opposite situation, with $R_1=0.116$ and $R_2=0.3096$. The value of $\mu_2=-3$ corresponds to the case where the distance between the spheres is maximum, with the centres of the spheres placed at $(-3,0)$ and $(3,0)$ respectively. The value of $\mu_2=2$ corresponds to the case where the distance between the spheres is minimum, with the centres of the spheres placed at $(-0.5,0)$ and $(0.5,0)$ respectively.

\begin{figure}[!tb]
	\centering
	\includegraphics[width=0.8\textwidth]{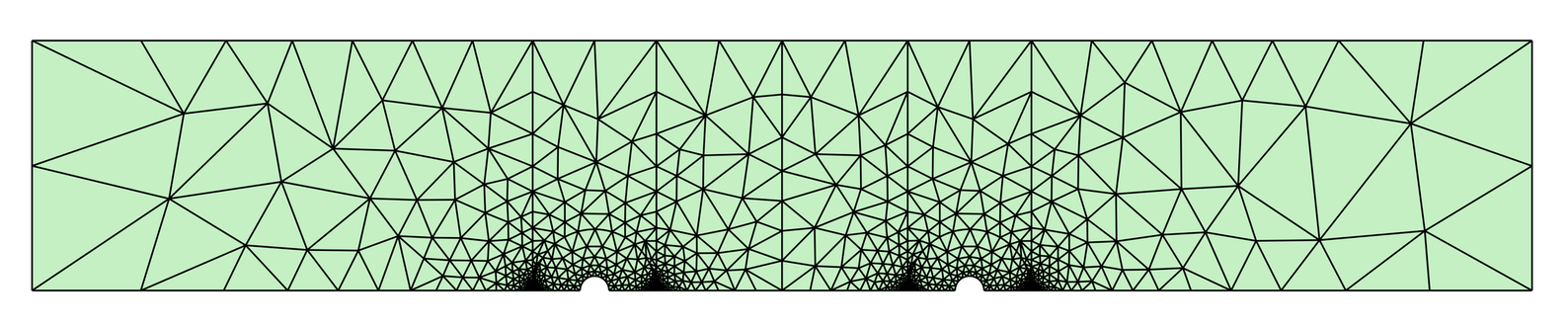}
	\caption{Axisymmetric flow around two micro-swimmers: Computational mesh.}
	\label{fig:swimmerMesh}
\end{figure}

Using the axial symmetry of the problem, the reference domain is chosen as $\Omega = \left([-L,L] \times [0,H]\right) \setminus \left( \mathcal{B}^+ \cup \mathcal{B}^- \right)$, where 
\begin{equation}
\mathcal{B}^{\pm} = \{ \bX \in \mathbb{R}^2 \; | \; \| \bX \pm \bX_0 \| \leq \Rref \},
\end{equation} 
where $L=6$, $H=2$, $\bX_0 = (1.5,0)$ and $\Rref = 0.116$.
Figure~\ref{fig:swimmerMesh} shows the triangular mesh of the reference domain used for this numerical example.
The mesh has 1,426 elements, leading to a system in the HDG global problem of 22,260 equations for a degree of approximation $k=4$.

On the left part of the boundary a Dirichlet boundary condition, corresponding to a horizontal velocity of magnitude one, is imposed. On the right part of the boundary a homogeneous Neumann boundary condition is imposed. On the surface of the two spheres a no-slip boundary condition is enforced and on the rest of the boundary a slip boundary condition is imposed. 

The geometric mapping used in this example is detailed in~\ref{sc:mappingSwimmer}.

\begin{figure}[!tb]
	\centering
	\subfigure[$m=1$]{\includegraphics[width=0.49\textwidth]{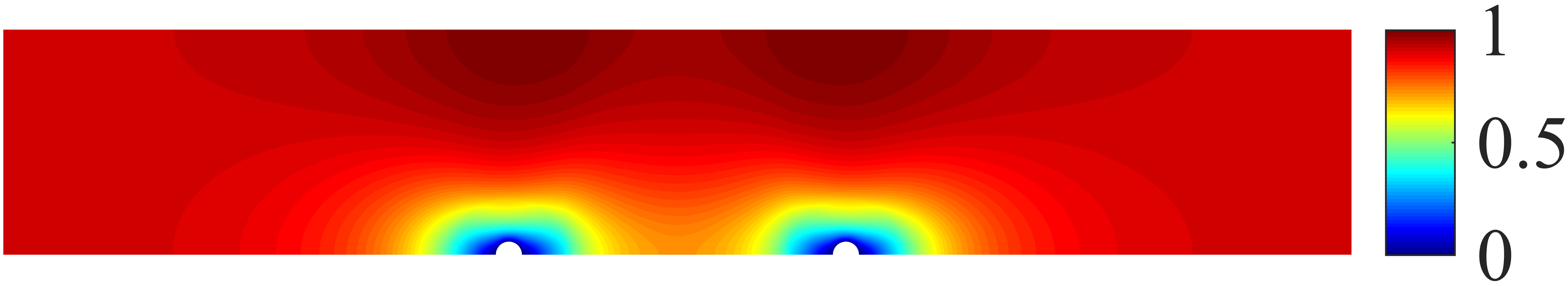}}
	\subfigure[$m=2$]{\includegraphics[width=0.49\textwidth]{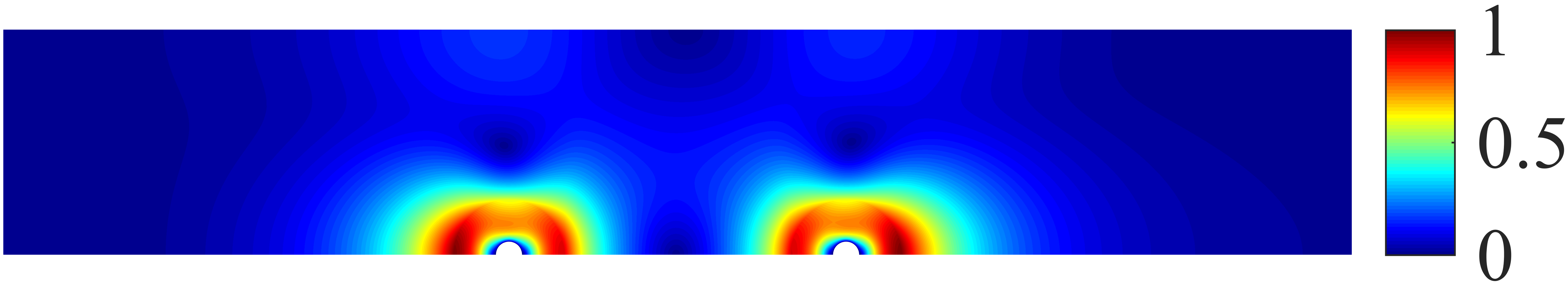}}
	\subfigure[$m=3$]{\includegraphics[width=0.49\textwidth]{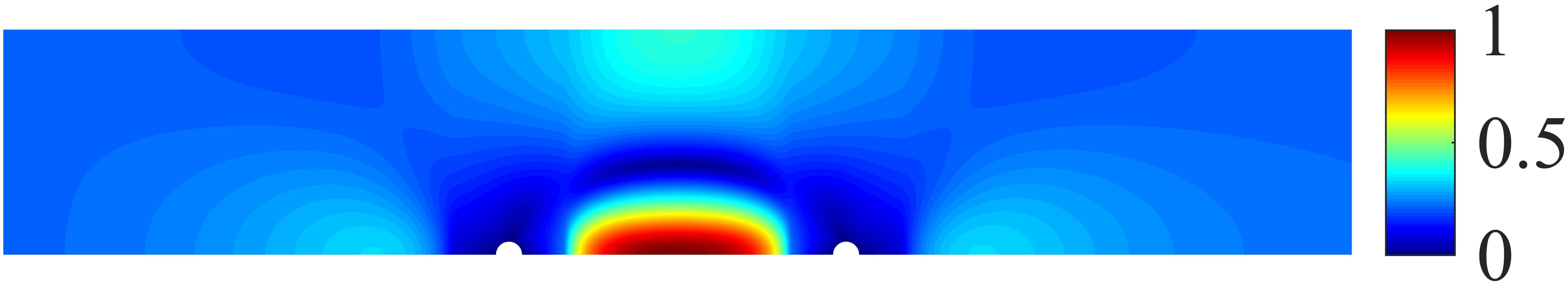}}
	\subfigure[$m=4$]{\includegraphics[width=0.49\textwidth]{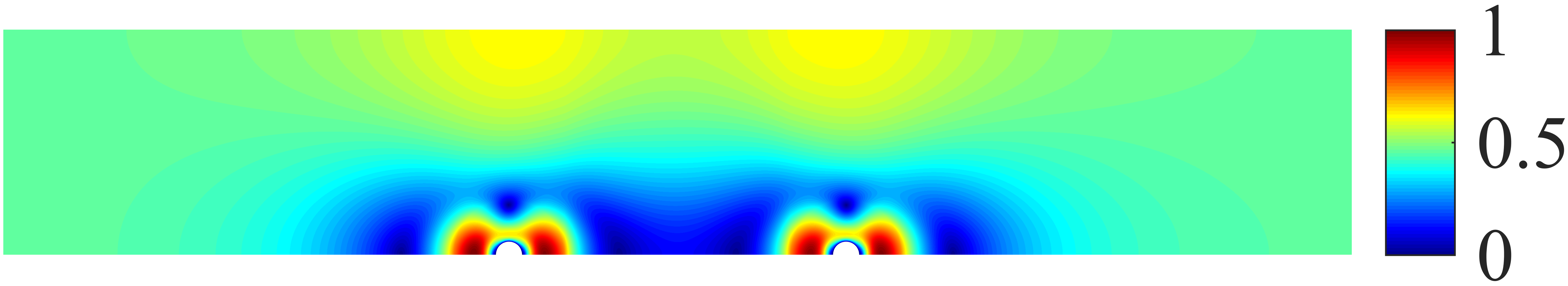}}
	\caption{Axisymmetric flow around two micro-swimmers: First four normalised spatial modes of the velocity field.}
	\label{fig:swimmerVModes}
\end{figure}

\begin{figure}[!tb]
	\centering
	\subfigure[$m=1$]{\includegraphics[width=0.49\textwidth]{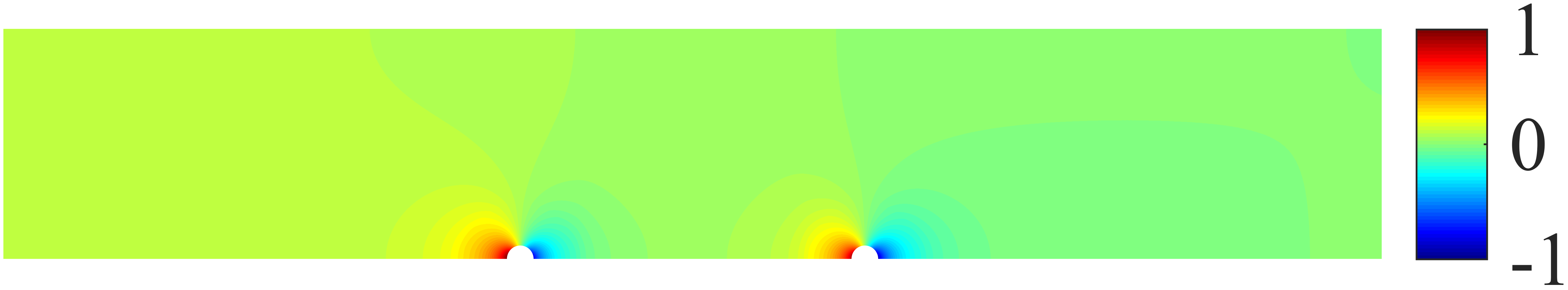}}
	\subfigure[$m=2$]{\includegraphics[width=0.49\textwidth]{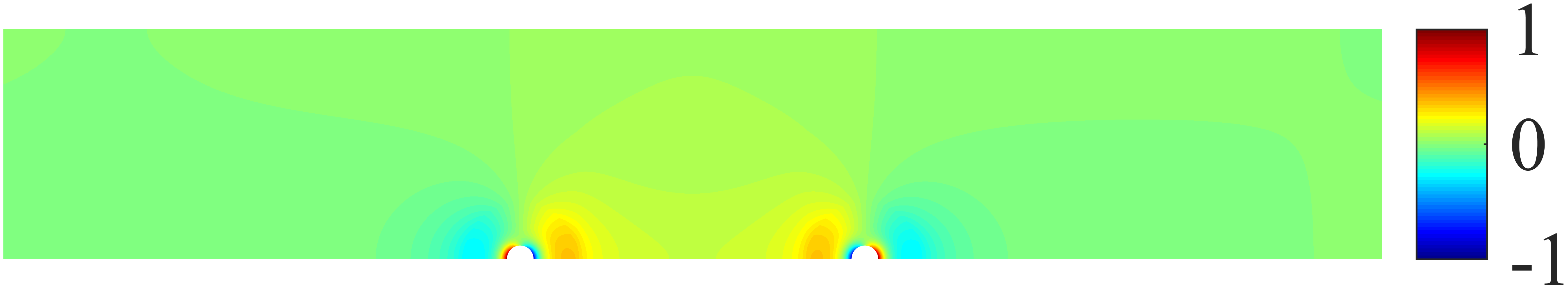}}
	\subfigure[$m=3$]{\includegraphics[width=0.49\textwidth]{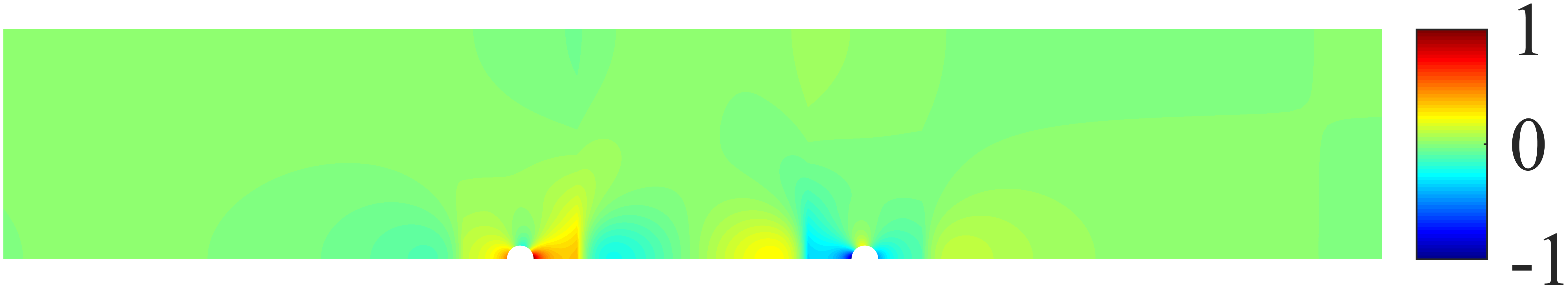}}
	\subfigure[$m=4$]{\includegraphics[width=0.49\textwidth]{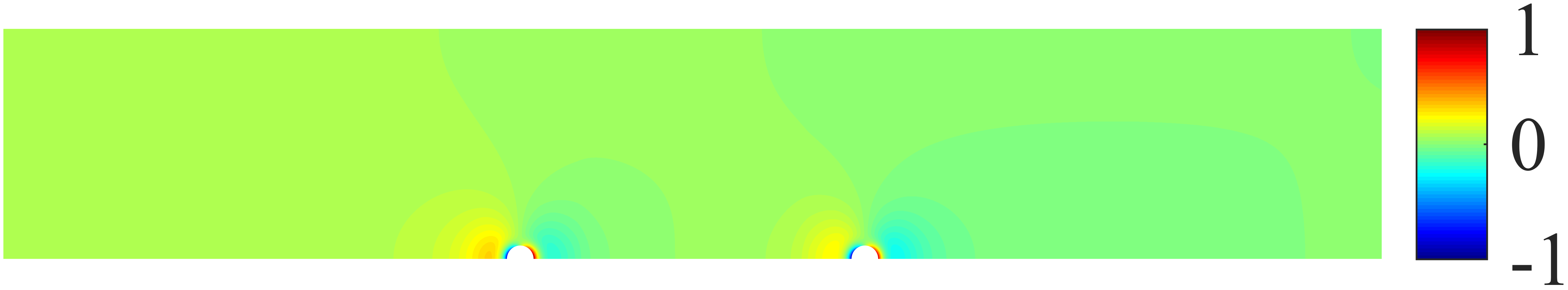}}
	\caption{Axisymmetric flow around two micro-swimmers: First four normalised spatial modes of the pressure field.}
	\label{fig:swimmerPModes}
\end{figure}

The first four spatial modes for the velocity and pressure computed with the proposed HDG-PGD are shown in figures~\ref{fig:swimmerVModes} and~\ref{fig:swimmerPModes}. 
The computation was performed using the mesh of figure~\ref{fig:swimmerMesh} with a degree of approximation $k=4$ for all the variables and with a mesh of 10,000 elements in each parametric dimension with also $k=4$. It is worth noting that the cost of the one-dimensional parametric problems is negligible when compared to the cost of the spatial iteration. Therefore, a large number of elements is used in the parametric dimension to ensure that the variation induced by the geometric parameters are captured with no a priori knowledge of the solution.

\begin{figure}[!tb]
	\centering
	\subfigure[\label{fig:swimmerParamModes1}]{\includegraphics[width=0.49\textwidth]{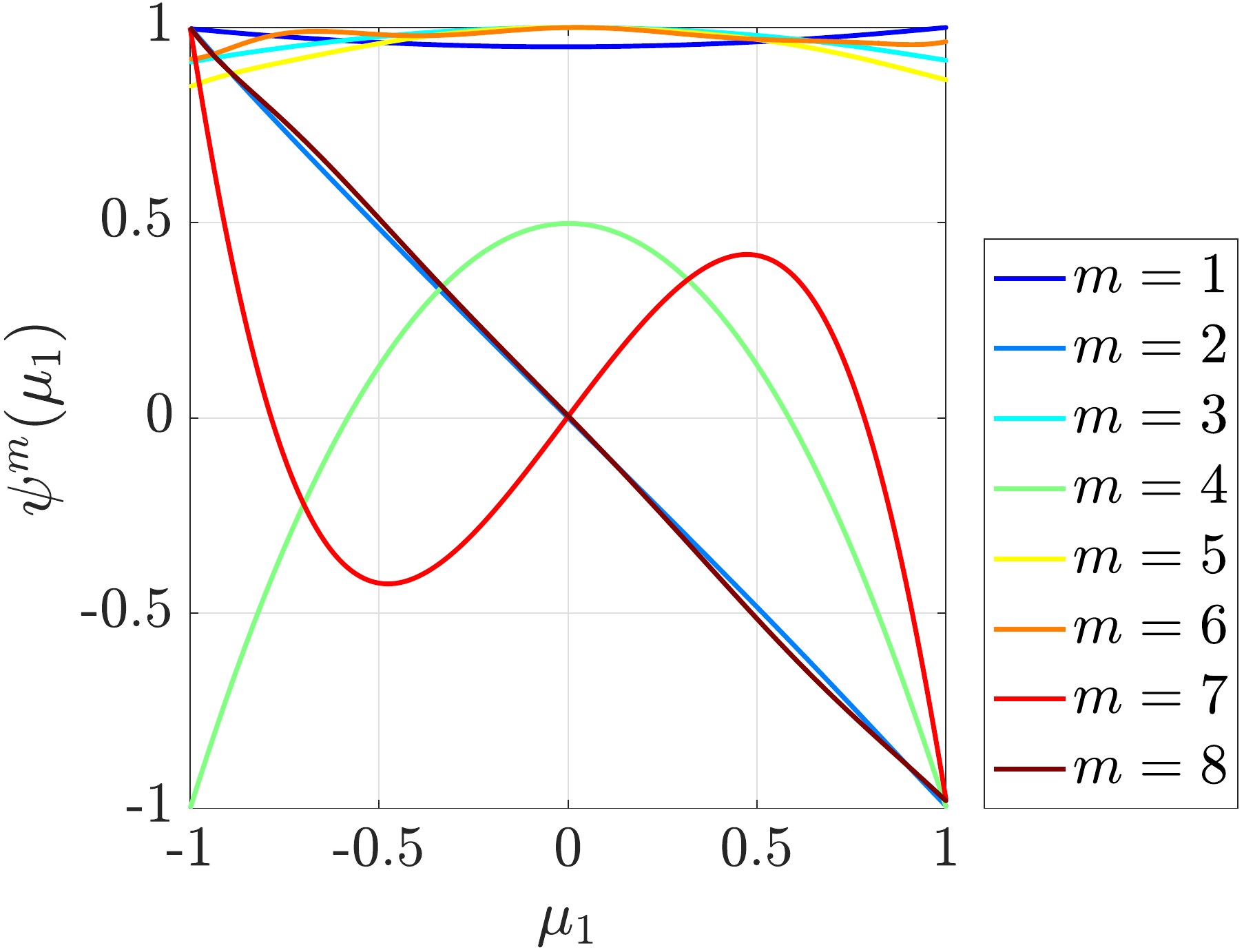}}
	\subfigure[\label{fig:swimmerParamModes2}]{\includegraphics[width=0.49\textwidth]{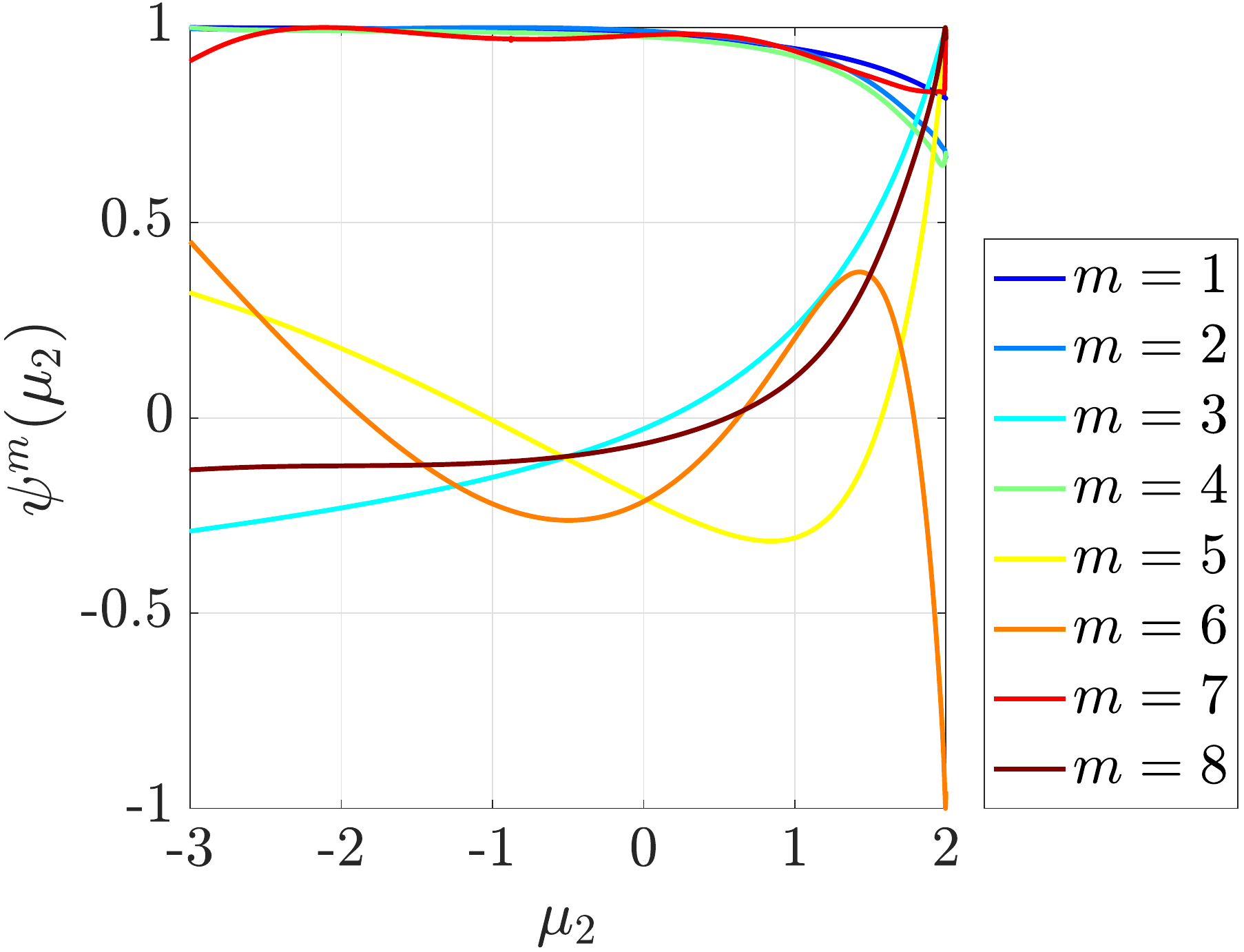}}
	\caption{Axisymmetric flow around two micro-swimmers: First eight normalised parametric modes.}
	\label{fig:swimmerParamModes}
\end{figure}

Figure~\ref{fig:swimmerParamModes} shows the first eight normalised parametric modes computed. 
Contrary to the previous examples, in this example there are more parametric modes that have an important influence over the whole range of values for both $\mu_1$ and $\mu_2$. For instance, in figure~\ref{fig:swimmerParamModes1} the first, third, fifth and six parametric modes have a normalised value near one for the whole range of values of $\mu_1$. A similar behaviour is observed for the second parameter $\mu_2$. In addition, the second parameter, corresponding to the distance between the spheres it can be observed that many of the modes have a much more relevant influence near $\mu_2 = 2$. This is expected as this configuration corresponds to the case where the distance between the spheres is minimum and therefore induces an important variation in the flow field because the first sphere will influence the flow that is reaching the second sphere.

\begin{figure}[!tb]
	\centering
	\includegraphics[width=0.45\textwidth]{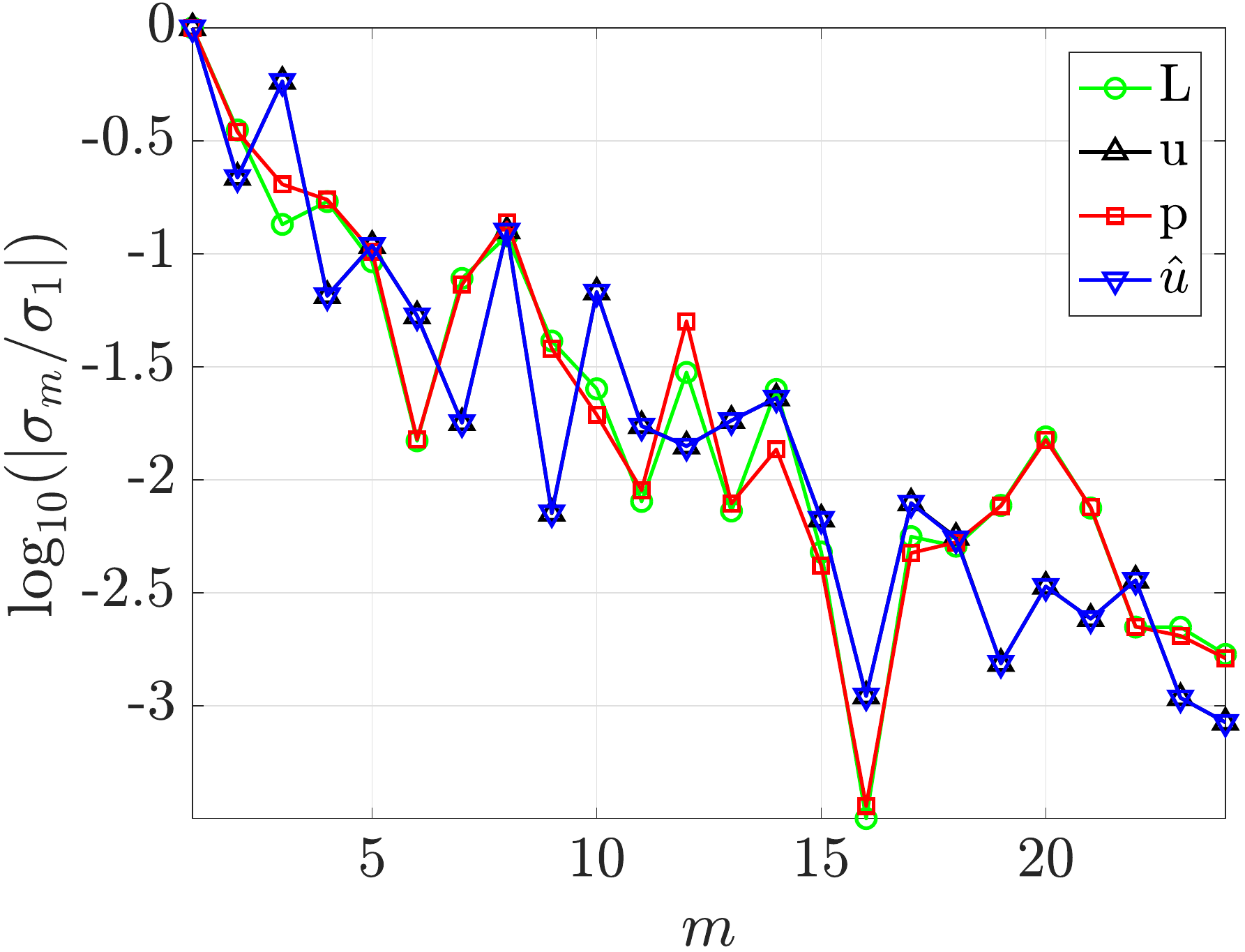}
	\caption{Axisymmetric flow around two micro-swimmers: Convergence of the mode amplitudes.}
	\label{fig:swimmerAmplitudes}
\end{figure}

The evolution of the relative amplitude of the modes is displayed in figure~\ref{fig:swimmerAmplitudes}. 
The results show that with 24 modes all the relative amplitude of the hybrid variable, used to check convergence, is below $10^{-3}$. A slower decrease of the relative amplitudes when compared with the previous examples can be observed. This is attributed to two factors. First, this problem considers two geometric parameters and, second, the range of variation of the distance is relatively high when compared to the minimum radius of the spheres. 

\begin{figure}[!tb]
	\centering
	\subfigure[$\mu_1=-1$, $\mu_2=-3$ \label{fig:swimmerCofigurations1V}]{\includegraphics[width=0.32\textwidth]{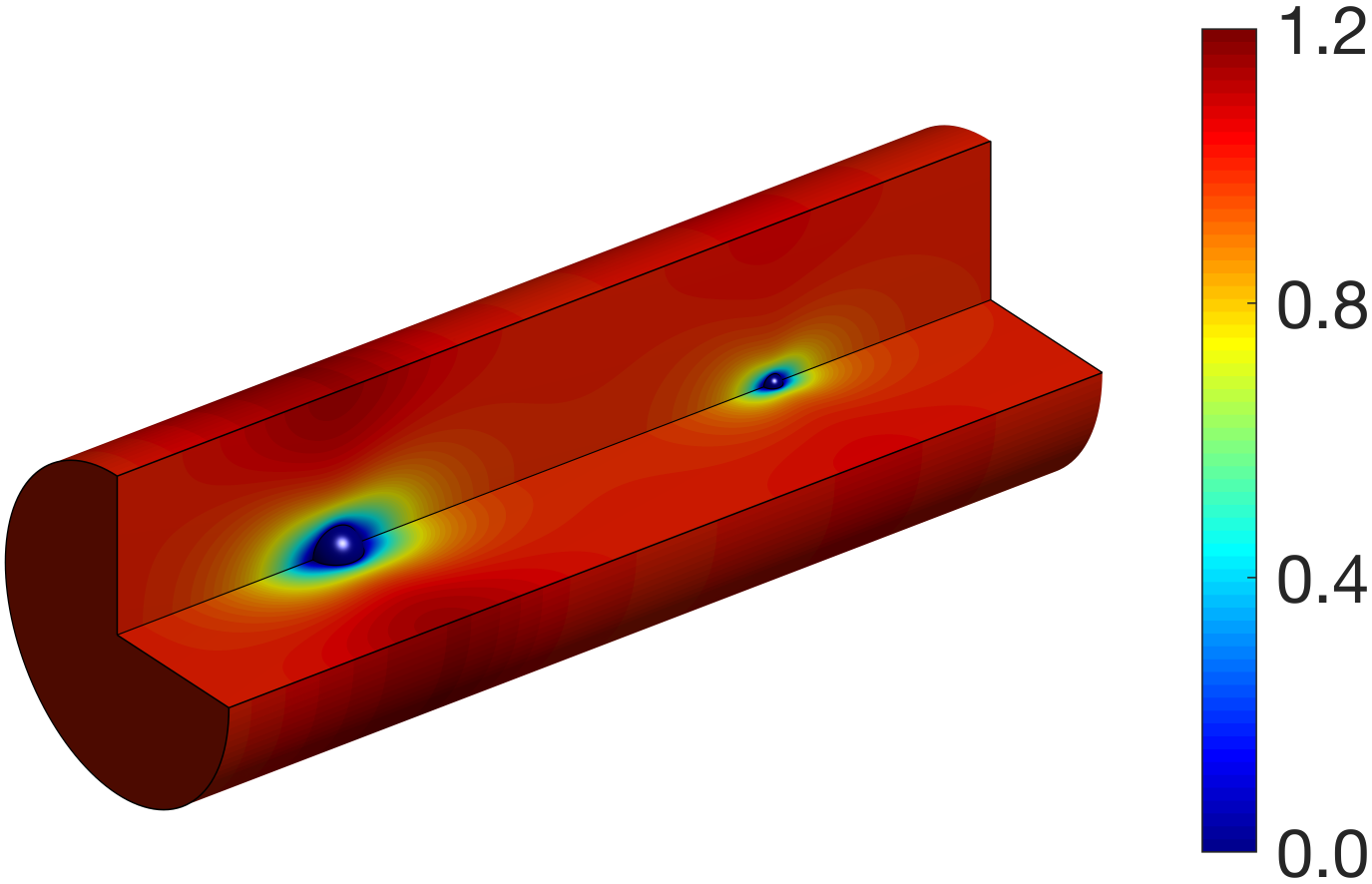}}
	\subfigure[$\mu_1=0$, $\mu_2=0$   \label{fig:swimmerCofigurations2V}]  {\includegraphics[width=0.32\textwidth]{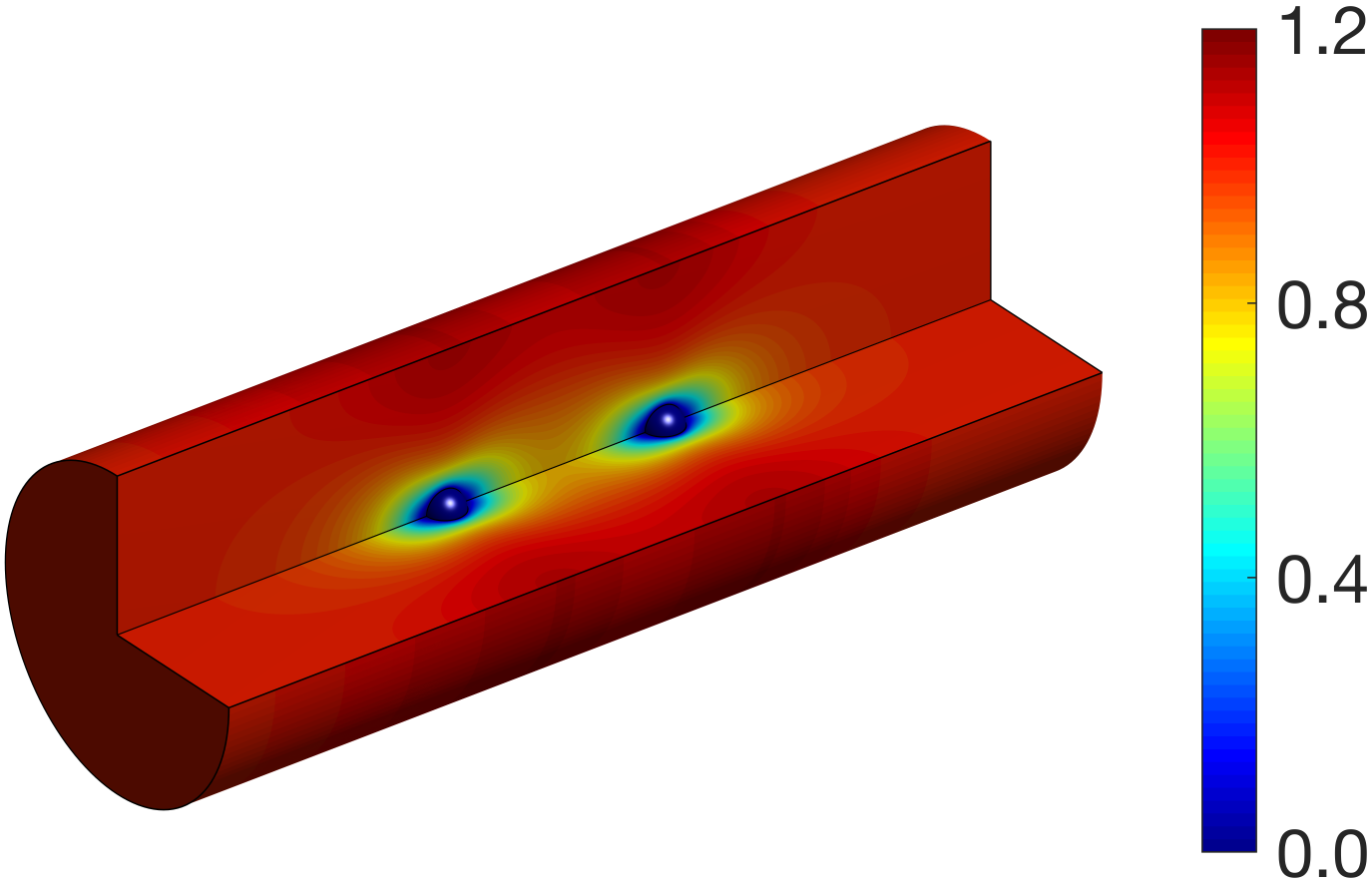}}
	\subfigure[$\mu_1=1$, $\mu_2=2$   \label{fig:swimmerCofigurations3V}]  {\includegraphics[width=0.32\textwidth]{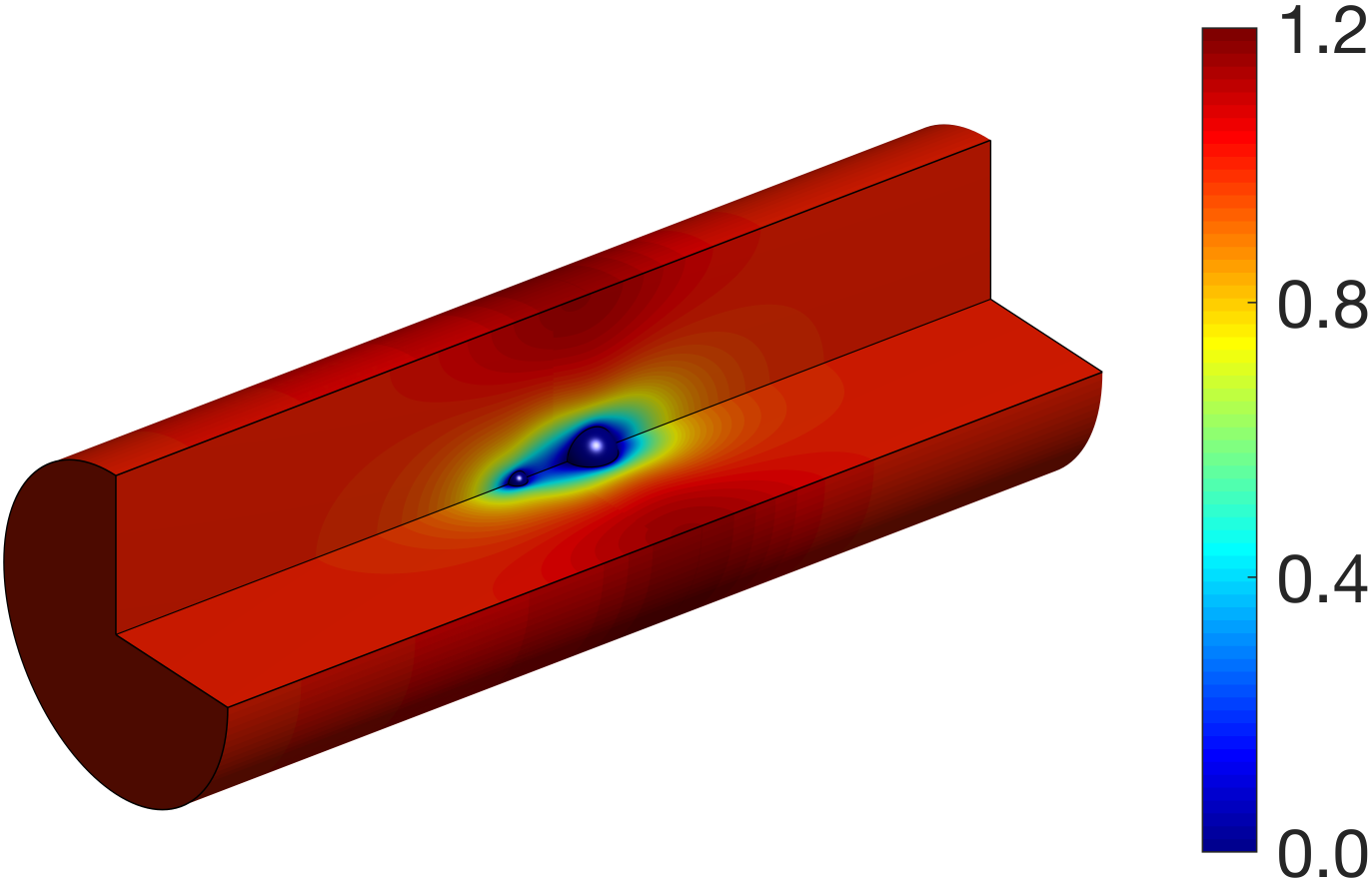}}
	\subfigure[$\mu_1=-1$, $\mu_2=-3$ \label{fig:swimmerCofigurations1P}]{\includegraphics[width=0.32\textwidth]{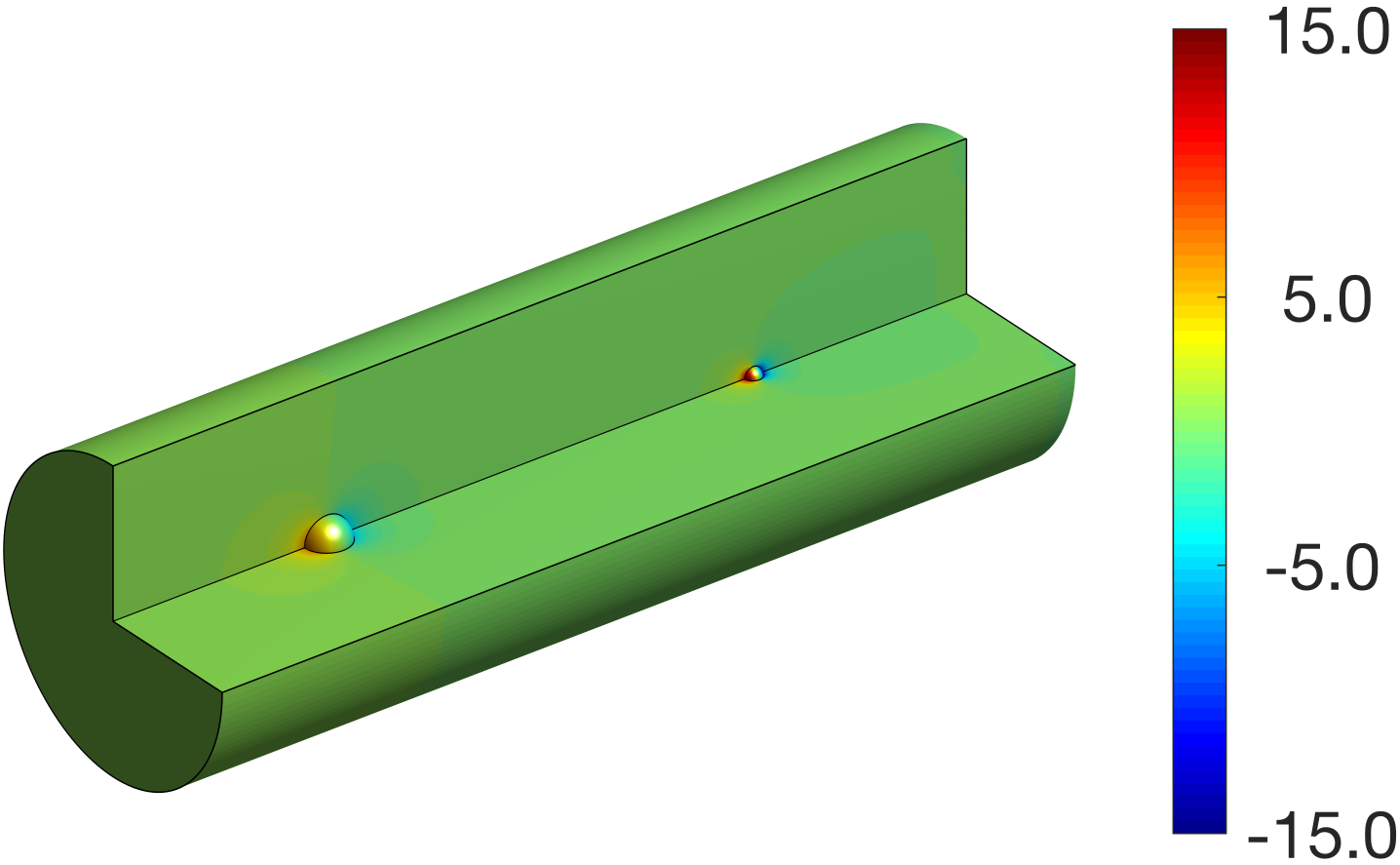}}
	\subfigure[$\mu_1=0$, $\mu_2=0$   \label{fig:swimmerCofigurations2P}]{\includegraphics[width=0.32\textwidth]{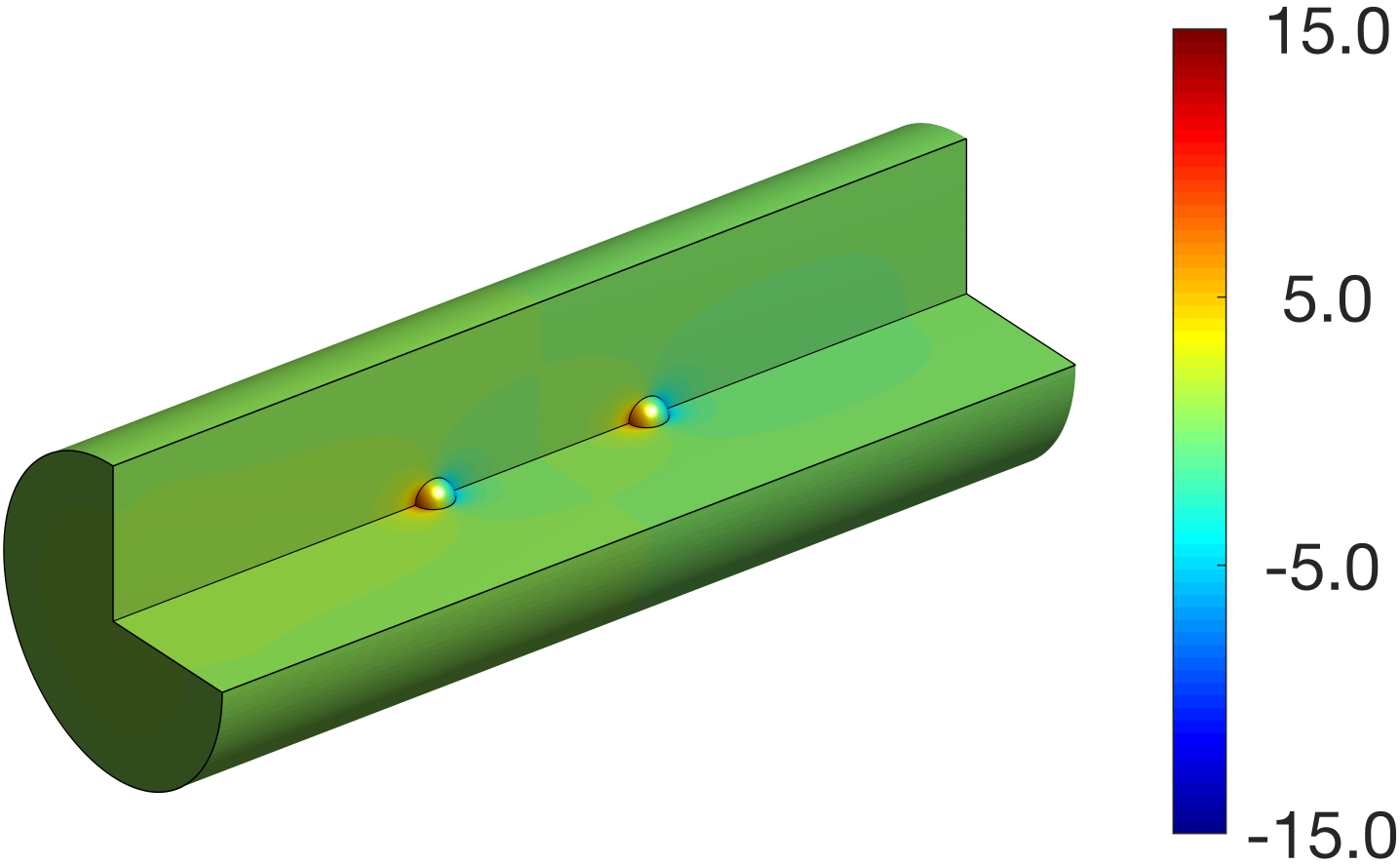}}
	\subfigure[$\mu_1=1$, $\mu_2=2$   \label{fig:swimmerCofigurations3P}]{\includegraphics[width=0.32\textwidth]{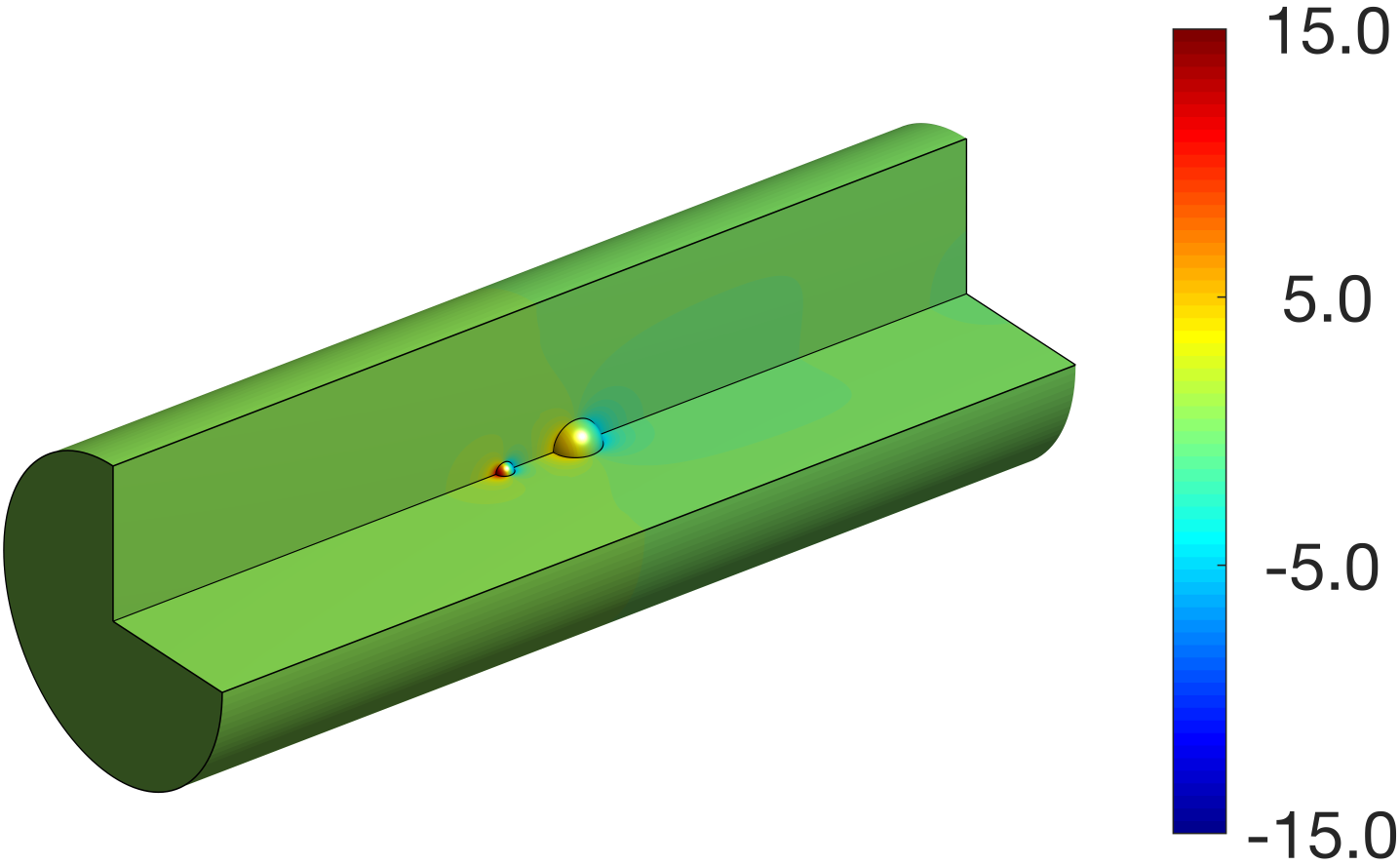}}
	\caption{Axisymmetric flow around two micro-swimmers: Velocity (top) and pressure (bottom) fields for three different geometric configurations.}
	\label{fig:swimmerCofigurations}
\end{figure}

To illustrate the variation in the geometry induced by the parameters as well as the different flow features that are induced by the geometric changes, figure~\ref{fig:swimmerCofigurations} shows the magnitude of the velocity and the pressure fields in the three dimensional domain for three different configurations.
The first configuration, shown in figures~\ref{fig:swimmerCofigurations1V} and ~\ref{fig:swimmerCofigurations1P}, corresponds to the case where the distance between the spheres is maximum and the sphere closer to the inflow boundary has maximum radius. The opposite scenario, with the distance between spheres is minimum and the sphere closer to the inflow boundary has minimum radius in shown in figures~\ref{fig:swimmerCofigurations3V} and~\ref{fig:swimmerCofigurations3P}. Finally, the configuration displayed in figures~\ref{fig:swimmerCofigurations2V} and ~\ref{fig:swimmerCofigurations2P} corresponds to the case when the distance between the spheres is half the maximum value and the radius of both spheres is the same.

\begin{figure}[!tb]
	\centering
	\subfigure[First sphere]{\includegraphics[width=0.45\textwidth]{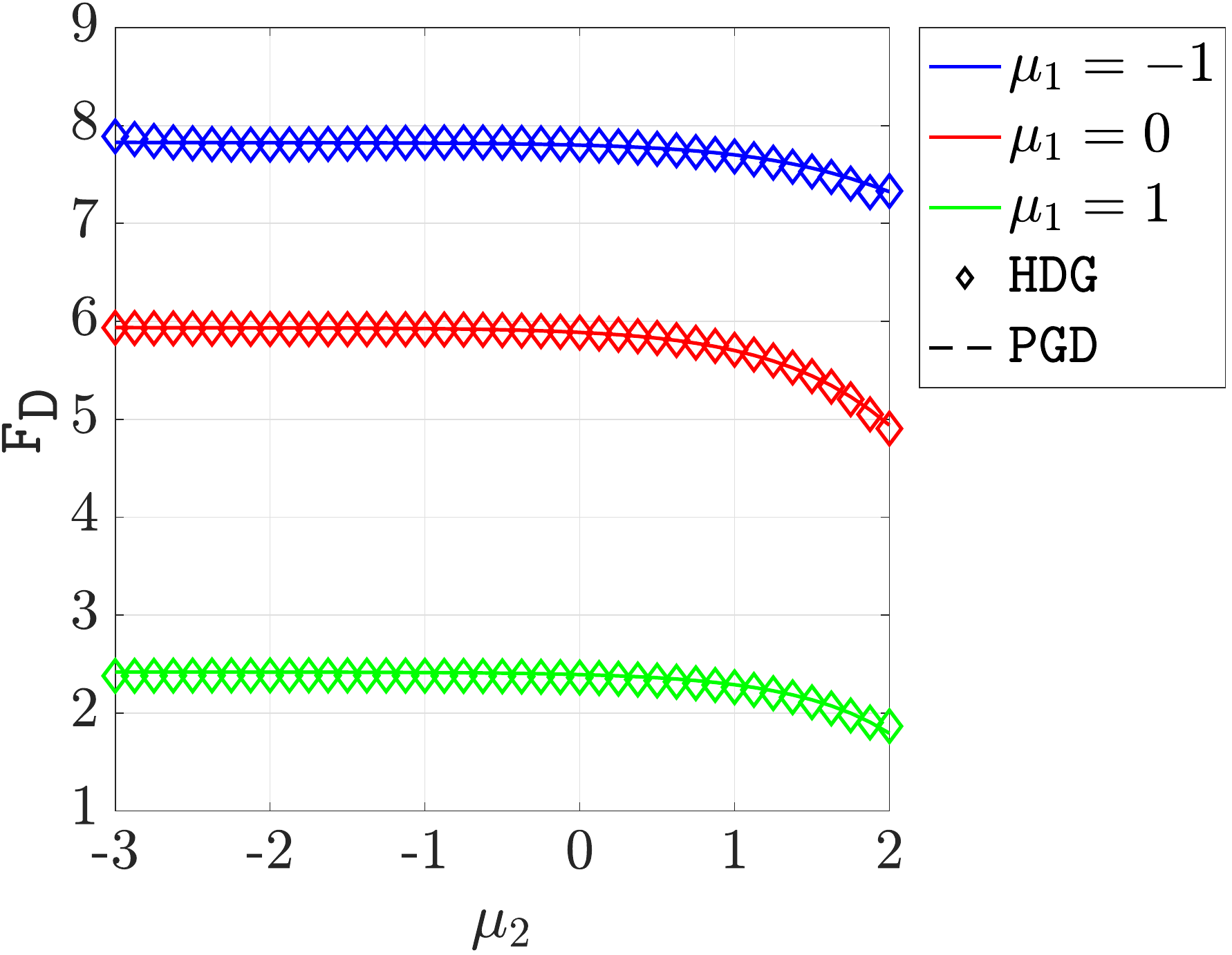}}
	\subfigure[Second sphere]{\includegraphics[width=0.45\textwidth]{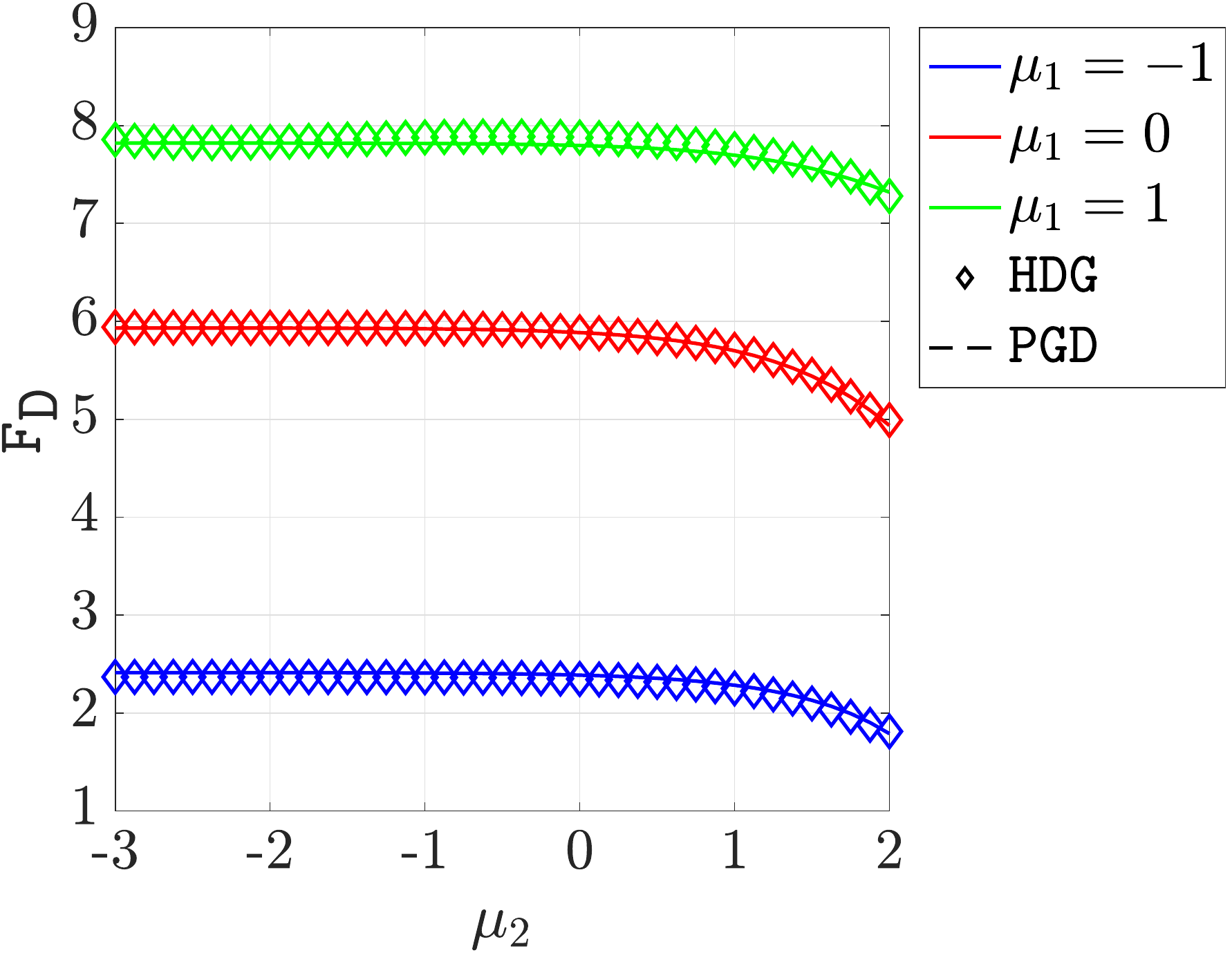}}
	\caption{Axisymmetric flow around two micro-swimmers: Comparison of the drag computed on the first and second sphere with the proposed HDG-PGD approach against a reference solution for different configurations.}
	\label{fig:swimmerDragComparison}
\end{figure}

To analyse the accuracy of the proposed approach, figure~\ref{fig:swimmerDragComparison} compares the drag force on the two spheres as a function of the $\mu_2$, controlling the distance between the spheres, and for three different configurations of the $\mu_1$, controlling the radius of both spheres.
The results obtained with the HDG-PGD approach are compared to the results of the standard HDG method on a reference mesh. Both solutions show an excellent agreement in all cases. 

\begin{figure}[!tb]
	\centering
	\subfigure[First sphere]{\includegraphics[width=0.32\textwidth]{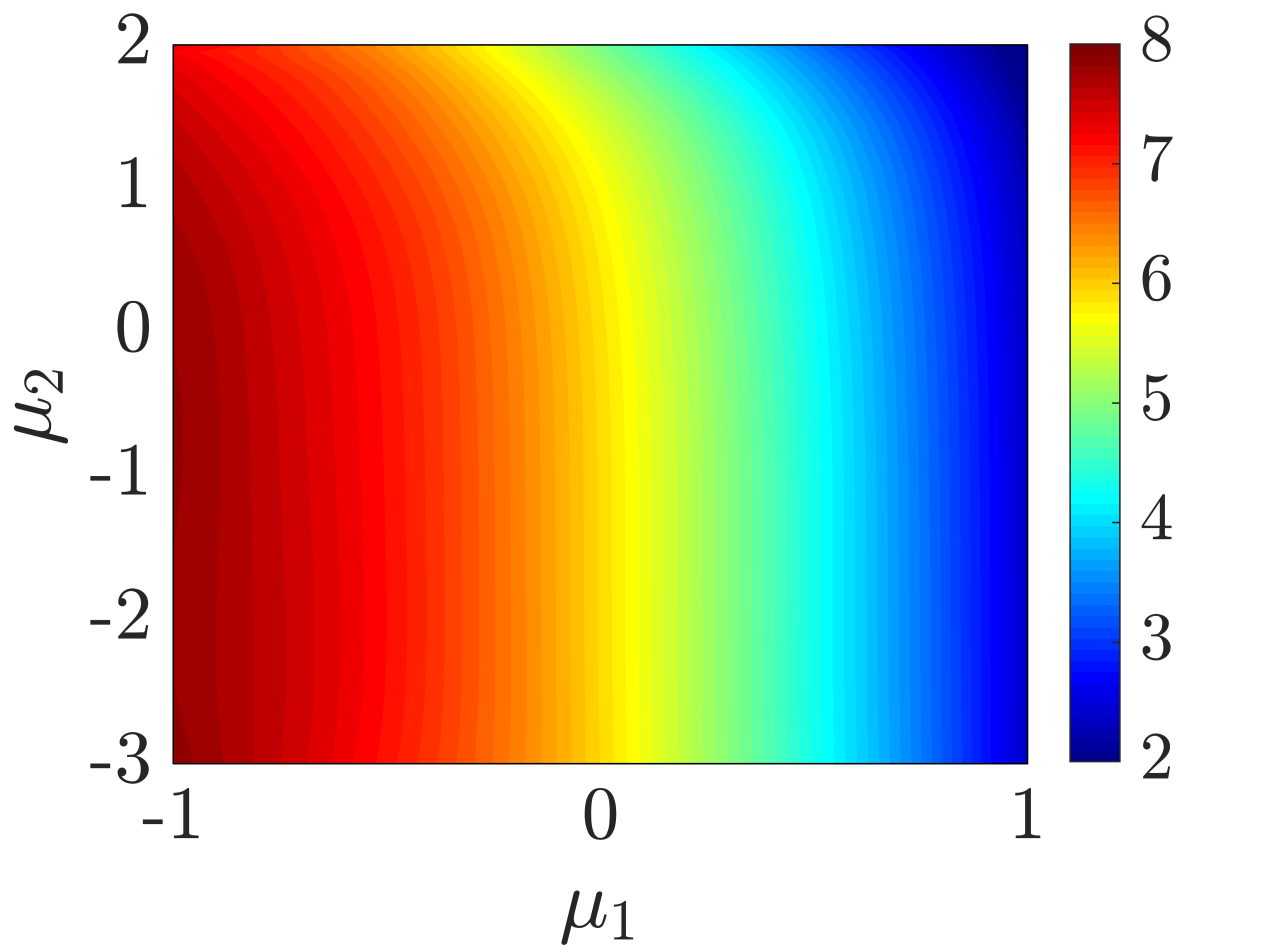}}
	\subfigure[Second sphere]{\includegraphics[width=0.32\textwidth]{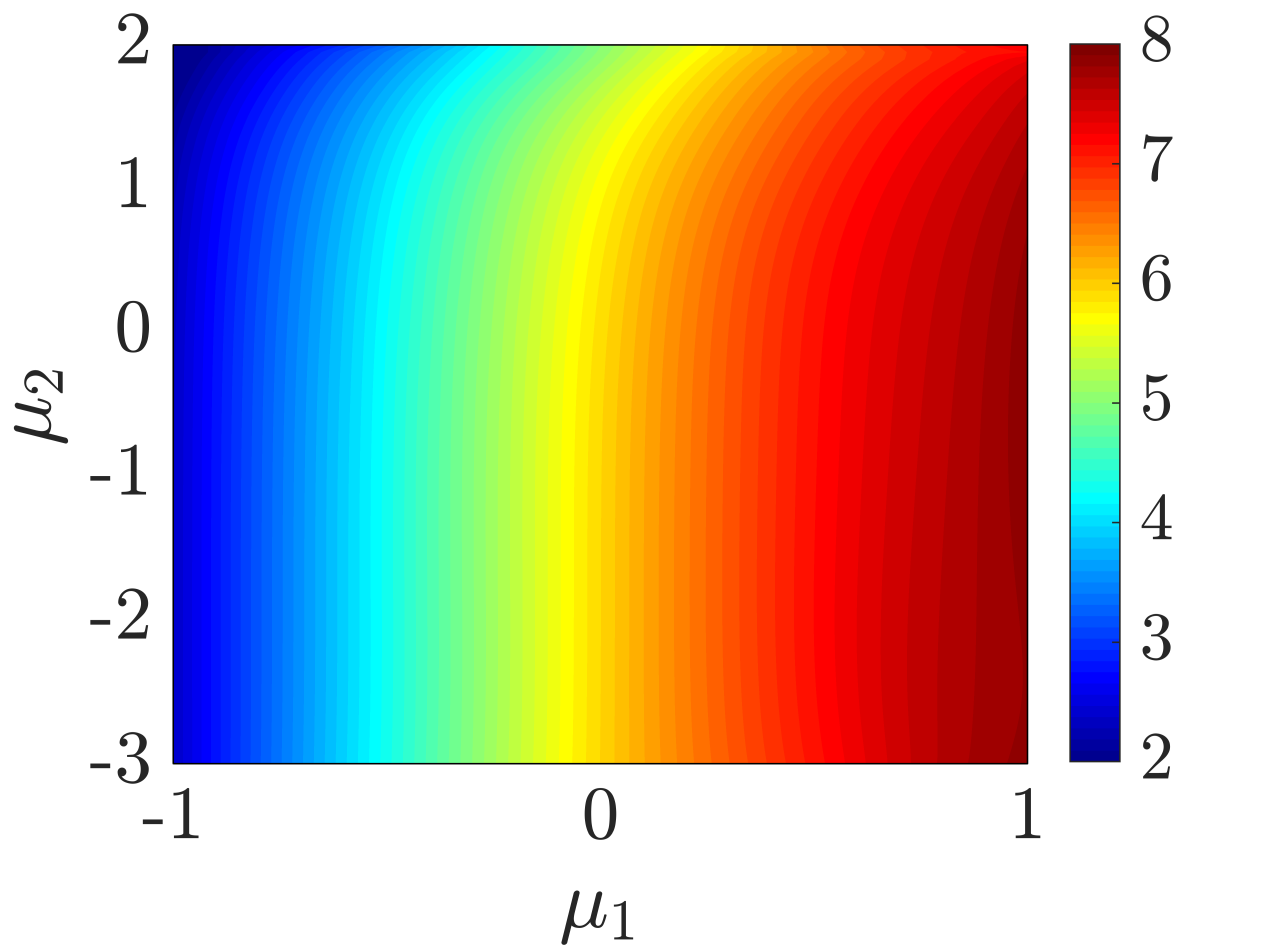}}
	\subfigure[Total]{\includegraphics[width=0.32\textwidth]{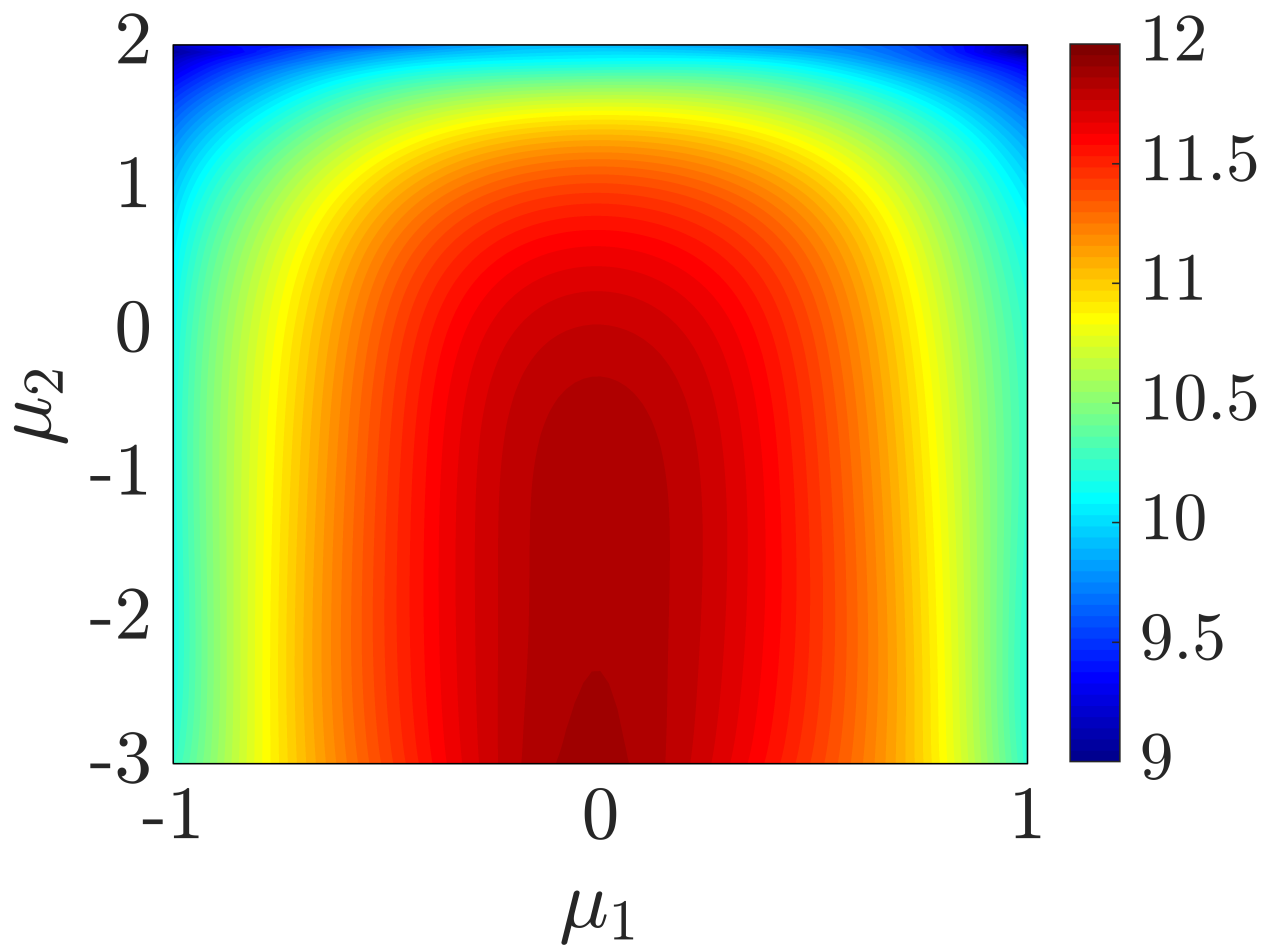}}
	\caption{Axisymmetric flow around two micro-swimmers: Drag force on the individual spheres and the total drag over the two spheres.}
	\label{fig:swimmerDrag}
\end{figure}

Finally, to stress the potential of the proposed approach, figure~\ref{fig:swimmerDrag} shows the drag force on the two spheres and the total drag as a function of both geometric parameters.
This figure shows that generalised solution computed with the HDG-PGD approach can be used to rapidly explore the whole space of parameters and used to find optimal strokes, of interest in many applications~\cite{alouges2009optimal}.

\subsection{Stokes flow around a sphere in a corrugated channel} \label{sc:channel}

The last example, inspired from the studies in~\cite{yang2017hydrodynamic,vieira2020second}, considers the flow past a sphere placed in a corrugated channel. 
The corrugated channel has a height of 1$\mu$m and the undulatory profile is defined by the expression
\begin{equation} \label{eq:geoDefChannel}
y = 
\begin{cases}
\frac{1}{2} \left( f_\omega + f_n \right) + \frac{1}{2} \left( f_\omega - f_n \right) \cos \left( \frac{16\pi x}{7L} \right) & \text{if} \quad \abs{x} < \frac{7}{16}L, \\
f_n & \text{if} \quad \frac{7}{16}L \leq \abs{x} \leq \frac{1}{2}L,
\end{cases}
\end{equation}
where $L=12.5\mu$m, $f_\omega=2\mu$m and the value of $f_n$ controls the oscillation of the boundary. A sphere of radius $R$, centred at the origin, is placed inside the corrugated channel. 

A Dirichlet boundary condition is imposed at one end of the channel, given by $\bu_D(\bm{x}) = \left\{ 64(x_2^2-1/4)(x_3^2-1/4), 0, 0 \right\}^T$, and a homogeneous Neumann boundary condition is imposed at the other end. A homogeneous Dirichlet boundary condition is on the rest of the boundary of the domain, corresponding to material walls. 

\begin{figure}[!tb]
	\centering
	\subfigure[Geometry \label{fig:channelGeo}]{\includegraphics[width=0.49\textwidth]{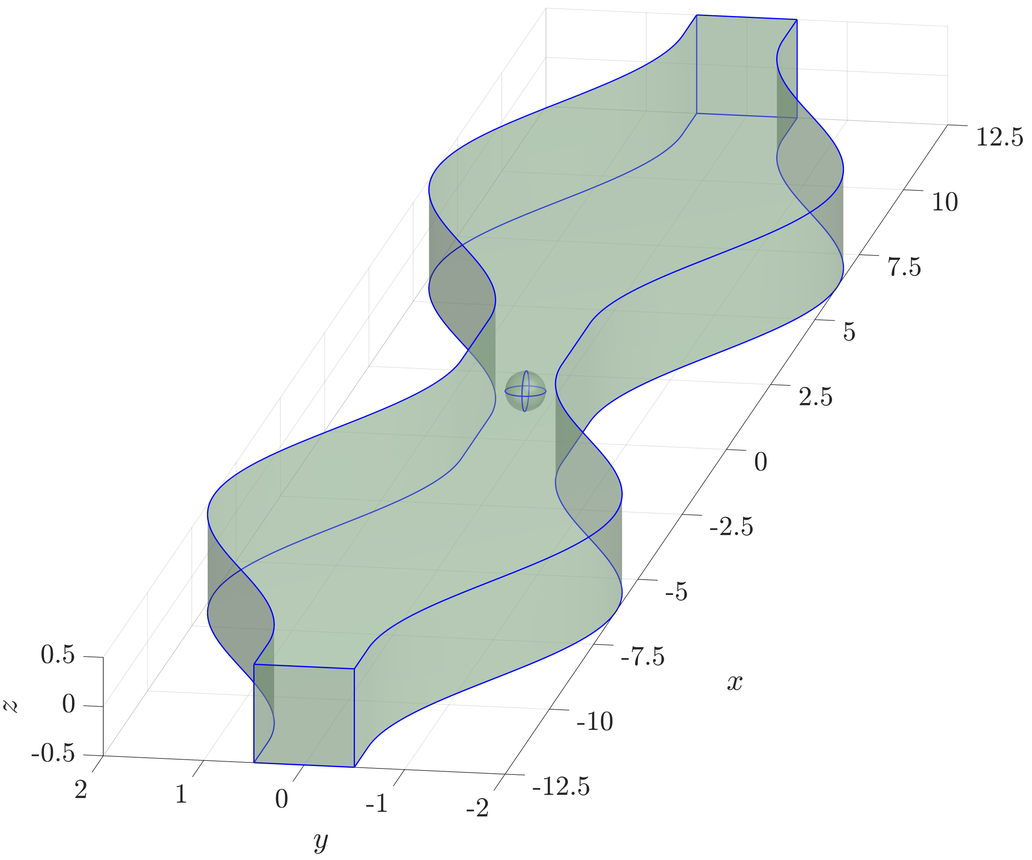}}
	\subfigure[Mesh		\label{fig:channelMesh}]{\includegraphics[width=0.49\textwidth]{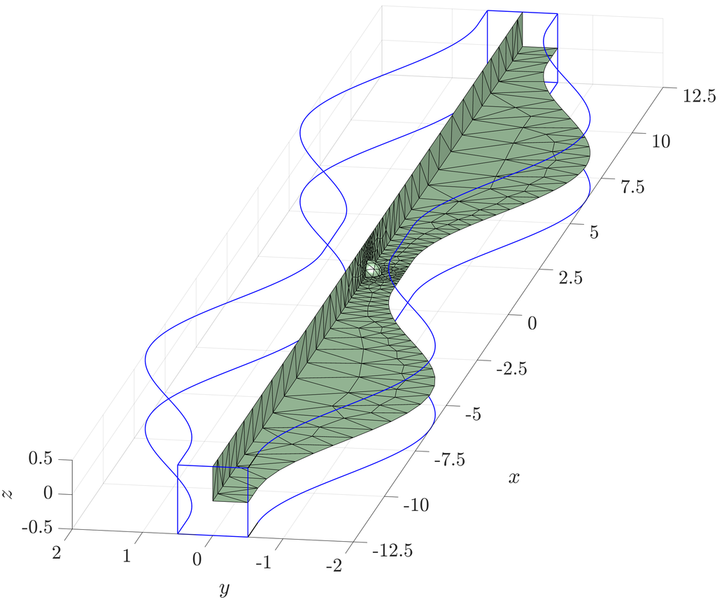}}
	\caption{Flow around a sphere in a corrugated channel: Geometry of the domain and computational mesh of a quarter of the domain.}
	\label{fig:channelGeoMesh}
\end{figure}

To demonstrate the applicability and potential of the proposed methodology in three dimensions, two geometric parameters are considered. The first parameter $\mu_1 \in [-1,1]$ is used to control the radius of the sphere, defined as $R(\mu_1) = (\mu_1+2)/10$. The second parameter $\mu_2 \in [0,2]$ controls the amplitude of the corrugated channel, given by $f_n = 1/2 + \mu_2$. The geometry of the reference domain, corresponding to $\mu_1=\mu_2=0$, is shown in Figure~\ref{fig:channelGeo}. 
Exploiting the symmetry of the problem, a mesh of a quarter of the domain is considered, with 2,191 tetrahedral elements, as depicted in  Figure~\ref{fig:channelMesh}. 

The geometric mapping used in this example is detailed in~\ref{sc:mappingCorrugated}.

\begin{figure}[!tb]
	\centering
	\subfigure[$m=1$]{\includegraphics[width=0.49\textwidth]{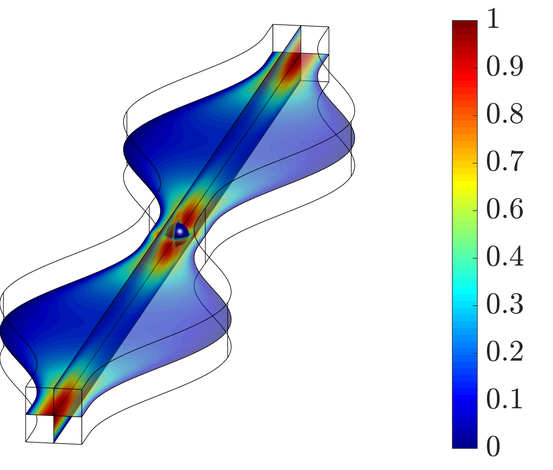}}
	\subfigure[$m=2$]{\includegraphics[width=0.49\textwidth]{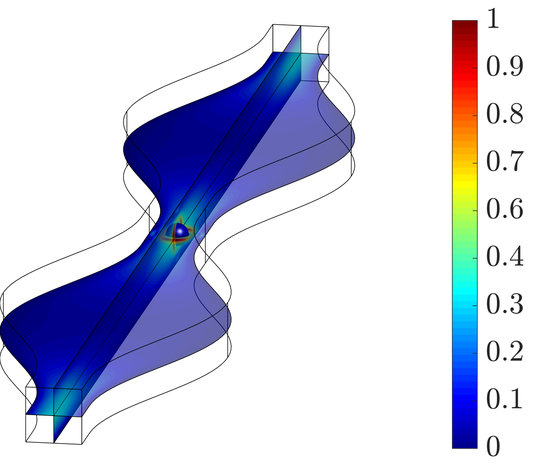}}
	\subfigure[$m=3$]{\includegraphics[width=0.49\textwidth]{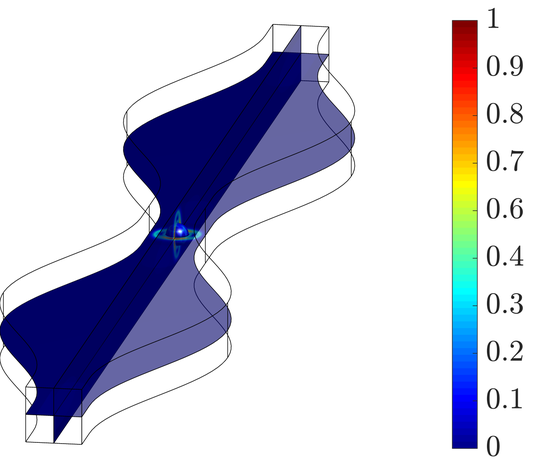}}
	\subfigure[$m=4$]{\includegraphics[width=0.49\textwidth]{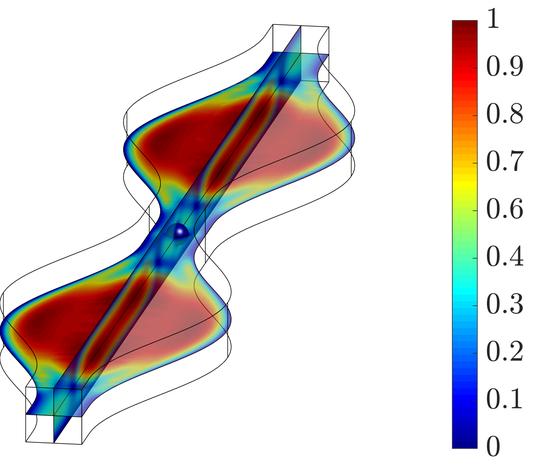}}
	\caption{Flow around a sphere in a corrugated channel: First four normalised spatial modes of the velocity field.}
	\label{fig:channelVModes}
\end{figure}

\begin{figure}[!tb]
	\centering
	\subfigure[$m=1$]{\includegraphics[width=0.49\textwidth]{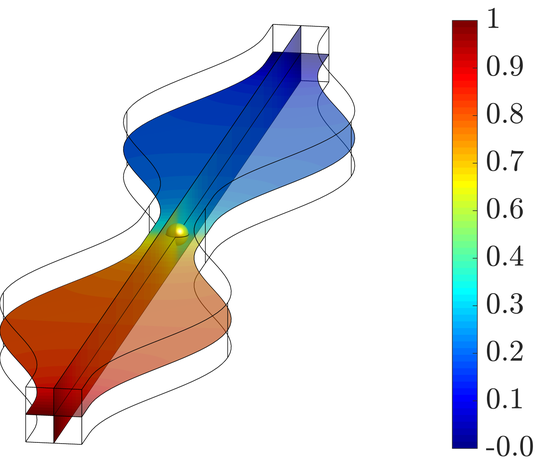}}
	\subfigure[$m=2$]{\includegraphics[width=0.49\textwidth]{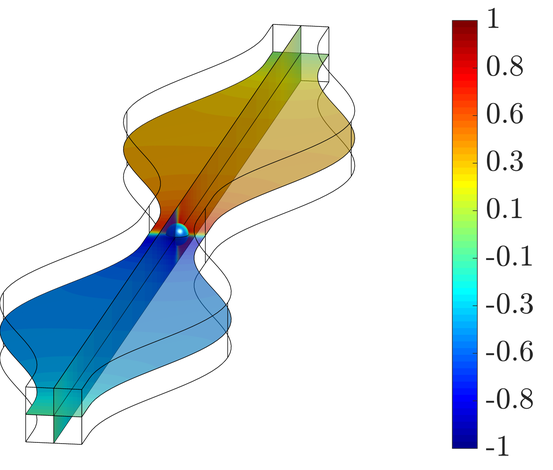}}
	\subfigure[$m=3$]{\includegraphics[width=0.49\textwidth]{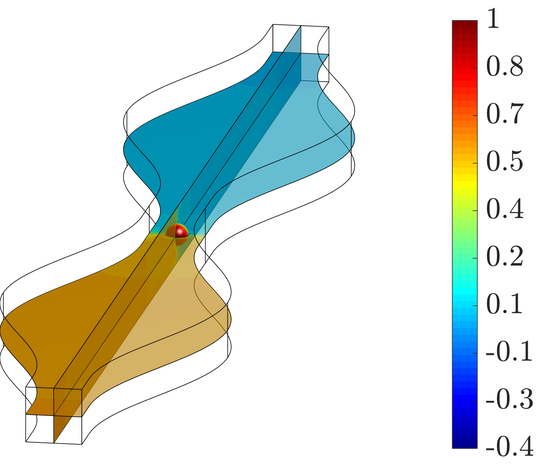}}
	\subfigure[$m=4$]{\includegraphics[width=0.49\textwidth]{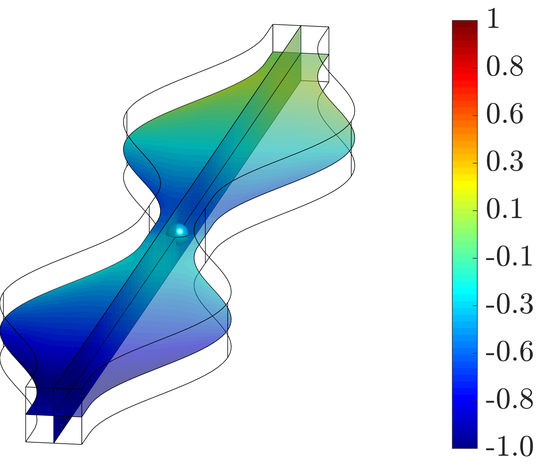}}
	\caption{Flow around a sphere in a corrugated channel: First four normalised spatial modes of the pressure field.}
	\label{fig:channelPModes}
\end{figure}

The first four spatial modes for the velocity and pressure computed with the proposed HDG-PGD are shown in figures~\ref{fig:channelVModes} and~\ref{fig:channelPModes}. 
The computation was performed using the mesh of figure~\ref{fig:channelGeo} with a degree of approximation $k=3$ for all the variables and with a mesh of 10,000 elements in each parametric dimension with also $k=3$. 

\begin{figure}[!tb]
	\centering
	\subfigure[\label{fig:channelParamModes1}]{\includegraphics[width=0.49\textwidth]{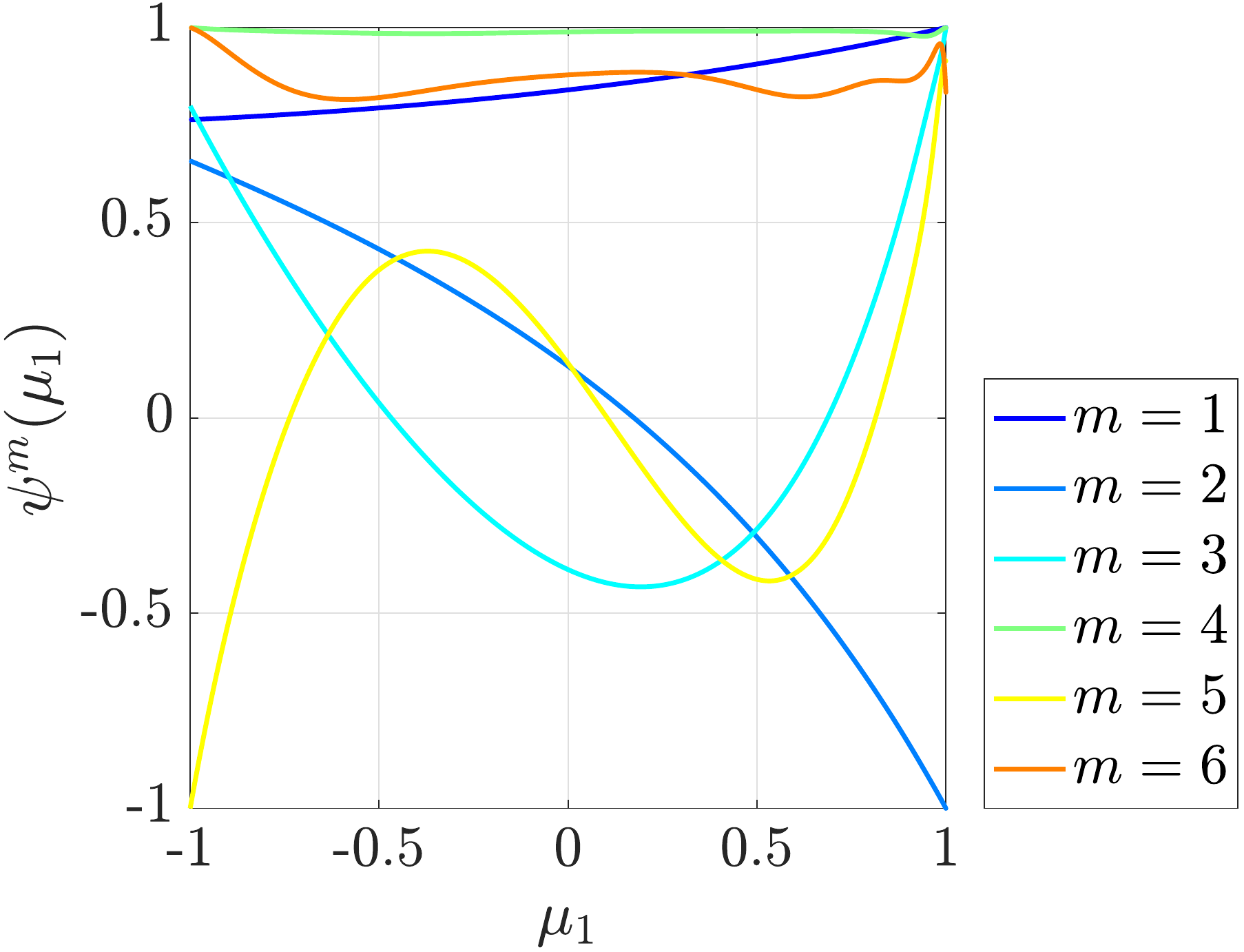}}
	\subfigure[\label{fig:channelParamModes2}]{\includegraphics[width=0.49\textwidth]{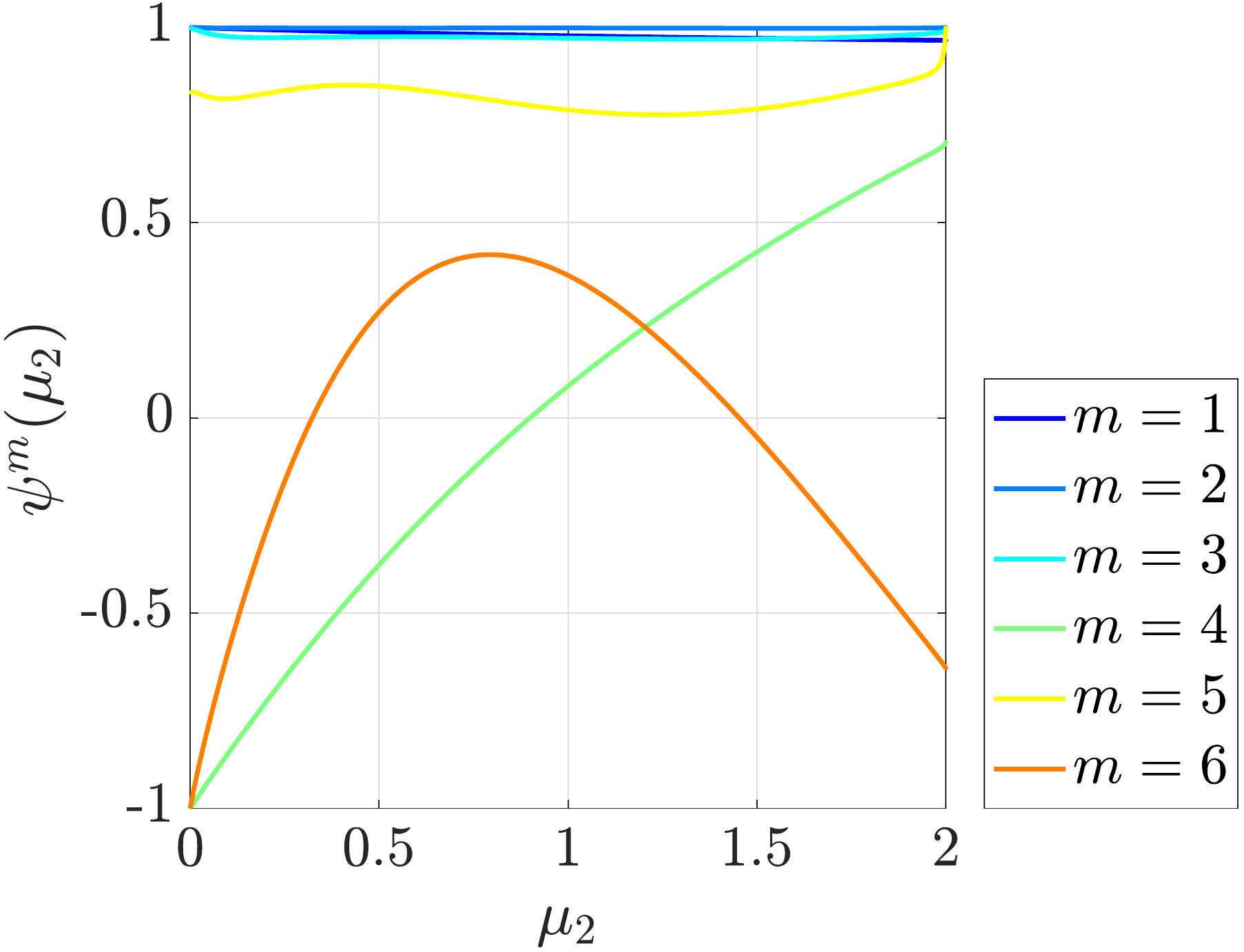}}
	\caption{Flow around a sphere in a corrugated channels: First six normalised parametric modes.}
	\label{fig:channelParamModes}
\end{figure}

Figure~\ref{fig:channelParamModes} shows the first six normalised parametric modes computed. 
Compared to previous examples, the results show that more modes have an influence over the whole range of parameters, illustrating the more complex nature of this three dimensional example.

\begin{figure}[!tb]
	\centering
	\includegraphics[width=0.45\textwidth]{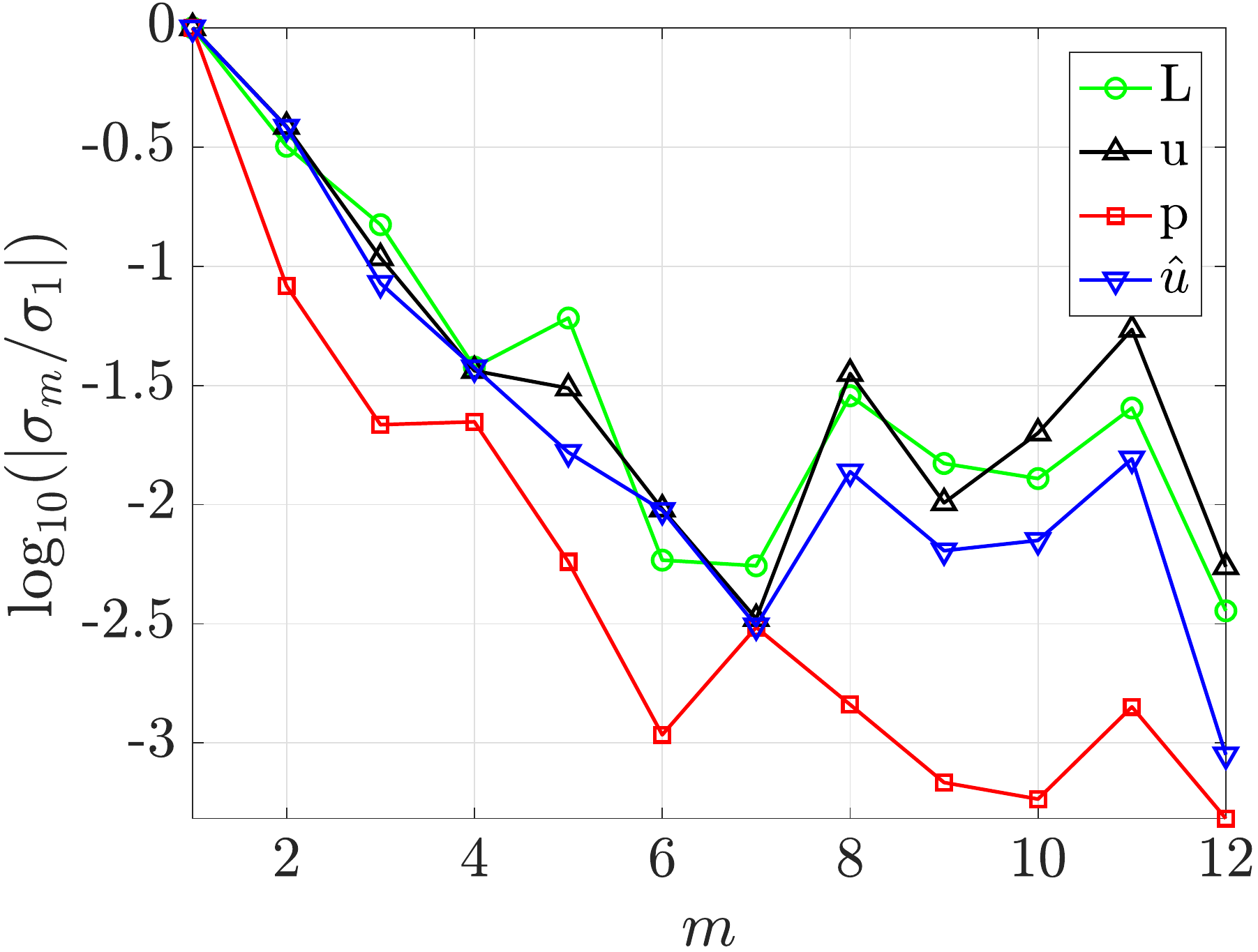}
	\caption{Flow around a sphere in a corrugated channel: Convergence of the mode amplitudes.}
	\label{fig:channelAmplitudes}
\end{figure}

The evolution of the relative amplitude of the modes is displayed in figure~\ref{fig:channelAmplitudes}. 
In this example, 12 modes are required to ensure the relative amplitude of the hybrid variable, used to check convergence, is below $10^{-3}$.  

\begin{figure}[!tb]
	\centering
	\subfigure[$\mu_1=-1$, $\mu_2=0$]{\includegraphics[width=0.32\textwidth]{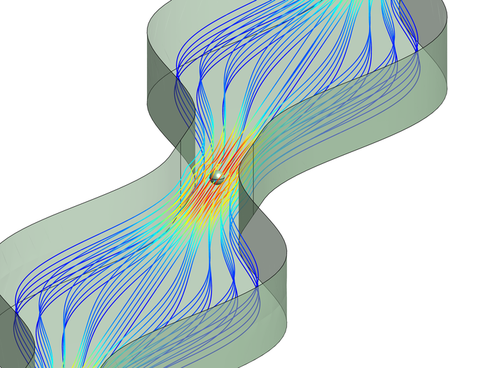}}
	\subfigure[$\mu_1= 0$, $\mu_2=1$]{\includegraphics[width=0.32\textwidth]{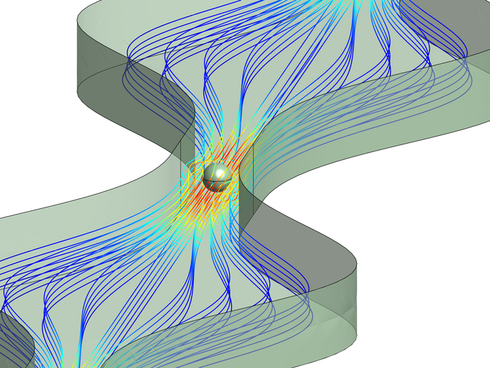}}
	\subfigure[$\mu_1= 1$, $\mu_2=2$]{\includegraphics[width=0.32\textwidth]{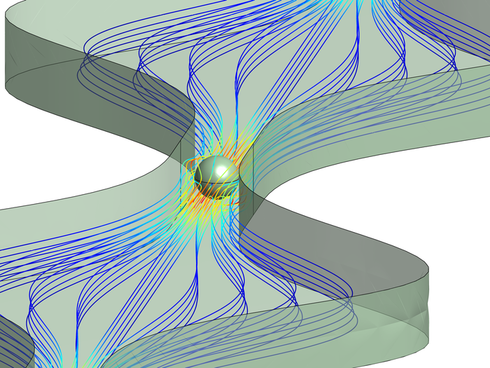}}
	\subfigure[$\mu_1=-1$, $\mu_2=0$]{\includegraphics[width=0.32\textwidth]{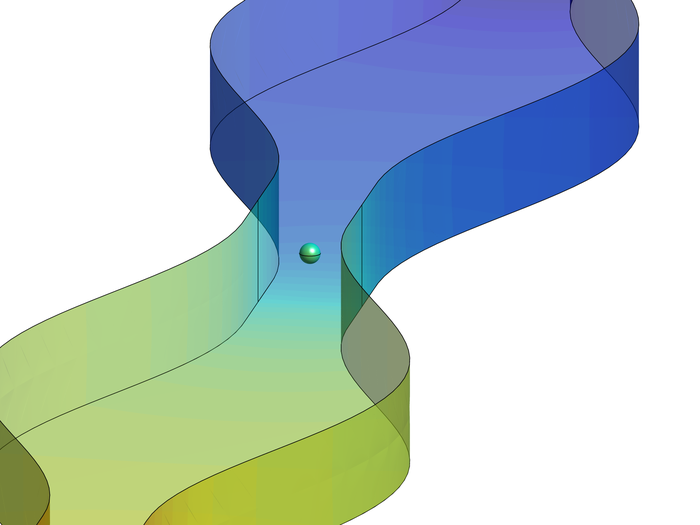}}
	\subfigure[$\mu_1= 0$, $\mu_2=1$]{\includegraphics[width=0.32\textwidth]{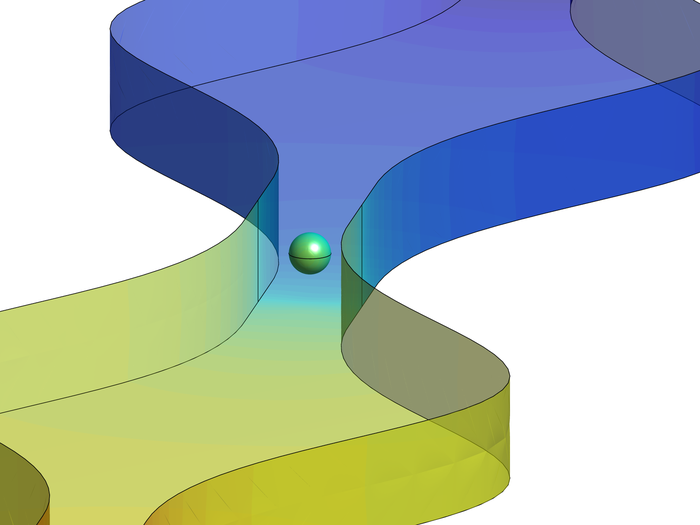}}
	\subfigure[$\mu_1= 1$, $\mu_2=2$]{\includegraphics[width=0.32\textwidth]{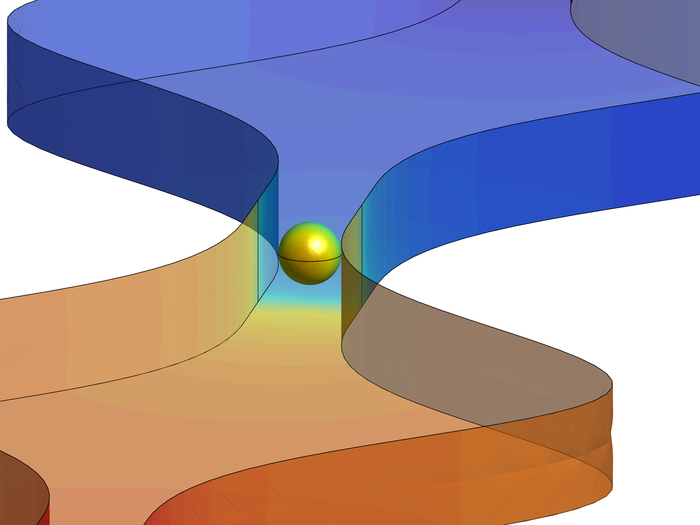}}
	\caption{Flow around a sphere in a corrugated channel: Velocity (top) and pressure (bottom) fields for three different geometric configurations.}
	\label{fig:channelCofigurations}
\end{figure}

Figure~\ref{fig:channelCofigurations} shows the magnitude of the velocity and the pressure fields in the channel for three different configurations.
The results illustrate the variation in the velocity and pressure fields as the amplitude of the channel and the radius of the sphere is increased.

\begin{figure}[!tb]
	\centering
	\subfigure[$\mu_1=-1$, $\mu_2=0$]{\includegraphics[width=0.32\textwidth]{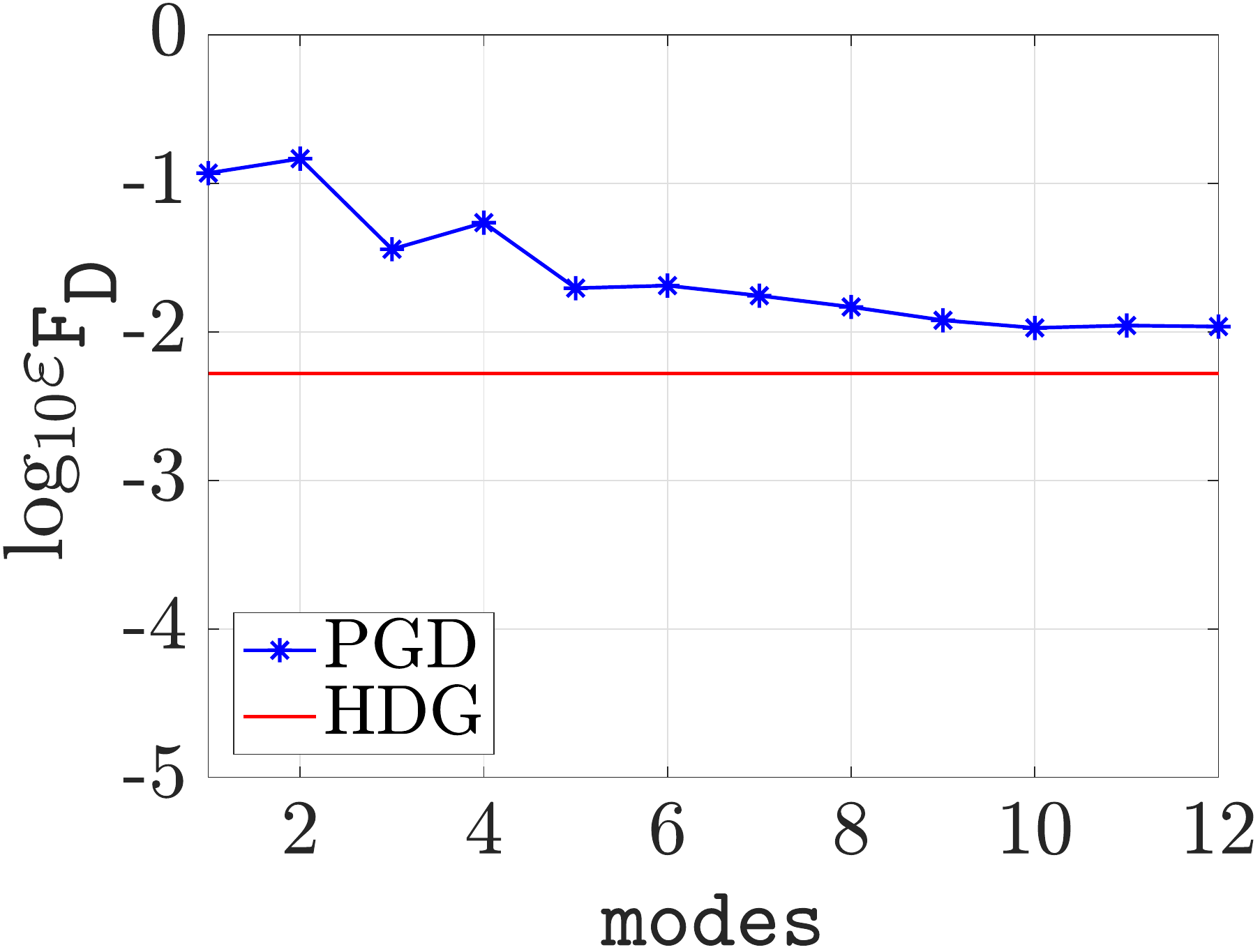}}
	\subfigure[$\mu_1= 0$, $\mu_2=1$]{\includegraphics[width=0.32\textwidth]{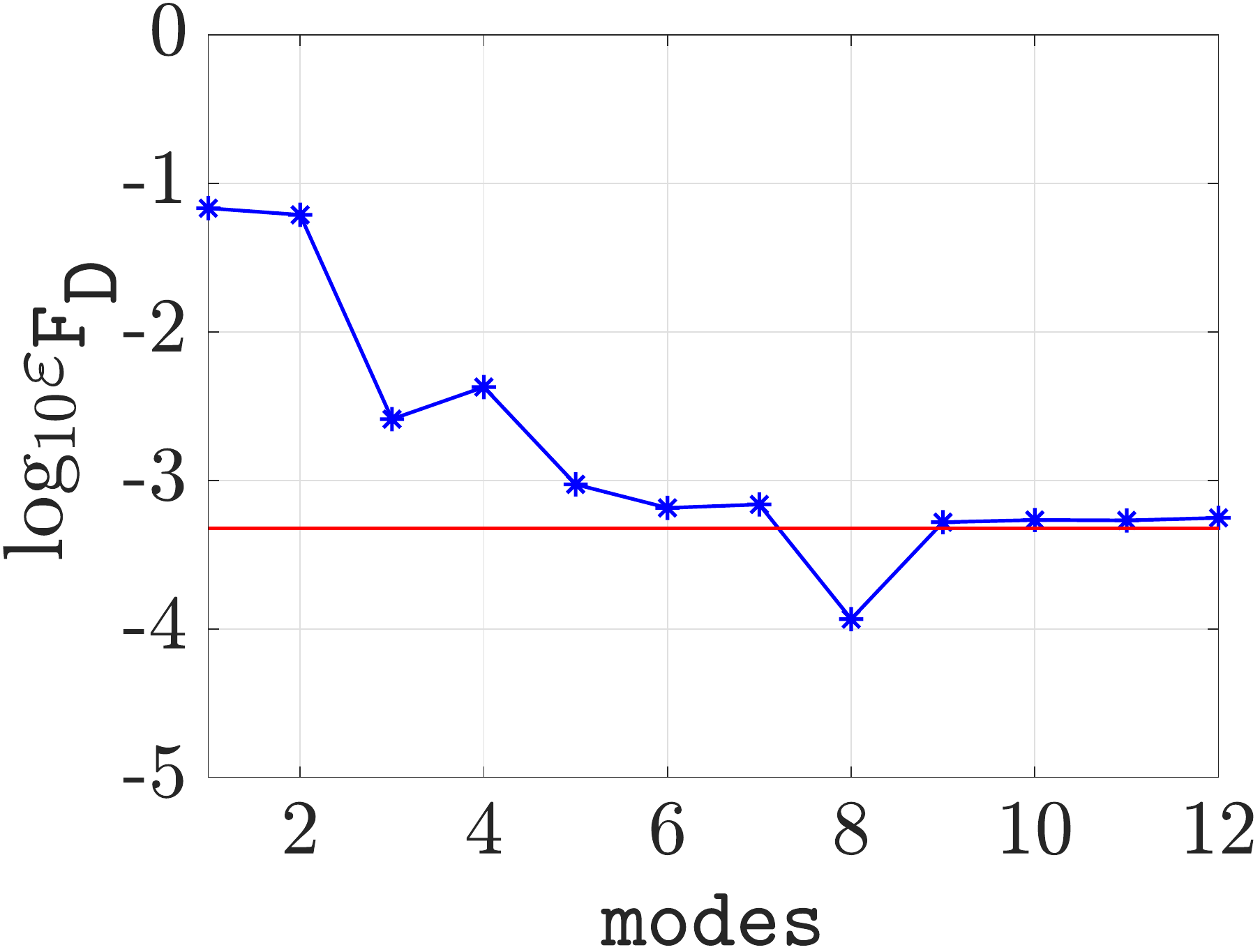}}
	\subfigure[$\mu_1= 1$, $\mu_2=2$]{\includegraphics[width=0.32\textwidth]{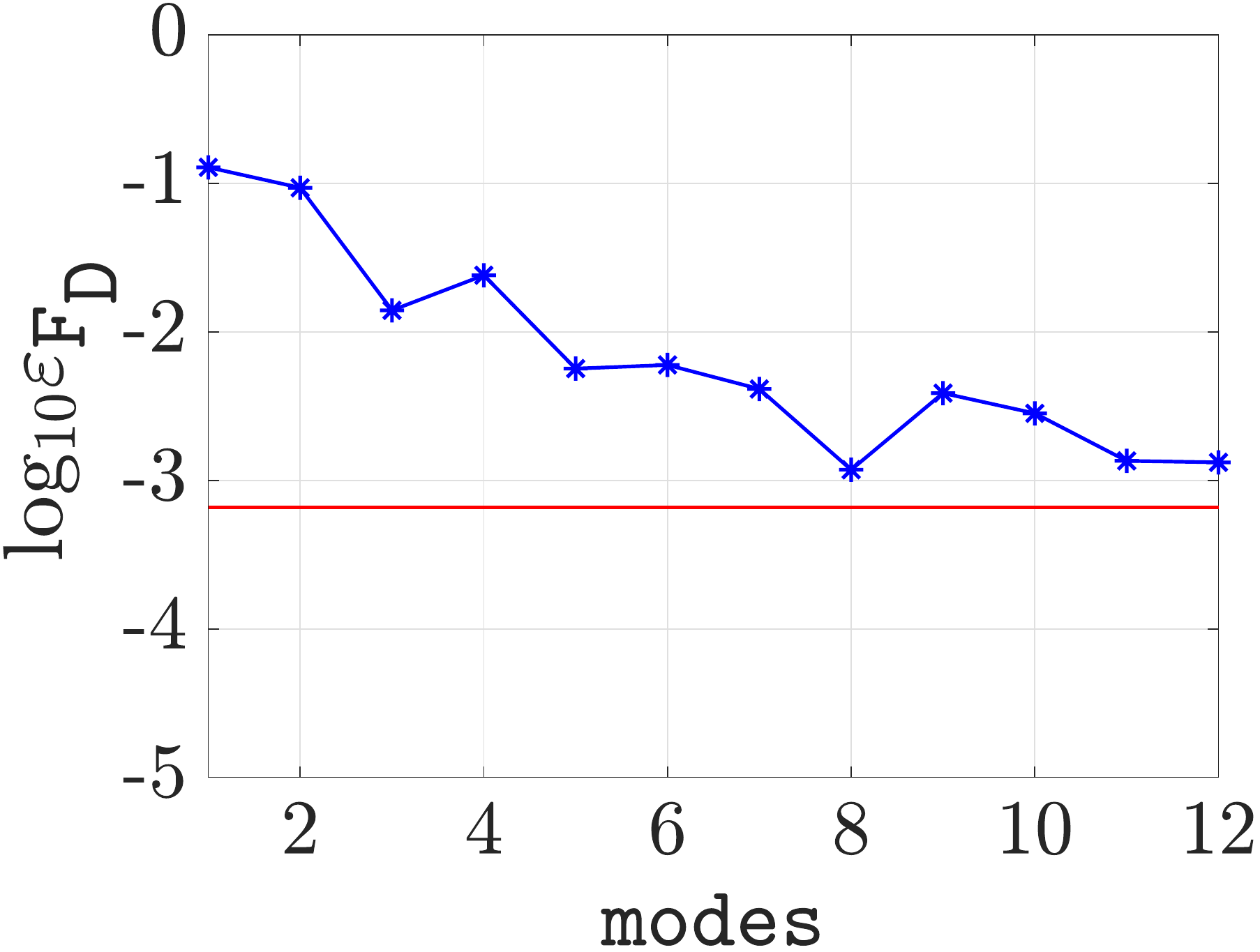}}
	\caption{Flow around a sphere in a corrugated channel: Evolution of the error on the drag force as the number of PGD modes is increased. The horizontal line denotes the reference error computed on a finer mesh with the standard HDG method.}
	\label{fig:channelDragError}
\end{figure}

To assess the accuracy of the computed generalised solution computed with the proposed approach, a reference solution is computed for the three configurations displayed in Figure~\ref{fig:channelCofigurations}. The reference solutions are computed on a much finer mesh with a standard HDG solver. As a quantity of interest, the drag on the sphere is measured. Figure~\ref{fig:channelDragError} shows the evolution of the error of the drag force as the number of PGD modes is increased.
To further analyse the accuracy of the computed generalised solution, the error of an HDG solution, computed in each configuration using the same spatial resolution as the one used in the HDG-PGD formulation is considered. The results show that the error of the HDG-PGD approach tends to the error of the HDG solution computed for each configuration, showing the ability of the proposed approach to accurately capture the solution for different geometric configurations.

\begin{figure}[!tb]
	\centering
	\subfigure[Drag]{\includegraphics[width=0.45\textwidth]{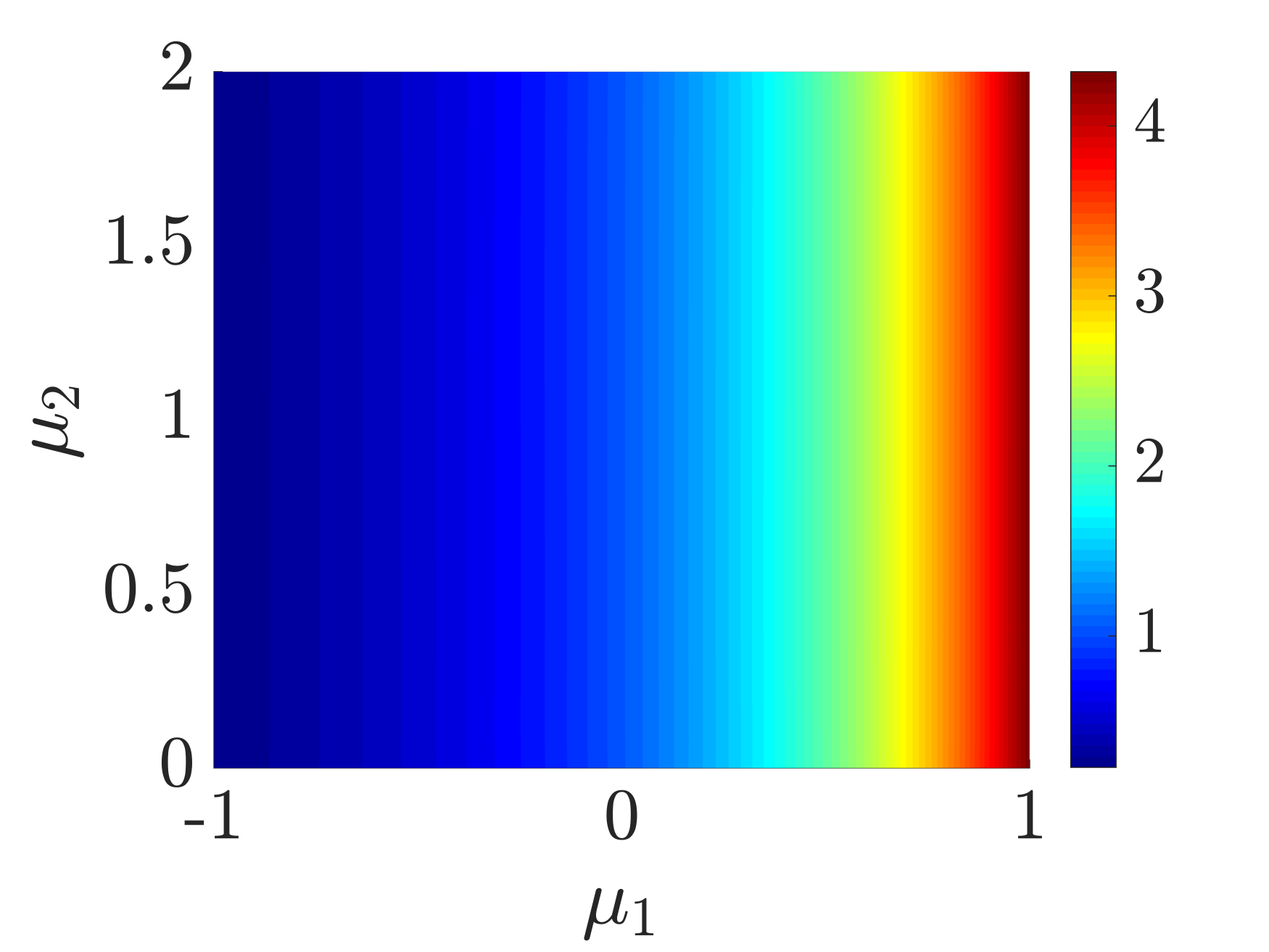}}
	\subfigure[Pressure drop]{\includegraphics[width=0.45\textwidth]{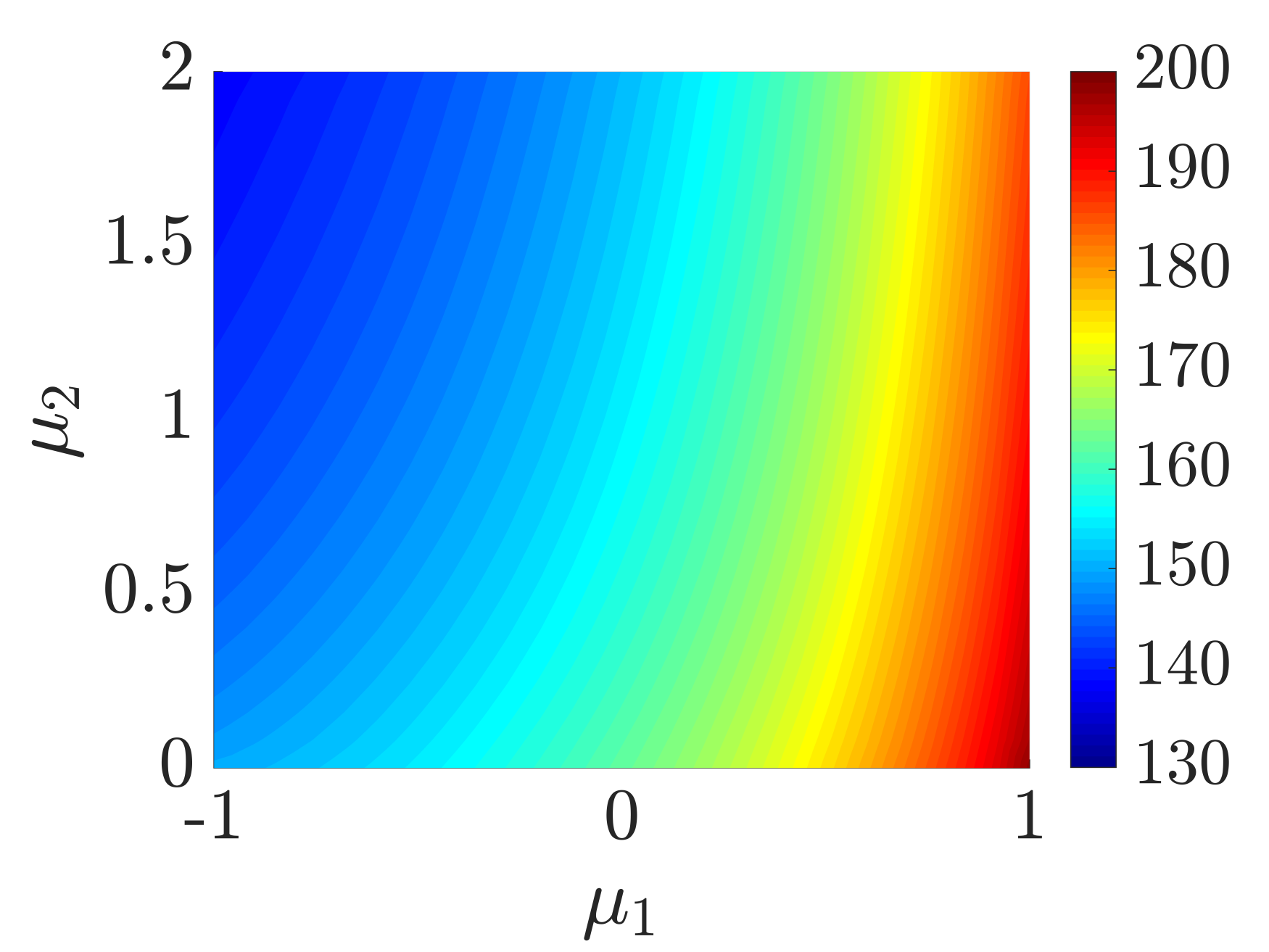}}
	\caption{Flow around a sphere in a corrugated channel: Drag force on the sphere and difference between the pressure at the inlet and the outlet.}
	\label{fig:channelQuantities}
\end{figure}

As mentioned in the previous example, the proposed approach provides a generalised solution that can be used to perform fast queries of different quantities of interest. To illustrate the potential of the developed HDG-PGD approach, Figure~\ref{fig:channelQuantities} shows the drag force on the sphere and the pressure drop, measured as the difference between the pressure at the inlet and outlet, as a function of the geometric parameters $\mu_1$ and $\mu_2$.
The results show that the drag force is not sensitive to the variation of the amplitude of the channel oscillation but very dependent on the radius of the sphere. In contrast, the pressure drop shows a dependency on both geometric parameters.

\section{Concluding remarks} \label{sc:conclusions}

A reduced order model approach based on the PGD and the HDG methods is being presented for the solution of geometrically parametrised Stokes flow problems. The mixed formulation characteristic of HDG methods is beneficial in the PGD context as all the terms of the weak formulation can be written in a separated form, without using to the memory intensive high-order PGD projection. The use of the HDG formulation also enables the use of equal order of approximation for all the variables circumventing the LBB condition. This is advantageous in the context of geometrically parametrised problems in complex domains as it enables the use of standard isoparametric formulations. In addition, the use of a DG formulation implies that no special treatment of the Dirichlet boundary conditions is required.

The optimal approximation properties of the proposed approach have been validated numerically using two and three dimensional test cases. In addition, the ability of the proposed approach to compute generalised solutions involving geometric parameters has been illustrated for problems relevant to the microfluidics community. The examples involve geometric parameters that involve substantial changes of the geometry and induce important changes in the flow features and the relevant quantities of interest.

\section*{Acknowledgements}

This work was partially supported by the European Union's Horizon 2020 research and innovation programme under the Marie Sk\l odowska-Curie Actions (Grant number: 675919) that financed the Ph.D. fellowship of L.B. and by the Spanish Ministry of Economy and Competitiveness (Grant number: DPI2017-85139-C2-2-R). M.G. and A.H. are also grateful for the support provided by the Spanish Ministry of Economy and Competitiveness through the Severo Ochoa programme for centres of excellence in RTD (Grant number: CEX2018-000797-S) and the Generalitat de Catalunya (Grant number: 2017-SGR-1278). R.S. also acknowledges the support of the Engineering and Physical Sciences Research Council (Grant number: EP/P033997/1).

\bibliographystyle{abbrv}
\bibliography{Ref-PGD-Geo}

\appendix

\section{Bilinear forms of the HDG-PGD weak formulation}  \label{sc:bilinearPGD}

The bilinear forms introduced in the spatial iteration are given by
\begin{equation} \label{eq:bilinearPGDSpatialLocal}
	\begin{aligned}	
	\mathcal{A}^k_{LL}(\de \fL, \fL)	 & := - \intE{\de \fL}{\nu^{-1} D^k \fL }, 	&	
	\mathcal{A}^k_{Lu}(\de \fL, \fU)	 & :=  \intE{\bA^k \DivX \de \fL}{ \fU },  		\\
	\mathcal{A}^k_{L\hu}(\de \fL, \fHU)& := \intBNoD{\bA^k \bn \cdot \de \fL}{ \fHU }, &
	\mathcal{A}^k_{uL}(\de \fU, \fL)	 & :=  \intE{\de \fU}{\bA^k \DivX \fL },  		\\
	\mathcal{A}_{uu}(\de \fU, \fU)	 & :=  \intBE{\de \fU}{\btau \fU },				&
	\mathcal{A}^k_{up}(\de \fU, \fP)	 & :=  \intE{\de \fU}{\bA^k \GradX \fP },  		\\
	\mathcal{A}_{u\hu}(\de \fU, \fHU)	 & := \intBNoD{\de \fU}{\btau \fHU }			&
	\mathcal{A}^k_{pu}(\de \fP, \fU)	 & :=   \intE{\bA^k \GradX \de \fP}{\fU},  		\\
	\mathcal{A}^k_{p\hu}(\de \fP, \fHU)& :=  \intBNoD{\de \fP}{\fHU \cdot \bA^k \bn},  &
	\mathcal{A}_{\rho p}(\de\fR, \fP)  & :=  \intBE{\de\fR}{|\partial\Omega_e|^{-1} \fP },						\\
	\mathcal{A}_{\rho \rho}(\de\fR, \fR) & :=  \de\fR \,\fR, 
\end{aligned}	
\end{equation}
for the HDG local problems and by
\begin{equation} \label{eq:bilinearPGDSpatialGlobal}
\begin{aligned}	
	\mathcal{A}^k_{\hu L}(\de \fHU, \fL)	 	& := \intBNoDS{\de \fHU}{\bA^k \bn \cdot \fL } - \intBS{\de \fHU}{\bA^k \bn \cdot \fL \bE },  		\\
	\mathcal{A}_{\hu u}(\de \fHU, \fU)     	& :=  \intBNoDS{\de \fHU}{\btau \fU} - \intBS{\de \fHU}{(\btau \fU) \!\cdot\! \bE} , 	\\
	\mathcal{A}^k_{\hu p}(\de \fHU, \fP)   	& :=  \intBNoDS{\de \fHU}{\fP \bA^k \bn }, \\
	\mathcal{A}_{\hu \hu}(\de \fHU, \fHU)		& := -\intBNoDS{\de \fHU}{\btau \fHU}  + \intBS{\de \fHU}{(\btau \fHU) \!\cdot\! \bE} ,  \\
	\mathcal{A}^k_{\hu \hu}(\de \fHU, \fHU)	& := \intBS{\de \fHU}{\fHU \cdot \bA^k \bD},
\end{aligned}	
\end{equation}
for the HDG global problems.

In addition, the following bilinear forms are introduced in the parametric iteration
\begin{equation} \label{eq:bilinearPGDParam}
\begin{aligned}	
\mathcal{A}^k_{\theta}(\de \psi, \psi)	& :=  \intI{\de \psi}{\theta^k \psi},  		\\
\mathcal{A}^k_{\vartheta}(\de \psi, \psi)	& :=  \intI{\de \psi}{\vartheta^k \psi}, 		\\
\mathcal{A}(\de \psi, \psi)				& :=  \intI{\de \psi}{\psi}.	
\end{aligned}	
\end{equation}

\section{Linear forms of the HDG-PGD weak formulation}  \label{sc:linearPGD}

The linear forms introduced in the spatial and parametric iterations are given by
\begin{equation} \label{eq:linearLocalPGDSpatial}
\begin{aligned}	
\mathcal{R}_L^m (\de \fL \psi)	:= & \sum_{k=1}^{\nadj} \sum_{l=1}^{\nDir} \intBD{\bA^k \bn \cdot \de \fL}{\bgD^l} \mathcal{A}^k_{\vartheta}(\psi, \lambda_D^l) 	\\
& {-} \sum_{i=1}^{m} \sum_{k=1}^{\ndet} \mathcal{A}^k_{LL}(\de \fL, \ampL^i \fL^i) \mathcal{A}^k_{\theta}(\psi, \psi^i)   \\
& {-} \sum_{i=1}^{m} \sum_{k=1}^{\nadj} \left\{ \mathcal{A}^k_{Lu}(\de \fL, \ampU^i \fU^i) - \mathcal{A}^k_{L\hu}(\de \fL, \ampHU^i \fHU^i)\right\} \mathcal{A}^k_{\vartheta}(\psi, \psi^i) 
\\
\mathcal{R}_u^m (\de \fU \psi)	 := & \sum_{k=1}^{\ndet} \sum_{l=1}^{\nSou} \intE{\de \fU}{D^k \bgS^l} \mathcal{A}^k_{\theta}(\psi, \lambda_S^l) 	\\
& {+}  \sum_{l=1}^{\nDir} \intBD{\de \fU}{\btau \bgD^l} \mathcal{A}(\psi, \lambda_D^l)  \\
& {-}  \sum_{i=1}^{m} \sum_{k=1}^{\nadj} \left\{ \mathcal{A}^k_{uL}(\de \fU, \ampL^i \fL^i) + \mathcal{A}^k_{up}(\de \fU, \ampP^i \fP^i) \right\} \mathcal{A}^k_{\vartheta}(\psi, \psi^i) \\ 
& {-}  \sum_{i=1}^{m} \left\{ \mathcal{A}_{uu}(\de \fU, \ampU^i \fU^i) - \mathcal{A}_{u\hu}(\de \fU, \ampHU^i \fHU^i) \right\} \mathcal{A}(\psi, \psi^i)
\\
\mathcal{R}_p^m (\de \fP \psi)	 := & \sum_{k=1}^{\nadj} \sum_{l=1}^{\nDir} \intBD{\de \fP}{\bgD^l \cdot \bA^k \bn} \mathcal{A}^k_{\vartheta}(\psi, \lambda_D^l)	\\
& {-}  \sum_{i=1}^{m} \sum_{k=1}^{\nadj} \left\{ \mathcal{A}^k_{pu}(\de \fP, \ampU^i \fU^i) - \mathcal{A}^k_{p\hu}(\de \fP, \ampHU^i \fHU^i) \right\} \mathcal{A}^k_{\vartheta}(\psi, \psi^i)
\\
\mathcal{R}_{\overline{p}}^m (\de\fR\psi) := & -\sum_{i=1}^{m} \left\{ \mathcal{A}_{\rho p}(\de\fR, \ampP^i \fP^i) - \mathcal{A}_{\rho \rho}(\de\fR, \ampR^i \fR^i) \right\} \mathcal{A}(\psi, \psi^i), 
\end{aligned}	
\end{equation}
for the HDG local problems and by
\begin{equation} \label{eq:linearGlobalPGDSpatial}
\begin{aligned}	
\mathcal{R}_{\hu}^m (\de \fHU \psi)	:=  & -  \sum_{l=1}^{\nNeu} \intBN{\de \fHU}{\bgN^l} \mathcal{A}(\psi, \lambda_N^l) \\
& {-}  \sum_{i=1}^{m} \left\{ \mathcal{A}_{\hu u}(\de \fHU, \ampU^i \fU^i)  +  \mathcal{A}_{\hu \hu}(\de \fHU, \ampHU^i \fHU^i) \right\} \mathcal{A}(\psi, \psi^i) \\
& {-}  \sum_{i=1}^{m} \sum_{k=1}^{\nadj}  \Bigl\{ \mathcal{A}^k_{\hu L}(\de \fHU, \ampL^i \fL^i) \mathcal{A}^k_{\vartheta}(\psi, \psi^i) \\
& {+} \left[ \mathcal{A}^k_{\hu p}(\de \fHU, \ampP^i \fP^i) + \mathcal{A}^k_{\hu \hu}(\de \fHU, \ampHU^i \fHU^i) \right] \mathcal{A}^k_{\vartheta}(\psi, \psi^i) \Bigr\} ,
\\
\mathcal{R}_{\rho}^m (\de \fR \psi)	 := & -\sum_{k=1}^{\nadj} \sum_{l=1}^{\nDir} \intBD{\de \fR}{\bgD^l \cdot \bA^k \bn} \mathcal{A}^k_{\vartheta}(\psi, \lambda_D^l)	\\
& {-}  \sum_{i=1}^{m} \sum_{k=1}^{\nadj} \mathcal{A}^k_{p\hu}(\de \fR, \ampHU^i \fHU^i) \mathcal{A}^k_{\vartheta}(\psi, \psi^i)
\end{aligned}	
\end{equation}
for the HDG global problems.

\section{Geometric mapping for the channel with two microswimmers}  \label{sc:mappingSwimmer}

\begin{figure}[!tb]
	\centering
	\subfigure[$\bm{\mathcal{M}}_{\mu_1}$]{\includegraphics[width=0.49\textwidth]{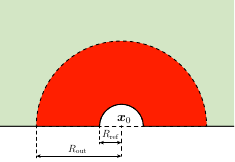}}
	\subfigure[$\bm{\mathcal{M}}_{\mu_2}$]{\includegraphics[width=0.49\textwidth]{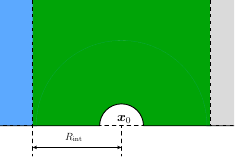}}
	\caption{Illustration of the piecewise nature of the mappings $\bm{\mathcal{M}}_{\mu_1}$ and $\bm{\mathcal{M}}_{\mu_2}$ detailed in equations~\eqref{eq:swimmerMapRadius} and~\eqref{eq:swimmerMapDistance} respectively in the vicinity of the sphere centred at $\bX_0$.}
	\label{fig:swimmeMapRadius}
\end{figure}

The mapping used in the example involving the flow around two microswimmers is designed as the composition of two mappings. The first mapping, $\bm{\mathcal{M}}_{\mu_1}$, is defined to account for the change of radius of the two spheres and it is written in the general separable expression of equation~\eqref{eq:displacementSep} with
\begin{equation} \label{eq:swimmerMapRadius}
\begin{aligned}
\bM^1_1(\bX) & = \left\{ 
\begin{split}
	\frac{1}{r} \bX_0^-  & \text{\quad if $\|\bX_0^-\| \leq \Rout$	}\\
	0 & \text{\quad otherwise}
\end{split}
\right.
\quad &
\psi^1_1(\mu_1) & = \dfrac{ \Rout (R^+(\mu_1) - \Rref)}{\Rout - \Rref}, 
\\
\bM^2_1(\bX) & = \left\{ 
\begin{split}
 \bX_0^- & \text{\quad if $\|\bX_0^-\| \leq \Rout$	}\\
 0 & \text{\quad otherwise}
\end{split}
\right.
\quad &
\psi^2_1(\mu_1) & = \dfrac{\Rout- R^+(\mu_1)}{\Rout - \Rref}, 
\\
\bM^3_1(\bX) & = \left\{ 
\begin{split}
\bX_0 & \text{\quad if $\|\bX_0^-\| \leq \Rout$	}\\
0 & \text{\quad otherwise}
\end{split}
\right.
\quad &
\psi^3_1(\mu_1) & = 1, 
\\
\bM^4_1(\bX) & = \left\{ 
\begin{split}
\frac{1}{r}  \bX_0^+ & \text{\quad if $\| \bX_0^+ \| \leq \Rout$	}\\
0 & \text{\quad otherwise}
\end{split}
\right.
\quad &
\psi^4_1(\mu_1) & = \dfrac{ \Rout (R^-(\mu_1) - \Rref)}{\Rout - \Rref}, 
\\
\bM^5_1(\bX) & = \left\{ 
\begin{split}
\bX_0^+  & \text{\quad if $\| \bX_0^+ \| \leq \Rout$	}\\
0 & \text{\quad otherwise}
\end{split}
\right.
\quad &
\psi^5_1(\mu_1) & = \dfrac{\Rout- R^-(\mu_1)}{\Rout - \Rref}, 
\\
\bM^6_1(\bX) & = \left\{ 
\begin{split}
-\bX_0 & \text{\quad if $\| \bX_0^+ \| \leq \Rout$	}\\
0 & \text{\quad otherwise}
\end{split}
\right.
\quad &
\psi^6_1(\mu_1) & = 1,
\end{aligned}
\end{equation}
where $\bX_0^\pm = \bX \pm \bX_0$, $\Rout=0.45$ and, as detailed in section~\ref{sc:swimmers}, $\bX_0 = (1.5,0)$ and $\Rref = 0.116$. The radius of the sphere centred at $\bX_0$ is defined as $R^+(\mu_1) = -0.0372\mu_1^2 + 0.0968\mu_1 + 0.25$ so that it takes value 0.116 for $\mu_1=-1$, 0.25 for $\mu_1=0$ and 0.3096 for $\mu_1=1$. The radius of the sphere centred at $-\bX_0$ is defined in terms of $R^+(\mu_1)$ in such a way that the total volume of the two spheres is maintained, namely $(R^+)^3 + (R^-)^3 = 1/32$. The piecewise nature of the mapping is illustrated in figure~\ref{fig:swimmeMapRadius}, in the vicinity of one of the spheres.

The second mapping, $\bm{\mathcal{M}}_{\mu_2}$, is defined to account for the change of distance between the spheres and it is written in the general separable expression of equation~\eqref{eq:displacementSep} with
\begin{equation} \label{eq:swimmerMapDistance}
\begin{aligned}
\bM^1_2(\bX) & = 
\begin{Bmatrix}
d(x) \\ 
0 
\end{Bmatrix} 
\qquad &
\psi^1_2(\mu_2) & = -x_0\mu_2/3, 
\\
\bM^2_2(\bX) & = \bX
\qquad &
\psi^2_2(\mu_2) & = 1,
\end{aligned}
\end{equation}
where the function $d(x)$ is given by
\begin{equation} \label{eq:fDistance}
d(x):= 
\begin{cases}
\displaystyle \frac{x+L}{x_0+\Rint-L} & \text{if }  x \in [-L,-x_0-\Rint]  \\
\displaystyle -1 & \text{if }  x \in [-x_0-\Rint,-x_0+\Rint]  \\
\displaystyle \frac{x}{x_0-\Rint} & \text{if }  x \in [-x_0+\Rint,x_0-\Rint]  \\
\displaystyle 1 & \text{if }  x \in [x_0-\Rint,x_0+\Rint]  \\
\displaystyle \frac{x-L}{x_0+\Rint-L} & \text{if }  x \in [x_0+\Rint,L]  
\end{cases},
\end{equation}
with $\Rint=0.47$ and, as detailed in section~\ref{sc:swimmers}, $L=6$. 

As illustrated in figure~\ref{fig:swimmeMapRadius} both mappings are defined in a piecewise form.
The mappings selected are only $\mathcal{C}^0$ on the artificial interfaces denoted by discontinuous lines in figure~\ref{fig:swimmeMapRadius}. Therefore, to facilitate the numerical integration of the terms involving the Jacobian and the adjoint of the mapping, the computational meshes selected are conforming with these interfaces, as it can be observed in the mesh displayed in figure~\ref{fig:swimmerMesh}. 

It is also worth noting that other mappings, with a smooth transition in the artificially created interfaces can be devised. Numerical experiments not reported here for brevity, demonstrate that the piecewise linear mapping described here results in a lower number of integration points required to ensure that errors due to the numerical integration are lower than the interpolation error. However, the choice of a smoother mapping circumvents the need to create meshes conforming with artificially created interfaces. In any case, as stressed in remark~\ref{rk:analyticalMapping}, this work focuses on the combination of the HDG and PGD formulations and for general geometries the general procedure described in~\cite{sevilla2020solution} is preferred, rather than the definition of analytical mappings.

\section{Geometric mapping for the corrugated channel}  \label{sc:mappingCorrugated}

Similarly to the previous example, the mapping used in the example involving the flow around a sphere in a corrugated channel is designed as the composition of two mappings. The first mapping, $\bm{\mathcal{M}}_{\mu_1}$, is defined to account for the change of radius of the sphere and it is written in the general separable expression of equation~\eqref{eq:displacementSep} with
\begin{equation} \label{eq:channelMapRadius}
\begin{aligned}
\bM^1_1(\bX) & = \left\{ 
\begin{split}
\frac{1}{r} \bX  & \text{\quad if $\|\bX\| \leq \Rout$	}\\
0 & \text{\quad otherwise}
\end{split}
\right.
\quad &
\psi^1_1(\mu_1) & = \dfrac{ \Rout (R(\mu_1) - \Rref)}{\Rout - \Rref}, 
\\
\bM^2_1(\bX) & = \left\{ 
\begin{split}
\bX & \text{\quad if $\|\bX\| \leq \Rout$	}\\
0 & \text{\quad otherwise}
\end{split}
\right.
\quad &
\psi^2_1(\mu_1) & = \dfrac{\Rout- R(\mu_1)}{\Rout - \Rref}, 
\\
\bM^3_1(\bX) & = \left\{ 
\begin{split}
\bX & \text{\quad if $\|\bX\| \leq \Rout$	}\\
0 & \text{\quad otherwise}
\end{split}
\right.
\quad &
\psi^3_1(\mu_1) & = 1,
\end{aligned}
\end{equation}
where $\Rout=0.4$ and $\Rref = 0.2$ and the radius of the sphere, centred at the origin, is defined as $R(\mu_1) = (\mu_1+2)/10$.

The second mapping, $\bm{\mathcal{M}}_{\mu_2}$, is defined to account for the change of amplitude in the undulatory part of the channel. It only affects the $y$ coordinate and, more precisely, only the definition of $f_n$ in equation~\eqref{eq:geoDefChannel}. More precisely, the profile of the channel is given byequation~\eqref{eq:geoDefChannel} with $f_n = 1/2 + \mu_2$.

\end{document}